\documentclass[11pt]{article}  
\title{Holomorphic Vector Bundles on the Quantum Projective Line}
\author{Paolo Aschieri, Masoud Khalkhali and Giovanni Landi}
\date{}

\usepackage[T1]{fontenc}
\usepackage[utf8]{inputenc}
\usepackage{lmodern}
\usepackage{amsmath,amssymb,amsthm,mathtools,mathrsfs}
\usepackage{microtype}
\usepackage[margin=1in]{geometry}
\usepackage{enumitem}
\usepackage[hidelinks]{hyperref}
\hypersetup{
  pdftitle={Holomorphic Vector Bundles on the Quantum Projective Line},
  pdfauthor={Paolo Aschieri, Masoud Khalkhali and Giovanni Landi}
}

\setlist{nosep}
\allowdisplaybreaks[2]
\numberwithin{equation}{section}

\newtheorem{theorem}{Theorem}[section]
\newtheorem{conjecture}[theorem]{Conjecture}
\newtheorem{proposition}[theorem]{Proposition}
\newtheorem{lemma}[theorem]{Lemma}
\newtheorem{corollary}[theorem]{Corollary}
\theoremstyle{definition}
\newtheorem{definition}[theorem]{Definition}
\theoremstyle{remark}
\newtheorem{remark}[theorem]{Remark}
\newtheorem{example}[theorem]{Example} 

\begin{document}
\maketitle

\begin{abstract}
We study Birkhoff--Grothendieck splitting for holomorphic vector bundles on
the quantum projective line in the Dolbeault-module framework. We prove that every polynomial
holomorphic bundle admitting a finite filtration by standard line objects
is a direct sum of standard lines, with uniquely determined splitting
integers. The proof uses a Dolbeault differential graded category and the
two-arrow Kronecker quiver. For a three-parameter rank-two extension
family, we determine the splitting types and distinguish the nonzero
splitting strata by the determinant of a quantum connecting map. For
general non-diagonal rank-two connections, we derive noncommutative
Riccati equations for invariant graph submodules and formulate the global
splitting problem in terms of holomorphic projectors. We also discuss
twisted positive Hochschild cocycles and formulate a variational principle
and a conformal-stability conjecture. On the algebraic side, we explain
why both the bare sheaf category of the quantum projective line and the
Artin--Zhang category of the quantum plane inherit classical
Birkhoff--Grothendieck splitting. We distinguish these results from
unrestricted standard-line splitting, which fails for polynomial gauges,
and outline a differential quantum GAGA programme with explicit attention
to regularity and equivariance.
\end{abstract}

\tableofcontents
\parskip = .75ex 

\section{Introduction}
\label{sec:introduction}

The Birkhoff--Grothendieck theorem classifies holomorphic vector bundles on
the Riemann sphere: every bundle $E$ of rank $r$ admits a decomposition
\begin{equation}
  E\simeq\mathcal O(k_1)\oplus\cdots\oplus\mathcal O(k_r),
  \qquad k_1\geq\cdots\geq k_r,
  \label{eq:intro-classical-BG}
\end{equation}
with uniquely determined integers $k_1,\ldots,k_r$
\cite{Birkhoff,Grothendieck,HazewinkelMartin}. By GAGA, its algebraic and
holomorphic formulations agree \cite{SerreGAGA}. The theorem can therefore be
expressed in terms of $\bar\partial$-operators, locally free sheaves, or
factorization of transition matrices over $\mathbb C[z,z^{-1}]$. In
noncommutative geometry these descriptions need not agree. A quantum
Birkhoff--Grothendieck problem must first specify the underlying geometry,
the category of bundles, and the allowed isomorphisms.

This paper studies the splitting problem from differential and
sheaf-theoretic viewpoints. On the differential side, we work with the
holomorphic structures of Khalkhali, Landi and van Suijlekom (KLvS)
on the standard Podle\'s sphere \cite{KLvS}. We prove splitting for polynomial
holomorphic bundles admitting a finite filtration by standard line objects,
and determine explicit rank-two splitting strata using the Kronecker quiver.
For general non-diagonal connections, we formulate the search for invariant
line summands as a noncommutative Riccati problem. On the sheaf side, we explain
why the underlying algebraic splitting theorem is already classical, and
identify the additional differential data that a comparison with KLvS geometry
must retain. We distinguish holomorphic splitting into possibly nonstandard
line objects from splitting into standard lines, and explain why retaining
arbitrary equivariance can obstruct the latter. The unrestricted holomorphic
splitting problem in the weaker sense, and the construction of a differential
comparison, remain open in the framework considered here.

The differential model starts from the fixed-point $*$-algebra
\[
  \mathcal A_q:=\mathcal A(\mathbb{CP}_q^1)
  =\mathcal A(\mathrm{SU}_q(2))_0,
\]
with its covariant Dolbeault calculus. The $\mathbb Z$-grading of
$\mathcal A(\mathrm{SU}_q(2))$ gives canonical line modules
\begin{equation}
  \mathcal L_n:=\mathcal A(\mathrm{SU}_q(2))_n,
  \qquad n\in\mathbb Z,
  \label{eq:intro-line-modules}
\end{equation}
which are finitely generated projective rank-one modules over $\mathcal A_q$.
Their natural bimodule structures give
\[
  \mathcal L_n^*\simeq\mathcal L_{-n},
  \qquad
  \mathcal L_n\otimes_{\mathcal A_q}\mathcal L_m\simeq\mathcal L_{n+m}.
\]
The canonical bimodule connection on $\mathcal L_1$
induces a distinguished
$\bar\partial$-connection $\overline\nabla_{(n)}$ on each $\mathcal L_n$.
With the convention
\begin{equation}
  \mathcal O_q(k):=(\mathcal L_{-k},\overline\nabla_{(-k)}),
  \label{eq:intro-Oq}
\end{equation}
one has
\begin{equation}
  \dim H^0\bigl(\mathcal O_q(k)\bigr)
  =\begin{cases}k+1,&k\geq0,\\0,&k<0.\end{cases}
  \label{eq:intro-H0}
\end{equation}
The dimensions are classical, whereas multiplication of sections recovers
the quantum plane.

A KLvS holomorphic bundle is a finitely generated projective left
$\mathcal A_q$-module $\mathcal E$ with a flat connection
\[
  \overline\nabla:\mathcal E\longrightarrow
  \Omega^{0,1}(\mathcal A_q)\otimes_{\mathcal A_q}\mathcal E.
\]
Since $\Omega^{0,2}(\mathcal A_q)=0$, every $\bar\partial$-connection is flat.
The classification problem thus concerns gauge equivalence rather than
integrability. Three questions must be distinguished: whether the underlying
module splits into canonical line modules; whether a decomposition can be
chosen to preserve the holomorphic structure; and whether every holomorphic
structure on a canonical line module is gauge equivalent to its distinguished KLvS
structure. Holomorphic splitting and line-bundle rigidity are separate
assertions. A decomposition into the objects $\mathcal O_q(k)$ would have
the form
\begin{equation}
  (\mathcal E,\overline\nabla)
  \simeq\mathcal O_q(k_1)\oplus\cdots\oplus\mathcal O_q(k_r),
  \qquad k_1\geq\cdots\geq k_r,
  \label{eq:intro-strong-BG}
\end{equation}
with uniquely determined splitting integers. This form cannot
hold for every polynomial holomorphic bundle under polynomial gauge
transformations: even the trivial line module has nonstandard structures,
as shown in Proposition~\ref{prop:polynomial-line-nonrigidity}. The preprint
\cite{GravemanEtAl} also constructs nonstandard line structures for an
enlarged completed gauge group, while \cite{BiswasGuinKumar} studies the
gauge equation through fixed-point methods. The regularity class and the
gauge group must therefore be part of every general classification statement.

Our splitting results concern a concrete subcategory: polynomial holomorphic
bundles with a finite filtration whose quotients are standard KLvS line
objects. We show that every such bundle is holomorphically isomorphic to a
direct sum of these line objects, with splitting integers uniquely determined
up to permutation
(Theorem~\ref{thm:line-filtered-splitting}). The proof uses the Dolbeault
differential graded category and the Beilinson object
$\mathcal T=\mathcal L_0\oplus\mathcal L_1$, whose endomorphism algebra yields
the two-arrow Kronecker quiver. A sufficiently negative twist reduces the
classification to finite-dimensional Kronecker representations; positivity
of the generic rank of projective summands selects the representations
corresponding to standard lines. This proves a splitting theorem for
line-filtered bundles, without assuming that every holomorphic bundle admits
such a filtration.

In particular, all upper-triangular rank-two extensions of standard lines
have split middle terms. This does not imply that their defining short exact
sequences split: the quantum Euler sequence is a basic example. For the
three-parameter extension family studied in
Section~\ref{sec:rank-two-kronecker}, the zero extension has middle term
$\mathcal L_2\oplus\mathcal L_{-2}$, while each nonzero extension has middle
term either $\mathcal L_0^{\oplus2}$ or
$\mathcal L_1\oplus\mathcal L_{-1}$. The determinant of the quantum connecting
map separates these two nonzero splitting strata. We also treat higher-rank
upper-triangular examples.

A general rank-two connection need not preserve a prescribed line filtration.
For such connections, we derive noncommutative Riccati equations expressing
invariance of graph submodules. Complementary invariant graphs yield a
holomorphic direct-sum decomposition. We give an explicit solvable
non-diagonal family and discuss a local analytic approach under stated
completion, boundedness, and regularity hypotheses. The global problem is
the existence of a holomorphic rank-one projector; identifying its image
with a canonical line module and reducing the induced line connection to
standard form are further classification questions.

Alongside bundle splitting, we examine the positive twisted Hochschild
cocycle associated with the KLvS complex structure. The modular Haar state
requires a twisted formulation of positivity and cohomology. Motivated by
Connes's positive-cocycle description of conformal geometry
\cite{ConnesNCG} and by conformal perturbations of the noncommutative torus
\cite{ConnesTretkoff,FathizadehKhalkhaliGB,ConnesMoscoviciModular}, we formulate
a twisted variational principle and a conformal-stability conjecture. These
questions concern the canonical complex structure and are logically separate
from the splitting theorems.

The algebraic viewpoint uses the sheaf approach to quantum principal bundles
of Aschieri, Fioresi and Latini and its differential extension with Weber
\cite{AschieriFioresiLatini,AschieriFioresiLatiniWeber}. Local coinvariants of a
principal sheaf $\mathcal F_q$ recover the base structure sheaf, and a
finite-dimensional left comodule $V$ of the structure Hopf algebra $H$ gives
an associated sheaf
\[
  \mathcal E_V(U)=\mathcal F_q(U)\square_H V.
\]
On cleft charts this is finite locally free. Retaining the comodule or
equivariant data defines a finer classification problem than classifying
the underlying locally free sheaf.

For the quantum projective line, the bare base sheaf on the standard
two-chart topology has local algebras
\begin{equation}
  \begin{gathered}
    \mathcal O_{q,\mathbb P^1}(U_0)=\Bbbk[z],\qquad
    \mathcal O_{q,\mathbb P^1}(U_1)=\Bbbk[w],\\
    \mathcal O_{q,\mathbb P^1}(U_0\cap U_1)=\Bbbk[z,z^{-1}],
    \qquad w=z^{-1}.
  \end{gathered}
  \label{eq:intro-local-sheaf}
\end{equation}
Here $\Bbbk$ is a field; the comparison with KLvS geometry concerns
$\Bbbk=\mathbb C$. These chart algebras and restriction maps are classical.
Consequently, Hazewinkel--Martin factorization \cite{HazewinkelMartin} gives
\begin{equation}
  \mathcal E\simeq\mathscr O_q(n_1)\oplus\cdots\oplus\mathscr O_q(n_r)
  \label{eq:intro-bare-splitting}
\end{equation}
for every bare locally free sheaf, uniquely up to permutation, where
$\mathscr O_q(n)\simeq\mathcal O(n)$ as underlying line sheaves. The
noncommutative structure remains visible in the total principal sheaf,
the calculus, and the ordered multiplication of sections:
\begin{equation}
  \bigoplus_{n\geq0}H^0\bigl(\mathbb P^1,\mathscr O_q(n)\bigr)
  \simeq\Bbbk\langle x_0,x_1\rangle/(x_0x_1-qx_1x_0),
  \label{eq:intro-quantum-section-ring}
\end{equation}
where the left-hand side carries the multiplication inherited from the
quantum total algebra.

Likewise, the quantum plane is a Zhang twist of the commutative polynomial
ring, and its Artin--Zhang category is equivalent to
$\operatorname{coh}(\mathbb P^1)$ \cite{ArtinVdB,ArtinZhang,ZhangTwist}.
Its vector-bundle objects therefore satisfy Birkhoff--Grothendieck; the
shift and section multiplication retain the deformation. Strict
$\operatorname{GL}_q(2)$-valued transition maps over the commutative Laurent
overlap impose additional restrictions on gluing and gauge data. Although
they admit a diagonal representative of every bare splitting type, they
do not encode arbitrary frame changes or extension data.

The remaining comparison problem motivates a differential quantum GAGA
programme. It requires the construction of source and target categories
with fixed regularity, gauge groups, and compatible differential structures,
and a comparison sending
$\mathcal O_q(n)$ to $\mathscr O_q(n)$. Exactness, compatibility with
cohomology, and tensor compatibility on suitable bimodule subcategories
are requirements to establish. Equivariant compatibility of a comparison
functor must be distinguished from equivariant splitting: the latter is
not expected for arbitrary retained representations. The classical
bare-sheaf theorem alone cannot provide this comparison or settle
unrestricted holomorphic splitting into possibly nonstandard line objects.

Section~\ref{sec:klv-module-picture} develops the differential framework,
the line objects $\mathcal L_n$, the twisted-positive cocycle questions, and the
regularity-dependent splitting and rigidity problems.
Section~\ref{sec:rank-two-test-objects} constructs
rank-two and higher-rank test objects.
Section~\ref{sec:rank-two-kronecker} determines the rank-two extension family
through the Kronecker model, and
Section~\ref{sec:all-line-filtered-examples-split} proves the line-filtered
splitting theorem. Section~\ref{sec:nc-riccati-splitting} treats general
non-diagonal connections through invariant graphs and projectors.
Finally, Section~\ref{sec:sheaf-approach} develops the sheaf and
noncommutative-projective viewpoints and sets out the differential comparison
programme and the obstruction to unrestricted equivariant splitting.

In the Dolbeault-module setting, the terms ``standard line module'',
``standard line bundle'' and ``standard line object'' refer to the
canonical modules $\mathcal L_n$ equipped with the holomorphic structures
$\overline\nabla_{(n)}$ introduced in
\cite{KLvS} and recalled in \eqref{eq:standard-KLvS-delbar}
in Section~\ref{sec:klv-module-picture}.
When regarded as a holomorphic object, $\mathcal L_n$ will denote this
module with its prescribed holomorphic structure unless another
structure is explicitly specified.

\section{Holomorphic structures and the KLvS quantum projective line}
\label{sec:klv-module-picture}

This section reviews the module-theoretic framework of  \cite{KLvS}. We recall complex differential calculi and
Dolbeault connections, using the noncommutative torus as a comparison example.
We then describe the standard Podle\'s sphere and formulate the splitting and
rigidity questions used in the subsequent sections.

There are four levels of structure that should not be conflated:
\begin{enumerate}[label=\textup{(\roman*)},leftmargin=2.3em]
\item a differential $*$-calculus and a Dolbeault bigrading on an algebra;
\item a flat $\bar\partial$-connection on a module;
\item a choice of regularity, such as polynomial, smooth Fr\'echet, or
      $C^*$-algebraic;
\item a sheaf or gluing description on a projective quantum base.
\end{enumerate}

\noindent
KLvS work primarily with the first three levels. The sheaf-theoretic
approach of 
\cite{AschieriFioresiLatini,AschieriFioresiLatiniWeber} treats the fourth.
Here \emph{algebraic} means over the polynomial coordinate algebra with its
algebraic differential calculus, rather than a locally free sheaf in a
noncommutative $\operatorname{Proj}$ category. A comparison of these
frameworks requires additional structure.

\subsection{Complex structures on involutive algebras}
\label{subsec:complex-star-algebras}

Let $\mathcal A$ be a unital complex $*$-algebra.  A differential $*$-calculus on
$\mathcal A$ is a graded differential $*$-algebra
\[
  \Omega^\bullet(\mathcal A)
  =\bigoplus_{k\geq 0}\Omega^k(\mathcal A),
  \qquad \Omega^0(\mathcal A)=\mathcal A,
\]
with a degree-one derivation
\[
  \mathrm d:\Omega^k(\mathcal A)\longrightarrow\Omega^{k+1}(\mathcal A),
  \qquad \mathrm d^2=0,
\]
satisfying
\[
  \mathrm d(\alpha\beta)
  =(\mathrm d\alpha)\beta+(-1)^{|\alpha|}\alpha\,\mathrm d\beta,
  \qquad
  \mathrm d(\alpha^*)=(\mathrm d\alpha)^* .
\]
No graded commutativity is assumed.  Thus both the product of forms and their
bimodule commutation rules are part of the geometry.

\begin{definition}[Complex differential calculus]
\label{def:complex-calculus}
A complex structure on the differential $*$-calculus
$\Omega^\bullet(\mathcal A)$ is a multiplicative bigrading
\[
  \Omega^n(\mathcal A)
  =\bigoplus_{p+q=n}\Omega^{p,q}(\mathcal A),
  \qquad
  \Omega^{p,q}(\mathcal A) \Omega^{r,s}(\mathcal A) \subseteq\Omega^{p+r,q+s}(\mathcal A),
  \qquad
  \bigl(\Omega^{p,q}(\mathcal A)\bigr)^*
  =\Omega^{q,p}(\mathcal A),
\]
together with operators
\[
  \partial:\Omega^{p,q}(\mathcal A)\longrightarrow\Omega^{p+1,q}(\mathcal A),
  \qquad
  \bar\partial:\Omega^{p,q}(\mathcal A)\longrightarrow\Omega^{p,q+1}(\mathcal A),
\]
such that
\[
  \mathrm d=\partial+\bar\partial,
  \qquad
  \partial^2=\bar\partial^2=0,
  \qquad
  \partial\bar\partial+\bar\partial\partial=0.
\]
\end{definition}

KLvS formulated the minimal bigraded version of this definition and observed that
additional integrability hypotheses may be needed in general.  The formulation
above incorporates integrability from the outset, as is now customary in
calculus-based treatments of noncommutative complex geometry; see
\cite{BeggsSmith,OBuachalla}.  The algebra of holomorphic elements is
\begin{equation}
  \mathcal O(\mathcal A)
  :=\ker\!\left(
       \bar\partial:\mathcal A\longrightarrow\Omega^{0,1}(\mathcal A)
     \right).
  \label{eq:holomorphic-elements}
\end{equation}
The Leibniz rule makes $\mathcal O(\mathcal A)$ a subalgebra of $\mathcal A$.  On a
compact connected classical complex manifold it consists of the global
holomorphic functions, and hence can be very small.  For the Riemann sphere and
for its KLvS deformation it is just $\mathbb C$.

\begin{example}[Classical and Lie-algebraic examples]
For a complex manifold $X$, the complexified de Rham calculus has the familiar
decomposition
\[
  \Omega^n(X;\mathbb C)
  =\bigoplus_{p+q=n}\Omega^{p,q}(X),
  \qquad
  \mathrm d=\partial+\bar\partial.
\]
This is the motivating example, but the definition does not require the algebra to
be commutative.

A useful algebraic construction
starts with a real Lie algebra $\mathfrak g$ acting on
the algebra
$\mathcal A$ by $*$-derivations.  Suppose that for its complexification,
\[
  \mathfrak g_{\mathbb C}
  =\mathfrak g^{1,0}\oplus\mathfrak g^{0,1},
  \qquad
  \overline{\mathfrak g^{1,0}}=\mathfrak g^{0,1},
\]
and that the two summands are Lie subalgebras.  The Chevalley--Eilenberg calculus
\[
  \Omega^k_{\mathfrak g}(\mathcal A)
  =\operatorname{Hom}_{\mathbb C}
      \bigl(\Lambda^k\mathfrak g_{\mathbb C},\mathcal A\bigr)
\]
then inherits the bigrading
\[
  \Omega^{p,q}_{\mathfrak g}(\mathcal A)
  =
  \operatorname{Hom}_{\mathbb C}
  \bigl(\Lambda^p\mathfrak g^{1,0}
        \otimes\Lambda^q\mathfrak g^{0,1},\mathcal A\bigr).
\]
This construction includes translation-invariant complex structures on
noncommutative tori and many calculi on quantum homogeneous spaces.
\end{example}

\begin{example}[The smooth noncommutative two-torus]
Let $A_\theta^\infty$ be the smooth noncommutative two-torus, the Fr\'echet
$*$-algebra of rapidly decreasing series in two generators $U_1$ and $U_2$,
\[
  a=\sum_{m,n\in\mathbb Z}a_{m,n}U_1^mU_2^n,
  \qquad
  U_1U_2=e^{2\pi i\theta}U_2U_1.
\]
The canonical derivations are normalized by
\[
  \delta_j(U_k)=2\pi i\,\delta_{jk}U_k,
  \qquad j,k=1,2.
\]
For $\tau$ in the upper half-plane
$\mathbb H=\{\tau\in\mathbb C:\operatorname{Im}\tau>0\}$,
\begin{equation}
  \partial_\tau
  =\frac{-\bar\tau\,\delta_1+\delta_2}{\tau-\bar\tau},
  \qquad
  \bar\partial_\tau
  =\frac{\tau\,\delta_1-\delta_2}{\tau-\bar\tau}.
  \label{eq:torus-dolbeault}
\end{equation}
This is the Dolbeault decomposition associated with the conformal class of the
constant metric determined by $\tau$; see
\cite{ConnesNCG,DiengSchwarz,KLvS}.  The only holomorphic elements are constants.

Let $\operatorname{tr}_\theta$ denote the canonical trace on $A_\theta^\infty$.  The fundamental cyclic
two-cocycle is
\begin{equation}
  \Psi(a_0,a_1,a_2)
  =\frac{i}{2\pi}\operatorname{tr}_\theta
     \left(
       a_0\bigl(\delta_1(a_1)\delta_2(a_2)
                -\delta_2(a_1)\delta_1(a_2)\bigr)
     \right),
  \label{eq:torus-cyclic-cocycle}
\end{equation}
whereas the positive Hochschild cocycle associated with the complex structure is
\begin{equation}
  \Phi_\tau(a_0,a_1,a_2)
  =\frac{2}{\pi}\operatorname{tr}_\theta
     \bigl(a_0\,\partial_\tau(a_1)\bar\partial_\tau(a_2)\bigr).
  \label{eq:torus-positive-cocycle}
\end{equation}
To compare different moduli in the fixed fundamental cyclic class, a
normalization is needed. Set $y=\operatorname{Im}\tau$ and
\begin{equation}
  \widehat\Phi_\tau:=y\, \Phi_\tau.
  \label{eq:normalized-torus-positive-cocycle}
\end{equation}
Indeed,
\[
 \widehat\Phi_\tau-\Psi=b\chi_\tau,
 \qquad
 \chi_\tau(a,b)=-\frac{y}{\pi}\operatorname{tr}_\theta
                    \bigl(a\,\partial_\tau\bar\partial_\tau b\bigr),
 \qquad B\chi_\tau=0.
\]
The last equality follows from invariance of the trace under the canonical
derivations. Thus $\widehat\Phi_\tau$ belongs to the fixed cyclic affine
class of $\Psi$, whereas $[\Phi_\tau]=y^{-1}[\Psi]$. We retain the
unnormalized $\Phi_\tau$ when describing the Dolbeault norm below; the
positive scalar factor does not change the encoded complex structure.
\end{example}

\noindent
The relation between \eqref{eq:torus-cyclic-cocycle} and
\eqref{eq:torus-positive-cocycle} is the model for the twisted construction on the
Podle\'s sphere.

\subsection{Holomorphic structures on projective modules}
\label{subsec:holomorphic-projective-modules}

We first use left modules, in accordance with KLvS.  The corresponding right-module
definitions are obtained by reversing the order of the Leibniz rule.

\begin{definition}[$\bar\partial$-connection]
\label{def:dbar-connection}
Let $\mathcal E$ be a left $\mathcal A$-module.  A $\bar\partial$-connection on
$\mathcal E$ is a complex-linear map
\[
  \overline\nabla:\mathcal E
  \longrightarrow
  \Omega^{0,1}(\mathcal A)\otimes_{\mathcal A}\mathcal E
\]
satisfying
\begin{equation}
  \overline\nabla(a\xi)
  =a\,\overline\nabla(\xi)+\bar\partial a\otimes\xi,
  \qquad a\in\mathcal A,\quad \xi\in\mathcal E.
  \label{eq:left-dbar-leibniz}
\end{equation}
It extends to
$\Omega^{0,k}(\mathcal A)\otimes_{\mathcal A}\mathcal E$ by
\begin{equation}
  \overline\nabla(\omega\otimes\xi)
  =\bar\partial\omega\otimes\xi
   +(-1)^k\omega\wedge\overline\nabla\xi.
  \label{eq:dbar-extension}
\end{equation}
The $(0,2)$-curvature is $F^{0,2}_{\overline\nabla}:=\overline\nabla^2$.
A holomorphic structure is a flat $\bar\partial$-connection,
$F^{0,2}_{\overline\nabla}=0$.  If $\mathcal E$ is finitely generated and projective,
the pair $(\mathcal E,\overline\nabla)$ is called a holomorphic vector bundle.
\end{definition}

This is the module-theoretic definition used in \cite{KLvS}; broader formulations
and refinements may be found in \cite{BeggsSmith}.  If $\mathcal E$ is finitely
generated projective, it always admits a $\bar\partial$-connection.  Indeed, choose
left $\mathcal A$-linear maps
\[
  \mathcal E\xrightarrow{\ i\ }\mathcal A^N
  \xrightarrow{\ \pi\ }\mathcal E,
  \qquad \pi \circ i=\operatorname{id}_{\mathcal E},
\]
and put
\begin{equation}
  \overline\nabla_0
  =(1\otimes\pi)\circ\bar\partial^{(N)}\circ i.
  \label{eq:grassmann-dbar}
\end{equation}
Any other $\bar\partial$-connection  is of the form
\begin{equation}
  \overline\nabla=\overline\nabla_0+\phi ,
  \qquad
  \phi \in
  \operatorname{Hom}_{\mathcal A}
  \bigl(\mathcal E,\Omega^{0,1}(\mathcal A)
                         \otimes_{\mathcal A}\mathcal E\bigr).
  \label{eq:affine-space-connections}
\end{equation}
The space of connections is therefore affine, with flatness imposing a
nonlinear constraint in general.

Two holomorphic structures on the same module are gauge equivalent when
\begin{equation}
  \overline\nabla_2=(1\otimes g^{-1})\circ\overline\nabla_1\circ g
  \label{eq:gauge-equivalence}
\end{equation}
for some $g\in\operatorname{Aut}_{\mathcal A}(\mathcal E)$.  More generally, a
holomorphic map from $(\mathcal E,\overline\nabla_{\mathcal E})$ to
$(\mathcal F,\overline\nabla_{\mathcal F})$ is an $\mathcal A$-linear map $T$ satisfying
\[
  (1\otimes T)\overline\nabla_{\mathcal E}
  =\overline\nabla_{\mathcal F}T.
\]

\begin{example}[The free-module Maurer--Cartan equation on the noncommutative torus]
On the smooth noncommutative torus $A_\theta^\infty$, for the matrix formulation temporarily use the right free module
$E=(A_\theta^\infty)^r$, written as column vectors.  A right
$\bar\partial$-connection has the form
\[
  \overline\nabla(e)=\bar\partial e + \omega \, e,
  \qquad
\omega \in M_r\bigl(\Omega^{0,1}(A_\theta^\infty)\bigr),
\]
and its curvature is
\begin{equation}
  F^{0,2}_{\overline\nabla}
  =\bar\partial \omega + \omega \wedge \omega.
  \label{eq:maurer-cartan}
\end{equation}
A right-module automorphism is given by $e\mapsto ge$, with
$g\in\operatorname{GL}_r(A_\theta^\infty)$.  Under the gauge action
$\overline\nabla\mapsto g^{-1}\overline\nabla \, g$,
\begin{equation}
  \omega\longmapsto
  \omega^g=g^{-1} \omega g+g^{-1}\bar\partial g.
  \label{eq:gauge-potential}
\end{equation}
\end{example}

In a differential calculus of complex dimension one,
$\Omega^{0,2}(\mathcal A)=0$.
Every $\bar\partial$-connection is then flat, so the classification problem
concerns gauge orbits rather than integrability. The regularity class remains
part of the data: algebraic, smooth and $C^*$-invertible gauge transformations
need not give the same quotient.

A flat structure gives the Dolbeault complex
\begin{equation}
  0\longrightarrow\mathcal E
  \xrightarrow{\overline\nabla}
  \Omega^{0,1}\otimes_{\mathcal A}\mathcal E
  \xrightarrow{\overline\nabla}
  \Omega^{0,2}\otimes_{\mathcal A}\mathcal E
  \longrightarrow\cdots .
  \label{eq:module-dolbeault-complex}
\end{equation}
Its cohomology is denoted $H^j(\mathcal E,\overline\nabla)$, and
\begin{equation}
  H^0(\mathcal E,\overline\nabla)=\ker\overline\nabla
  \label{eq:holomorphic-sections}
\end{equation}
is the space of holomorphic sections.  It is naturally a left
$\mathcal O(\mathcal A)$-module.

\paragraph{Bimodule connections and tensor products.}

The canonical line modules considered below are bimodules.  Tensoring their connections
requires a replacement for the ordinary flip.  A bimodule connection on an
$\mathcal A$-bimodule $\mathcal E$ consists of a left connection $\nabla$ and an
$\mathcal A$-bimodule isomorphism
\[
  \sigma_{\mathcal E}:
  \mathcal E\otimes_{\mathcal A}\Omega^1(\mathcal A)
  \longrightarrow
  \Omega^1(\mathcal A)\otimes_{\mathcal A}\mathcal E
\]
such that
\begin{equation}
  \nabla(\xi a)
  =\nabla(\xi)a+
   \sigma_{\mathcal E}(\xi\otimes\mathrm d a).
  \label{eq:twisted-right-leibniz}
\end{equation}
For bimodules $(\mathcal E_j,\nabla_j,\sigma_j)$, the tensor-product connection is
\begin{equation}
  \nabla_{1\otimes 2}
  =\nabla_1\otimes 1
   +(\sigma_1\otimes 1)(1\otimes\nabla_2).
  \label{eq:tensor-product-connection}
\end{equation}
See, for example, \cite{BeggsSmith,CarotenutoOBuachalla} for the general formalism.  In a
noncommutative calculus, the curvature of \eqref{eq:tensor-product-connection} is
not automatically the sum of the curvatures: mixed terms can occur.  For the corresponding Dolbeault connections in complex
dimension one, this obstruction to flatness disappears because
$\Omega^{0,2}(\mathcal A)=0$.  KLvS construct the required twisted flips explicitly for their
canonical line modules.

\subsection{Noncommutative elliptic curves and Heisenberg modules}
\label{subsec:heisenberg-modules}

The noncommutative torus supplies explicit nonfree holomorphic modules and
a genus-one comparison. Rieffel constructed its projective modules, and the
holomorphic theory was developed by Polishchuk and Schwarz and by Dieng and
Schwarz \cite{RieffelModules,PolishchukSchwarz,DiengSchwarz}.
Unlike the genus-zero setting, noncommutative elliptic curves have
indecomposable higher-rank holomorphic bundles.

For this subsection it is convenient to use right modules and the normalization of
\cite{PolishchukSchwarz}.  Set
\begin{equation}
  \delta_\tau\left(
     \sum_{m,n}a_{m,n}U_1^mU_2^n
  \right)
  =2\pi i\sum_{m,n}(m\tau+n)a_{m,n}U_1^mU_2^n,
  \qquad \operatorname{Im}\tau<0.
  \label{eq:delta-tau-PS}
\end{equation}
This convention is related to \eqref{eq:torus-dolbeault} by replacing
the upper-half-plane parameter there by $-\tau$ and rescaling the
Dolbeault operator by a nonzero scalar.  A holomorphic structure on a right
$A_\theta^\infty$-module $E$ is an operator satisfying
\begin{equation}
  \overline\nabla(ea)=\overline\nabla(e)a+e\,\delta_\tau(a).
  \label{eq:right-torus-holomorphic}
\end{equation}
Let $m,n\in\mathbb Z$, with $m\neq 0$ and $n+m\theta>0$.  The Heisenberg module
\begin{equation}
  E_{n,m}(\theta)
  =\mathcal S(\mathbb R\times\mathbb Z/m\mathbb Z)
  \label{eq:heisenberg-module}
\end{equation}
has right action
\begin{align}
  (fU_1)(x,\alpha)
  &=f\left(x-\frac{n+m\theta}{m},\alpha-1\right),
  \label{eq:heisenberg-U1}\\
  (fU_2)(x,\alpha)
  &=e^{2\pi i(x-\alpha n/m)}f(x,\alpha).
  \label{eq:heisenberg-U2}
\end{align}
Its numerical invariants are
\begin{equation}
  \deg E_{n,m}=m,
  \qquad
  \operatorname{rk}E_{n,m}=n+m\theta,
  \qquad
  \mu(E_{n,m})=\frac{m}{n+m\theta} \, ,
  \label{eq:heisenberg-invariants}
\end{equation}
where $\mu$ denotes the slope of the module \cite{PolishchukClassification}.
When $\gcd(n,m)=1$, the module is basic.  The standard one-parameter family of
holomorphic structures is
\begin{equation}
  \overline\nabla_z f
  =\frac{\partial f}{\partial x}
   +2\pi i\bigl(\tau\mu(E_{n,m})x+z\bigr)f,
  \qquad z\in\mathbb C.
  \label{eq:heisenberg-delbar}
\end{equation}
For a basic module, two members of this family are gauge equivalent precisely when
\begin{equation}
  z'-z\in
  \frac{1}{\operatorname{rk}E_{n,m}}
  \bigl(\mathbb Z+\tau\mathbb Z\bigr).
  \label{eq:heisenberg-gauge}
\end{equation}
Thus even a fixed basic projective module carries an elliptic family of standard
holomorphic structures.
The cohomology of \eqref{eq:heisenberg-delbar} is explicit.  Under our
standing assumption $\operatorname{rk}E_{n,m}>0$, if $m>0$ (equivalently,
$\mu(E_{n,m})>0$), then
\begin{equation}
  H^1(E_{n,m},\overline\nabla_z)=0,
  \qquad
  \dim H^0(E_{n,m},\overline\nabla_z)=m,
  \label{eq:heisenberg-H0}
\end{equation}
and a basis is formed by the theta vectors
\begin{equation}
  \phi_\alpha^z(x,\beta)
  =
  \exp\!\left(
    2\pi i\left[-\frac{\tau\mu(E_{n,m})}{2}x^2-zx\right]
  \right)\delta_\alpha(\beta),
  \qquad
  \alpha\in\mathbb Z/m\mathbb Z.
  \label{eq:theta-vectors}
\end{equation}
For negative degree (and positive rank), $H^0$ vanishes and $H^1$ has dimension $-\deg(E)$.
Polishchuk later proved that every holomorphic vector bundle on a
noncommutative two-torus is obtained by successive extensions of standard
bundles \cite{PolishchukClassification}.  Equivalently, the holomorphic category
is the heart of a tilted $t$-structure in the derived category of an elliptic
curve. Extensions and indecomposable objects are therefore intrinsic to this
genus-one theory, in contrast with the Birkhoff--Grothendieck decomposition
on the projective line $\mathbb{CP}^1$.

\subsection{The KLvS complex structure on the standard Podle\'s sphere}
\label{subsec:podles-klv}

Let $0<q<1$.  The coordinate Hopf $*$-algebra
$\mathcal A(\mathrm{SU}_q(2))$ is generated by $a,c$ with relations
\begin{equation}
  ac=qca,\qquad
  ac^*=qc^*a,\qquad
  cc^*=c^*c,
  \label{eq:SUq-relations-one}
\end{equation}
and
\begin{equation}
  a^*a+c^*c=1,
  \qquad
  aa^*+q^2cc^*=1.
  \label{eq:SUq-relations-two}
\end{equation}
Equivalently, the matrix
\[
  \begin{pmatrix}a&-qc^*\\ c&a^*\end{pmatrix}
\]
is unitary \cite{WoronowiczSUq2}.  The right $\mathrm U(1)$-action of the quantum
Hopf fibration gives a $\mathbb Z$-grading.  We use the   convention
\begin{equation}
  \deg a=\deg c=-1,
  \qquad
  \deg a^*=\deg c^*=1.
  \label{eq:grading-convention}
\end{equation}
The fixed-point algebra
\begin{equation}
  \mathcal A_q
  :=\mathcal A(\mathbb{CP}_q^1)
  =\mathcal A(\mathrm{SU}_q(2))_0
  \label{eq:quantum-hopf-fibration}
\end{equation}
is the coordinate algebra of the standard Podle\'s sphere
\cite{PodlesSpheres}.  It is generated by
\begin{equation}
  B_-=ac^*,\qquad
  B_+=ca^*,\qquad
  B_0=cc^*,
  \label{eq:podles-generators}
\end{equation}
with $B_0^*=B_0$, $B_+^*=B_-$ and relations
\begin{align}
  B_-B_0&=q^2B_0B_-,
  &
  B_+B_0&=q^{-2}B_0B_+,
  \label{eq:podles-relations-one}\\
  B_-B_+&=q^2B_0(1-q^2B_0),
  &
  B_+B_-&=B_0(1-B_0).
  \label{eq:podles-relations-two}
\end{align}

\paragraph{Canonical line modules.}

For $n\in\mathbb Z$ let
\begin{equation}
  \mathcal L_n
  :=\mathcal A(\mathrm{SU}_q(2))_n
  =\{x \in \mathcal A(\mathrm{SU}_q(2)) : \deg x=n\}.
  \label{eq:line-modules}
\end{equation}
Then $\mathcal L_0=\mathcal A_q$, each $\mathcal L_n$ is an
$\mathcal A_q$-bimodule and a finitely generated projective left module of rank
one, and
\begin{equation}
  \mathcal L_n^*\simeq\mathcal L_{-n},
  \qquad
  \mathcal L_n\otimes_{\mathcal A_q}\mathcal L_m
  \simeq\mathcal L_{n+m}.
  \label{eq:line-tensor-law}
\end{equation}
The multiplication map induces the second isomorphism.  Explicit module generating sets
are
\begin{equation}
  \{c^{*\mu}a^{*\,n-\mu}:0\leq\mu\leq n\},
  \qquad n\geq 0,
  \label{eq:positive-line-generators}
\end{equation}
and
\begin{equation}
  \{a^\mu c^{|n|-\mu}:0\leq\mu\leq |n|\},
  \qquad n\leq 0.
  \label{eq:negative-line-generators}
\end{equation}
For $n\neq0$, these are nonfree projective modules; $\mathcal L_0=\mathcal A_q$
is free. Their winding number is $-n$ in the KLvS convention; see \cite{BrzezinskiMajidMonopole,KLvS}.

\paragraph{The covariant Dolbeault calculus.}

Let $\mathcal U_q(\mathfrak{su}(2))$ be the Hopf $*$-algebra paired with
$\mathcal A(\mathrm{SU}_q(2))$. Its generators $K,K^{-1},E,F$ satisfy
\[
 KK^{-1}=K^{-1}K=1,\qquad
 K^{\pm1}E=q^{\pm1}EK^{\pm1},\qquad
 K^{\pm1}F=q^{\mp1}FK^{\pm1},
\]
\[
 [E,F]=\frac{K^2-K^{-2}}{q-q^{-1}}.
\]
The $*$-structure is simply $K^*=K$, $E^*=F$, $F^*=E$.
The bilinear pairing $\langle \cdot , \cdot \rangle $ between $\mathcal U_q(\mathfrak{su}(2))$ and $\mathcal A(\mathrm{SU}_q(2))$
gives canonical left and right $\mathcal U_q(\mathfrak{su}(2))$-module algebra
structures on $\mathcal A(\mathrm{SU}_q(2))$,
\[
h \rhd x := x_{(1)} \langle h , x_{(2)} \rangle , \qquad  x \lhd h := \langle h , x_{(1)} \rangle x_{(2)} ,
\]
where $\Delta(x)=x_{(1)}\otimes x_{(2)}$ is Sweedler notation for the
coproduct, with summation understood.
The finite-dimensional type-1 irreducible $\mathcal U_q(\mathfrak{su}(2))$-modules
are labelled by spins $l\in\frac12\mathbb N_0$. We denote the spin-$l$
module by $V_l$; it has dimension $2l+1$.

The quantum tangent vectors of the
three-dimensional left-covariant calculus on $\mathcal A(\mathrm{SU}_q(2))$ may be chosen as
\begin{equation}
  X_z=\frac{1-K^4}{1-q^{-2}},
  \qquad
  X_-=q^{-1/2}FK,
  \qquad
  X_+=q^{1/2}EK.
  \label{eq:quantum-tangent-vectors}
\end{equation}
The vertical component $X_z$ disappears upon restriction to the fixed-point
algebra, and one obtains the unique left-covariant two-dimensional calculus on
the standard Podle\'s sphere \cite{PodlesCalculus,KLvS}:
\begin{equation}
  \mathrm d f
  =(X_+\triangleright f)\omega_+
   +(X_-\triangleright f)\omega_-,
  \label{eq:podles-differential}
\end{equation}
where
\[
\omega_z = a^*\,da + c^*\,dc,
\qquad
\omega_- = c^*\,da^* - q a^*\,dc^*,
\qquad
\omega_+ = a\,dc - q c\,da.
\]
The complex structure is
\begin{equation}
  \partial f=(X_+\triangleright f)\omega_+,
  \qquad
  \bar\partial f=(X_-\triangleright f)\omega_-.
  \label{eq:podles-dolbeault}
\end{equation}
As bimodules,
\begin{equation}
  \Omega^{1,0}(\mathcal A_q)\simeq\mathcal L_2 \, \omega_+,
  \qquad
  \Omega^{0,1}(\mathcal A_q)\simeq\mathcal L_{-2} \, \omega_-.
  \label{eq:podles-one-forms}
\end{equation}
Moreover,
\[
  \omega_\pm\wedge\omega_\pm=0,
  \qquad
  \omega_-\wedge\omega_+
  +q^{-2}\omega_+\wedge\omega_-=0,
\]
and
\begin{equation}
  \Omega^2(\mathcal A_q)
  =\mathcal A_q\,\omega,
  \qquad
  \omega:=\omega_-\wedge\omega_+,
  \qquad
  \omega f=f\omega.
  \label{eq:podles-two-forms}
\end{equation}
In particular,
\begin{equation}
  \Omega^{0,2}(\mathcal A_q)=0,
  \label{eq:podles-zero-zero-two}
\end{equation}
so every $\bar\partial$-connection on an $\mathcal A_q$-module is flat.

\paragraph{The canonical connection and its holomorphic part.}

Projecting the three-dimensional calculus to the horizontal forms gives a connection on $\mathcal L_n$:
\begin{equation}
  \nabla\phi
  =q^{-n-2}\omega_+(X_+\triangleright\phi)
   +q^{-n+2}\omega_-(X_-\triangleright\phi),
  \qquad \phi\in\mathcal L_n.
  \label{eq:KLvS-full-connection}
\end{equation}
With the notation
\[
  [n]_q=\frac{q^{-n}-q^n}{q^{-1}-q},
\]
its curvature is the equivariant two-form
\begin{equation}
  \nabla^2
  =-q^{-n-1}[n]_q\,\omega_+\wedge\omega_-.
  \label{eq:KLvS-curvature}
\end{equation}
This has type $(1,1)$.  Therefore
\begin{equation}
  \overline\nabla_{(n)}\phi
  :=q^{-n+2}\omega_-(X_-\triangleright\phi)
  \label{eq:standard-KLvS-delbar}
\end{equation}
is flat and defines the \emph{standard KLvS holomorphic structure} on
$\mathcal L_n$ \cite{KLvS}.  For every $n$ there is a twisted flip
\begin{equation}
  \Phi_{(n)}:
  \mathcal L_n\otimes_{\mathcal A_q}\Omega^{0,1}(\mathcal A_q)
  \xrightarrow{\ \simeq\ }
  \Omega^{0,1}(\mathcal A_q)\otimes_{\mathcal A_q}\mathcal L_n,
  \label{eq:KLvS-twisted-flip}
\end{equation}
which makes $\overline\nabla_{(n)}$ a bimodule connection.  Under
\eqref{eq:line-tensor-law}, the tensor product as in \eqref{eq:tensor-product-connection} of the standard structures on
$\mathcal L_n$ and $\mathcal L_m$ is the standard structure on
$\mathcal L_{n+m}$.

\paragraph{Holomorphic functions and sections.}

At the polynomial level KLvS prove
\begin{equation}
  \ker\!\left(
    \bar\partial:\mathcal A_q\longrightarrow\Omega^{0,1}(\mathcal A_q)
  \right)
  =\mathbb C.
  \label{eq:podles-holomorphic-functions}
\end{equation}
They also formulate an analytic domain in the $C^*$-completion,
\[
  \operatorname{Dom}(\bar\partial)
  =\{f\in C(\mathrm{SU}_q(2)):
      \|F\triangleright f\|<\infty\},
\]
and show that the holomorphic elements of
$\operatorname{Dom}(\bar\partial)\cap C(\mathbb{CP}_q^1)$ are again the
constants.  This illustrates a general point: at the $C^*$-level,
$\bar\partial$ is unbounded, so its domain is part of the holomorphic data.

For the spaces of holomorphic sections, they obtain, for $n>0$,
\begin{equation}
  H^0(\mathcal L_n,\overline\nabla_{(n)})=0,
  \qquad
  H^0(\mathcal L_{-n},\overline\nabla_{(-n)})
  \simeq\mathbb C^{n+1}.
  \label{eq:KLvS-H0}
\end{equation}
A basis of the second space is represented by
\begin{equation}
  a^\mu c^{n-\mu},
  \qquad 0\leq\mu\leq n.
  \label{eq:KLvS-holomorphic-basis}
\end{equation}
It is convenient to normalize the line-bundle notation by
\begin{equation}
  \mathcal O_q(k)
  :=(\mathcal L_{-k},\overline\nabla_{(-k)}).
  \label{eq:Oq-notation}
\end{equation}
Then
\begin{equation}
  \dim H^0\bigl(\mathcal O_q(k)\bigr)
  =
  \begin{cases}
    k+1,& k\geq 0,\\
    0,& k<0,
  \end{cases}
  \label{eq:Oq-H0}
\end{equation}
exactly as on the classical $\mathbb{CP}^1$.  Multiplication of sections gives
\begin{equation}
  \bigoplus_{k\geq 0}H^0\bigl(\mathcal O_q(k)\bigr)
  \simeq
  \mathbb C\langle a,c\rangle/(ac-qca),
  \label{eq:quantum-homogeneous-coordinate-ring}
\end{equation}
the homogeneous coordinate algebra of the quantum plane.  The dimensions and
Hilbert series are undeformed, while the multiplication is deformed.  This
relates the KLvS description by modules with connections to the
noncommutative projective-geometric viewpoint of twisted homogeneous coordinate rings
\cite{ArtinTateVdB,ArtinVdB}.

\subsection{Twisted positive fundamental cocycles and a variational principle}

The Haar state $h$ on $\mathcal A(\mathrm{SU}_q(2))$ and on its subalgebra $\mathcal A_q$
is not a trace; let $\sigma$ denote its modular
automorphism.  With
\[
  \int_h b\,\omega_-\wedge\omega_+:=h(b),
\]
the functional
\begin{equation}
 \tau_q(a_0,a_1,a_2)
  =\frac12\int_h a_0\,\mathrm d a_1\wedge\mathrm d a_2
  \label{eq:twisted-fundamental-cocycle}
\end{equation}
is a nontrivial $\sigma$-twisted cyclic two-cocycle, whereas
\begin{equation}
  \varphi_q(a_0,a_1,a_2)
  =\int_h a_0\,\partial a_1\,\bar\partial a_2
  \label{eq:twisted-positive-cocycle}
\end{equation}
is a positive $\sigma$-twisted Hochschild cocycle.  KLvS prove that
$\varphi_q$ and $\tau_q$ are cohomologous in twisted Hochschild cohomology.  Thus, as in the classical case \cite{ConnesCuntz},  the
complex structure is represented by a positive cocycle in the same twisted
Hochschild class as the fundamental cocycle.  The twist records the modularity of
the Haar state and cannot be discarded \cite{KLvS,SchmudgenWagner}.

\paragraph{Connes's variational principle and a conjecture for the quantum sphere.}

The positive cocycle associated with a complex structure has a variational
interpretation which is particularly transparent in complex dimension one.
Let $\mathcal A$ be a unital $*$-algebra and let
$Z^2_+(\mathcal A,\mathcal A^*)$ denote the convex cone of positive Hochschild
$2$-cocycles.  Thus $\Phi\in Z^2_+(\mathcal A,\mathcal A^*)$ when
\[
  \bigl\langle a_0\otimes a_1,b_0\otimes b_1\bigr\rangle_\Phi
  :=\Phi\bigl(b_0^*a_0,a_1,b_1^*\bigr)
\]
defines a positive sesquilinear form on $\mathcal A\otimes\mathcal A$.
For an oriented compact surface, Connes showed that the positive cocycles in the
cyclic affine class of the fundamental cocycle form a convex set whose extreme
points encode conformal, or equivalently complex, structures.  More precisely, a
positive element of the universal $2$-forms exposes the point corresponding to the
complex structure determined by the associated metric
\cite[Chapter~VI, Section~2]{ConnesNCG}.

For the noncommutative two-torus of
Section~\ref{subsec:complex-star-algebras}, consider first the square complex
structure, with parameter $i$, and put
\[
  \partial_i=\frac12(\delta_1-i\delta_2),
  \qquad
  \bar\partial_i=\frac12(\delta_1+i\delta_2),
\]
and consider the positive Hochschild cocycle
\begin{equation}
  \Phi_i(a_0,a_1,a_2)
  =\frac{2}{\pi}\operatorname{tr}_\theta
       \bigl(a_0\,\partial_i(a_1)\bar\partial_i(a_2)\bigr)
  =\frac{1}{2\pi}\operatorname{tr}_\theta
       \bigl(a_0(\delta_1-i\delta_2)(a_1)
                    (\delta_1+i\delta_2)(a_2)\bigr).
  \label{eq:connes-square-positive-cocycle}
\end{equation}
Let $\tau_2$ be the cyclic cocycle representing the oriented fundamental class,
with the orientation chosen so that \eqref{eq:connes-square-positive-cocycle}
belongs to its cyclic affine class.  If $b$ and $B$ denote the Hochschild and
Connes coboundary operators, define
\[
  \mathscr C_\theta^+
  :=Z^2_+(A_\theta^\infty,(A_\theta^\infty)^*)
       \cap\bigl(\tau_2+b(\ker B)\bigr).
\]
The affine slice $\tau_2+b(\ker B)$ fixes the cyclic orientation class,
not merely the Hochschild cohomology class.
For convenience write $U=U_1$ and $V=U_2$.
\begin{theorem}[Connes's variational characterization]
\label{thm:connes-positive-minimum}
For the positive universal $2$-form
\begin{equation}
  G_0=dU(dU)^*+dV(dV)^*,
  \label{eq:flat-torus-positive-metric}
\end{equation}
the linear functional
\[
  \mathcal E_{G_0}(\Phi)
  :=\langle G_0,\Phi\rangle
  =\Phi(1,U,U^*)+\Phi(1,V,V^*)
\]
has a unique minimum on $\mathscr C_\theta^+$.  Its minimizer is the cocycle
$\Phi_i$ in \eqref{eq:connes-square-positive-cocycle}, and
\begin{equation}
  \mathcal E_{G_0}(\Phi)\geq 4\pi,
  \qquad
  \mathcal E_{G_0}(\Phi)=4\pi
  \quad\Longleftrightarrow\quad
  \Phi=\Phi_i.
  \label{eq:connes-four-pi-minimum}
\end{equation}
In particular, $\Phi_i$ is an exposed, and hence extreme, point of this
positive affine slice.
\end{theorem}

\noindent
Indeed,
\[
  \Phi_i(1,U,U^*)=2\pi,
  \qquad
  \Phi_i(1,V,V^*)=2\pi,
\]
which gives the value $4\pi$. The uniqueness assertion identifies the complex
structure through the metric element $G_0$. A constant metric with parameter
$\tau\in\mathbb H$ as in \eqref{eq:torus-dolbeault} similarly selects the
normalized cocycle $\widehat\Phi_\tau$ in the fixed fundamental cyclic affine
class. For $\tau=i$, the normalized and unnormalized cocycles agree.
An orientation convention must be fixed consistently in the cyclic class
and in the positive cocycle.

For the quantum two-sphere
with cocycles $\tau_q$ and $\varphi_q$ as in
\eqref{eq:twisted-fundamental-cocycle} and \eqref{eq:twisted-positive-cocycle},
KLvS prove that $\varphi_q$ is positive in the $\sigma$-twisted sense and that
\begin{equation}
  \tau_q-\varphi_q=b_\sigma\psi_q,
  \qquad
  \psi_q(a,b)=\frac12\int_h a\,\partial\bar\partial b.
  \label{eq:KLvS-transgression-variational}
\end{equation}
The twisted cyclic condition must again be retained.  This is especially important
here, since for the modular automorphism $\sigma$ the corresponding twisted
Hochschild cohomology in degree two is trivial, whereas $\tau_q$ is nontrivial in
twisted cyclic cohomology.  Thus the phrase ``in the same Hochschild class'' would
be too weak and would lead to a degenerate variational problem.

Let $B_\sigma$ denote the twisted Connes operator and let
$Z^2_{\sigma,+}(\mathcal A_q,\mathcal A_q^*)$ be the cone of $\sigma$-twisted
positive Hochschild $2$-cocycles.  The natural analogue of Connes's affine slice is
\begin{equation}
  \mathscr C_{q,+}^{\mathrm{cov}}
  :=Z^2_{\sigma,+}(\mathcal A_q,\mathcal A_q^*)^{\mathrm{cov}}
    \cap\bigl(\tau_q+b_\sigma(\ker B_\sigma)\bigr),
  \label{eq:twisted-positive-affine-slice}
\end{equation}
where the superscript indicates covariance under the natural
$\mathcal U_q(\mathfrak{su}(2))$-symmetry.  Equivalently, this is the covariant
positive part of the affine family having the same twisted cyclic orientation
class as $\tau_q$.  For the normalized twisted $(b,B)$-complex, the KLvS
transgression satisfies
\[
  (B_\sigma\psi_q)(a)
  =\psi_q(1,a)+\psi_q(\sigma(a),1)
  =\frac12\int_h\partial\bar\partial a=0,
\]
where the last equality follows from the twisted Stokes property.  Hence
\eqref{eq:KLvS-transgression-variational} places $\varphi_q$ in the affine
family \eqref{eq:twisted-positive-affine-slice}.

Let
\[
  G_q\in\Omega^2_{+,\mathrm{univ}}(\mathcal A_q)
\]
be a positive, covariant universal $2$-form in Connes's sense.  We require the
associated Hermitian form on universal $1$-forms to descend to the KLvS calculus
and to give the standard invariant Hermitian metric compatible with
$\Omega^1(\mathcal A_q)=\Omega^{1,0}(\mathcal A_q)\oplus\Omega^{0,1}(\mathcal A_q)$.  A scale must still be fixed.  For a
positive presentation
$G_q=\sum_\alpha dx_\alpha(dx_\alpha)^*$, set
\begin{equation}
  \mathcal E_q(\Phi)
  :=\langle G_q,\Phi\rangle
  =\sum_\alpha\Phi(1,x_\alpha,x_\alpha^*).
  \label{eq:podles-positive-functional}
\end{equation}

\begin{conjecture}[Twisted variational characterization of the KLvS complex structure]
\label{conj:KLvS-variational-characterization}
For every $0<q<1$, the functional $\mathcal E_q$ has a unique minimum on
$\mathscr C_{q,+}^{\mathrm{cov}}$, and its minimizer is the KLvS cocycle
$\varphi_q$ in \eqref{eq:twisted-positive-cocycle}.  Equivalently,
\begin{equation}
  \mathcal E_q(\Phi)
  \geq \mathcal E_q(\varphi_q)
  \qquad
  \text{for all }\Phi\in\mathscr C_{q,+}^{\mathrm{cov}},
  \label{eq:KLvS-variational-inequality}
\end{equation}
with equality if and only if $\Phi=\varphi_q$. Thus $\varphi_q$ is an exposed,
and therefore extreme, point of the covariant twisted-positive affine slice.
\end{conjecture}

A stronger version would remove the covariance restriction in
\eqref{eq:twisted-positive-affine-slice}.  The numerical minimum in
\eqref{eq:KLvS-variational-inequality} depends on the normalization of $G_q$; the
normalization-independent content is the ratio
\[
  \frac{\mathcal E_q(\Phi)}{\mathcal E_q(\varphi_q)}\geq 1.
\]
Once a canonical scale for the invariant metric has been fixed, one may further
ask whether the minimum is undeformed under $q$-deformation, or more generally
compute its precise $q$-dependence and its classical limit.  This conjecture
concerns rigidity of the complex structure on the base quantum sphere. It
should not be confused with the classification of holomorphic structures
on a fixed line module $\mathcal L_n$ under a specified gauge group;
polynomial line-rigidity fails by
Proposition~\ref{prop:polynomial-line-nonrigidity}.

\paragraph{Conformal changes and the persistence of the positive cocycle.}

Conformal changes distinguish the positive cocycle from the volume weight.
On an oriented classical surface, a conformal change
\[
  g_h=e^{-h}g
\]
does not alter the complex structure.  Equivalently, the Dirichlet energy of
functions is conformally invariant in real dimension two, since
\[
  \sqrt{\det g_h}\,g_h^{-1}
  =\sqrt{\det g}\,g^{-1}.
\]
Consequently, the positive Hochschild cocycle
\[
  \Phi_J(f^0,f^1,f^2)
  =\frac{i}{\pi}\int_M f^0\,\partial_J f^1\wedge\bar\partial_J f^2
\]
depends on the complex structure $J$, or equivalently on the conformal class,
but not on a particular metric in that class.

The construction of curved noncommutative two-tori implements exactly this
separation \cite{ConnesTretkoff,ConnesMoscoviciModular,FathizadehKhalkhaliGB}.
Fix the complex parameter $\tau\in\mathbb H$ as in \eqref{eq:torus-dolbeault} and the corresponding positive
Hochschild cocycle $\Phi_\tau$ of \eqref{eq:torus-positive-cocycle}.  For a
self-adjoint element $h\in A_\theta^\infty$, set
\begin{equation}
  \varphi_h(a) :=\operatorname{tr}_\theta(ae^{-h}).
  \label{eq:conformal-weight-torus}
\end{equation}
The functional $\varphi_h$ is a positive KMS weight, generally nontracial.  Its modular operator on the algebra and the corresponding
one-parameter automorphism group are
\begin{equation}
  \Delta_h(a)=e^{-h}ae^h,
  \qquad
  \sigma_t^h(a)=e^{ith}ae^{-ith}.
  \label{eq:modular-automorphism-torus}
\end{equation}
Let $\mathcal H_h^0$ be the GNS Hilbert space of $\varphi_h$, and let
$\mathcal H_\tau^{1,0}$ be the Hilbert space of $(1,0)$-forms obtained from the
positive cocycle $\Phi_\tau$.  The crucial point is that
$\mathcal H_\tau^{1,0}$ is not changed by the Weyl factor.  One uses the same
algebraic Dolbeault operator
\begin{equation}
  \partial_{\tau,h}:A_\theta^\infty\subset\mathcal H_h^0
  \longrightarrow \mathcal H_\tau^{1,0},
  \qquad
  \partial_{\tau,h}(a)=\partial_\tau(a),
  \label{eq:unchanged-dolbeault-conformal}
\end{equation}
and changes only its Hilbert-space adjoint.

\begin{theorem}[Conformal invariance of the positive cocycle]
\label{thm:positive-cocycle-conformal-invariance}
Let
\[
  \Delta_{\tau,h}
  :=\partial_{\tau,h}^{\,*}\partial_{\tau,h}
\]
be the Laplacian on functions for the conformal weight $\varphi_h$.  Then, for
$a\in A_\theta^\infty$,
\begin{equation}
  \langle a,\Delta_{\tau,h}a\rangle_{\varphi_h}
  =\|\partial_\tau a\|_{\mathcal H_\tau^{1,0}}^2
  =\frac{\pi}{2}\,\Phi_\tau(1,a,a^*).
  \label{eq:dirichlet-form-conformal-invariance}
\end{equation}
In particular, as algebraic data, the Dirichlet form, the Hilbert space of
$(1,0)$-forms and the positive Hochschild cocycle are independent of $h$.
The Weyl factor changes the volume weight, the adjoint of $\partial_\tau$, and
so the Laplacian as an operator, but it does not change the complex
structure encoded by $\Phi_\tau$.
\end{theorem}

Indeed, the first equality in
\eqref{eq:dirichlet-form-conformal-invariance} is the definition of the
adjoint, and the second follows from the inner product determined by
$\Phi_\tau$.  After transporting the operator to the flat Hilbert space, one
finds that $\Delta_{\tau,h}$ is antiunitarily equivalent to
\begin{equation}
  k\Delta_{\tau,0}k , \qquad k:=e^{h/2}.
  \label{eq:conformal-laplacian-k-delta-k}
\end{equation}
Thus the same algebraic quadratic form is represented by different Laplace
operators because the $L^2$-structure on functions has changed.

This suggests that, in complex dimension one, a metric in the conformal class
determined by $\tau$ is naturally packaged by the pair
\begin{equation}
  \bigl(\Phi_\tau,\varphi_h\bigr),
  \label{eq:metric-as-positive-cocycle-and-weight}
\end{equation}
where $\Phi_\tau$ records the complex, or conformal, structure and
$\varphi_h$ records the volume density.  The noncentrality of $h$ is precisely
what turns the latter into modular data.  Curvature, the determinant of the
Laplacian and the Polyakov anomaly depend on $h$, even though $\Phi_\tau$ does
not. In particular, Connes and Moscovici prove that the scale-invariant
determinant assumes its extreme value only at the flat metric
\cite{ConnesMoscoviciModular}.

Replacing the trace in $\Phi_\tau$ by the conformal weight gives a cochain
\begin{equation}
  \widetilde\Phi_{\tau,h}(a_0,a_1,a_2)
  :=\frac{2}{\pi}\,
  \varphi_h\bigl(a_0\partial_\tau(a_1)\bar\partial_\tau(a_2)\bigr)
  \label{eq:naive-weighted-positive-cochain}
\end{equation}
which is generally not an ordinary Hochschild cocycle; ordinary positivity
is not automatic either. The Leibniz rules give, with
$X=a_0\partial_\tau(a_1)\bar\partial_\tau(a_2)$,
\[
 (b\widetilde\Phi_{\tau,h})(a_0,a_1,a_2,a_3)
 =\frac{2}{\pi}\,\varphi_h(Xa_3-a_3X).
\]
This need not vanish because $\varphi_h$ is generally nontracial. A
metric-dependent cocycle obtained by absorbing the weight therefore
requires a compatible modular or twisted construction.

The preceding theorem concerns Dirichlet forms and spectral triples. A
stronger variational question asks whether the same cocycle minimizes the
functional associated with a conformally perturbed universal metric.  Let
$G_{\tau,h}\in\Omega^2_+(A_\theta^\infty)$ denote a canonical positive
universal $2$-form representing the metric obtained from the flat metric of
modulus $\tau$ by the Weyl factor $e^{-h}$.  The word \emph{canonical} is
important here, because in the noncommutative case there are ordering choices
in placing the factor $e^{-h}$ inside a universal metric form.
For this variational question we use the normalization
\eqref{eq:normalized-torus-positive-cocycle} and set
\begin{equation}
  \mathscr C_{\theta,\tau}^+
  :=\bigl\{\Phi\in Z^2_+(A_\theta^\infty,(A_\theta^\infty)^*):
             \Phi-\widehat\Phi_\tau\in b(\ker B)\bigr\}.
  \label{eq:fixed-positive-cyclic-class-tau}
\end{equation}
This is the positive part of the fixed fundamental cyclic affine class;
it agrees with the class containing $\widehat\Phi_\tau$ for every $\tau$.

\begin{conjecture}[Conformal stability of Connes's minimizer]
\label{conj:conformal-stability-connes-minimizer}
For every self-adjoint $h\in A_\theta^\infty$, the functional
\[
  \Phi\longmapsto\langle G_{\tau,h},\Phi\rangle
\]
on $\mathscr C_{\theta,\tau}^+$ has the same unique minimizer
$\widehat\Phi_\tau$ as the flat metric.  Equivalently, Connes's variational extraction of the complex
structure is invariant under noncommutative Weyl perturbations.
\end{conjecture}

For $h$ scalar, the conjecture reduces immediately to the flat result, since
the metric functional is merely rescaled.  Its content concerns genuinely
noncentral Weyl factors.  The curved-torus results prove the persistence of
the Dolbeault quadratic form in
\eqref{eq:dirichlet-form-conformal-invariance}, but do not by themselves
identify and solve the corresponding minimization problem for a
noncommutative universal metric element.

The same distinction is useful for the quantum projective line.  It suggests
that a future conformal perturbation of the KLvS geometry should retain the
KLvS twisted-positive cocycle as the complex-structure datum and place the
Weyl factor in a separate KMS weight.  Since the Haar functional is already
modular, a further Weyl perturbation changes the modular data by a cocycle
perturbation.  A precise conformal-stability statement would therefore require
a canonical comparison of the corresponding twisted cyclic complexes.  In
that form it becomes a conformal refinement of
Conjecture~\ref{conj:KLvS-variational-characterization}.

\subsection{The quantum Birkhoff--Grothendieck problem}
\label{subsec:quantum-BG-conjecture}

The classical Birkhoff--Grothendieck theorem states that every holomorphic vector
bundle $E$ on $\mathbb{CP}^1$ has a decomposition
\begin{equation}
  E\simeq
  \mathcal O(k_1)\oplus\cdots\oplus\mathcal O(k_r),
  \qquad
  k_1\geq\cdots\geq k_r,
  \label{eq:classical-BG}
\end{equation}
and that the ordered integers $k_j$ are uniquely determined
\cite{Birkhoff,Grothendieck,HazewinkelMartin}.  The algebraic and holomorphic
versions agree by GAGA \cite{SerreGAGA}.  For example, tensoring the tangent Euler sequence on the classical
projective line by $\mathcal O(-1)$ gives
\begin{equation}
  0\longrightarrow\mathcal O(-1)
  \longrightarrow\mathcal O^{\oplus 2}
  \longrightarrow\mathcal O(1)
  \longrightarrow 0
 \label{eq-Euler-Kronecker} .
 \end{equation}
The sequence is non-split, although its middle term is already a direct sum of line bundles.
Notice that the Birkhoff--Grothendieck theorem classifies the
middle bundle, not the way it may occur in a particular extension sequence.

For the quantum projective line one must specify the category before stating a
conjecture.  We shall use either the polynomial algebra $\mathcal A_q$ or a fixed
smooth Fr\'echet completion $\mathcal A_q^\infty$ stable under the
$\mathcal U_q(\mathfrak{su}(2))$-action and the Dolbeault operators.  Write
\[
  \operatorname{Hol}_{\mathrm{fgp}}(\mathcal A_q^{\mathscr R})
\]
for the category of finitely generated projective left modules in the regularity
class $\mathscr R\in\{\mathrm{pol},\infty\}$, equipped with flat
$\bar\partial$-connections of the same regularity.  A $C^*$-algebraic category
requires explicit domains, closedness and compatibility with gauge
transformations, and is therefore not included without further hypotheses.

Three logically distinct questions are involved:
\begin{enumerate}[label=\textup{(\roman*)},leftmargin=2.3em]
\item Does the underlying projective module decompose into canonical rank-one
      modules $\mathcal L_n$?
\item Can such a decomposition be chosen to intertwine the
      $\bar\partial$-connections?
\item Is every holomorphic structure on $\mathcal L_n$ gauge equivalent to the
      standard KLvS structure?
\end{enumerate}
The third question is independent of the first two.  Since
$\Omega^{0,2}(\mathcal A_q)=0$, every perturbation of a line-module
$\bar\partial$-connection is flat; the issue is entirely gauge-theoretic.

\begin{conjecture}[Holomorphic splitting conjecture]
\label{conj:weak-BG}
For $\mathscr R=\mathrm{pol}$, and separately for an appropriate smooth category
$\mathscr R=\infty$, every object
$(\mathcal E,\overline\nabla)$ of
$\operatorname{Hol}_{\mathrm{fgp}}(\mathcal A_q^{\mathscr R})$ is
holomorphically isomorphic to a finite direct sum
\begin{equation}
  (\mathcal E,\overline\nabla)
  \simeq
  \bigoplus_{j=1}^r
  (\mathcal L_{n_j},\overline\nabla_j),
  \label{eq:weak-BG}
\end{equation}
where each $\overline\nabla_j$ is a holomorphic structure on the canonical line module
$\mathcal L_{n_j}$.  The resulting rank-one summands are unique up to permutation
and holomorphic isomorphism.
\end{conjecture}

This is the form of the Birkhoff--Grothendieck problem that does not presuppose
line-bundle rigidity.  It remains meaningful even if a fixed topological line
module carries several gauge-inequivalent holomorphic structures.

\begin{proposition}[Failure of line-bundle rigidity for polynomial gauges]
\label{prop:polynomial-line-nonrigidity}
The polynomial coordinate algebra $\mathcal A_q$ of the quantum projective line has unit group
$\mathcal A_q^\times=\mathbb C^\times$. Consequently, for every nonzero
$\eta\in\Omega^{0,1}(\mathcal A_q)$, the flat connection
\[
 \overline\nabla_\eta(a)=\bar\partial a+a\eta,
 \qquad a\in\mathcal L_0=\mathcal A_q,
\]
is not polynomially gauge equivalent to the standard connection
$\bar\partial$.
\end{proposition}

\begin{proof}
Put $z=B_0$, $x=B_-$ and $y=B_+$. The relations
\eqref{eq:podles-relations-one}--\eqref{eq:podles-relations-two} give the
rank-one generalized Weyl algebra over $\mathbb C[z]$ with
\[
 xz=q^2zx,\qquad yz=q^{-2}zy,\qquad
 yx=z(1-z),\qquad xy=q^2z(1-q^2z).
\]
Its usual normal forms are
$\mathbb C[z]$, $\mathbb C[z]x^k$ and $\mathbb C[z]y^k$ for $k\geq1$.
They can also be obtained by embedding the algebra into
$\mathbb C(z)[t,t^{-1};\sigma]$, where
$\sigma(z)=q^2z$, by $x\mapsto t$ and
$y\mapsto z(1-z)t^{-1}$. In particular, the algebra is a domain and is
$\mathbb Z$-graded by $\deg x=1$, $\deg y=-1$, $\deg z=0$.
For nonzero elements the highest and lowest grading degrees are additive
under multiplication. If $uv=1$, their degree widths must therefore
both vanish, so $u$ and $v$ are homogeneous. A nonzero-degree homogeneous
product with an opposite-degree element contains a nonconstant factor
from $x^ky^k$ or $y^kx^k$ and cannot equal $1$. Degree-zero units are
units of $\mathbb C[z]$, hence nonzero scalars.

Every automorphism of the left regular module $\mathcal A_q$ is right
multiplication by a unit. Scalar automorphisms leave
$\overline\nabla_\eta$ unchanged. The connection has the left Leibniz
rule and is flat because $\Omega^{0,2}(\mathcal A_q)=0$; for
$\eta\neq0$ it is distinct from the standard connection. For example,
one may take $\eta=\bar\partial B_0\neq0$.
\end{proof}

This elementary polynomial obstruction depends on the gauge group.
The broader uniqueness question was raised in \cite{KLvS}.
Enlarging the gauge group changes the problem: an invertible element of a
smooth or $C^*$-completion need not be polynomial. One must specify the
completion and require both the gauge and its inverse to preserve the
Dolbeault domain before comparing the resulting structures.

A distinct uniqueness result holds in the equivariant category. D\'iaz
Garc\'ia, Krutov, \'{O} Buachalla, Somberg and Strung construct covariant
holomorphic structures on finitely generated relative Hopf modules over
irreducible quantum flag manifolds with Heckenberger--Kolb calculi, and
prove uniqueness for simple relative Hopf modules
\cite{DiazEtAl,DiazEtAlCorrection}. For the canonical line objects on the
Podle\'s sphere this gives covariant rigidity. It does not classify
arbitrary non-equivariant $\bar\partial$-connections.

\begin{remark}[Completed gauges and nonstandard line structures]
\label{rem:recent-preprints}
Graveman, La Rue, MacArthur, Pesin and Wei \cite{GravemanEtAl} work with
an enlarged gauge group consisting of invertible elements of
$C(\mathbb{CP}_q^1)$ for which both the element and its inverse lie in
$\operatorname{Dom}(\bar\partial)$. Their Theorem~3.27 gives a criterion
for a polynomial coefficient determined by
$f\in\mathbb C[B_0]$ to be removable by such a gauge. In their terminology,
$f=B_0/(1-q^2)$ has a defective spot and provides an obstruction even
for this completed gauge group. Their results therefore concern a
stronger equivalence relation than polynomial gauge equivalence; the
obstruction also persists for any smooth gauge group contained in their
admissible completed group. They further construct infinitely many
pairwise non-gauge-equivalent holomorphic line structures. The related
work of Biswas, Guin and Kumar \cite{BiswasGuinKumar} develops a
fixed-point criterion for gauge equivalence. These results reinforce the
need to distinguish splitting into holomorphic rank-one summands from
identifying each summand with a standard KLvS line. The proofs of the
splitting results below do not depend on these analytic constructions.
\end{remark}

\begin{remark}[Standard-line decompositions and their scope]
\label{rem:standard-line-decomposition-scope}
A decomposition into the standard KLvS lines has the form
\begin{equation}
  (\mathcal E,\overline\nabla)
  \simeq
  \mathcal O_q(k_1)\oplus\cdots\oplus\mathcal O_q(k_r),
  \qquad
  k_1\geq\cdots\geq k_r.
  \label{eq:strong-BG}
\end{equation}
This is not a conjecture for all polynomial holomorphic bundles:
Proposition~\ref{prop:polynomial-line-nonrigidity} already excludes that assertion in
rank one. The weak splitting conjecture asks instead for holomorphic
rank-one summands whose line structures may be nonstandard. In this
paper a standard-line decomposition is proved for the polynomial
bundles admitting a finite holomorphic filtration with standard KLvS
line quotients; see Theorem~\ref{thm:line-filtered-splitting}. No universal
line-rigidity assertion in an unspecified smooth category is assumed.
\end{remark}

The dimension formula \eqref{eq:Oq-H0} gives a cohomological test whenever
\eqref{eq:strong-BG} holds. To twist a left module, we tensor on the left
by a standard line bimodule, using its KLvS bimodule connection. This
convention does not require a right-module structure on the given bundle.
For a standard-line decomposition,
\begin{equation}
  h_{\mathcal E}(m)
  :=
  \dim H^0\bigl(\mathcal O_q(m)\otimes_{\mathcal A_q}\mathcal E\bigr)
  =
  \sum_{j=1}^r\max\{m+k_j+1,0\}.
  \label{eq:Hilbert-function-splitting}
\end{equation}
Consequently, the second finite difference
\[
  h_{\mathcal E}(m)-2h_{\mathcal E}(m-1)+h_{\mathcal E}(m-2)
\]
counts the multiplicity of the integer $-m$ among the $k_j$. In
particular, a standard splitting type is unique whenever such a
decomposition exists. This formula cannot be applied to arbitrary
nonstandard line structures without computing their cohomology.

The following sections develop the Dolbeault-module analysis. Its comparison
with the sheaf-theoretic approach requires a quantum GAGA statement relating
locally free holomorphic sheaves to projective modules with
$\bar\partial$-connections; no such comparison is assumed here.

\section{Rank-two test objects for a quantum Grothendieck theorem}
\label{sec:rank-two-test-objects}

We construct holomorphic bundles whose presentations are not diagonal and
examine the distinction between an extension class and the isomorphism class
of its middle term. Upper-triangular connections provide the main examples;
Clifford and Bott projectors give complementary descriptions of the underlying
projective modules. Sections~\ref{sec:rank-two-kronecker}
and~\ref{sec:all-line-filtered-examples-split} determine the decompositions of
the resulting line-filtered objects.

\paragraph{Conventions for the line modules.}

We retain the algebra $\mathcal A_q$, the line bimodules $\mathcal L_n$,
and the connections $\overline\nabla_{(n)}$ of
Section~\ref{subsec:podles-klv}. Unless stated otherwise, $\mathcal L_n$
includes this holomorphic structure, so that
$\mathcal O_q(k)=\mathcal L_{-k}$. Multiplication identifies
$\mathcal L_n\otimes_{\mathcal A_q}\mathcal L_m$ with $\mathcal L_{n+m}$
and intertwines these connections through the KLvS twisted flip
\cite{KLvS,KhalkhaliMoatadelro,KhalkhaliMoatadelroCoordinate}. We also use
\[
  \Omega^{0,1}(\mathcal A_q)
  \simeq\mathcal L_{-2}\omega_-
  \simeq\omega_-\mathcal L_{-2}.
\]
The vanishing $\Omega^{0,2}(\mathcal A_q)=0$ makes all the
$\bar\partial$-connections below flat.

\subsection{A general upper-triangular construction}
\label{subsec:upper-triangular-construction}

Fix integers $a,b$ and set $\mathcal E_{a,b}=\mathcal L_a\oplus\mathcal L_b$.
For an $\mathcal A_q$-linear map
\[
  \beta:\mathcal L_b\longrightarrow
  \Omega^{0,1}(\mathcal A_q)\otimes_{\mathcal A_q}\mathcal L_a,
\]
define
\[
  \overline\nabla_\beta
  =\begin{pmatrix}
     \overline\nabla_{(a)}&\beta\\
     0&\overline\nabla_{(b)}
   \end{pmatrix},
  \qquad
  \overline\nabla_\beta\begin{pmatrix}\xi\\\eta\end{pmatrix}
  =\begin{pmatrix}\overline\nabla_{(a)}\xi+\beta(\eta)\\
                  \overline\nabla_{(b)}\eta\end{pmatrix}.
\]
The $\mathcal A_q$-linearity of $\beta$ ensures the left Leibniz rule,
and the curvature vanishes because its  values lie in
$\Omega^{0,2}(\mathcal A_q)\otimes_{\mathcal A_q}\mathcal E_{a,b}$.
The first summand is a holomorphic subbundle, with the second as quotient:
\begin{equation}
  0\longrightarrow(\mathcal L_a,\overline\nabla_{(a)})
  \longrightarrow(\mathcal E_{a,b},\overline\nabla_\beta)
  \longrightarrow(\mathcal L_b,\overline\nabla_{(b)})\longrightarrow0.
  \label{triang-extension}
\end{equation}
The underlying module is a direct sum, but its displayed splitting need
not intertwine the connections. The off-diagonal term $\beta$ measures the
defect of the chosen splitting; its cohomology class will measure the
obstruction to splitting the specified extension holomorphically.

\paragraph{The defect of a chosen splitting.}
Throughout, the end terms
$\mathcal L_a$ and $\mathcal L_b$ carry their standard holomorphic
structures, and all maps are taken in the same category as in the
preceding subsections.
Write the extension in \eqref{triang-extension} as
\[
0\longrightarrow\mathcal L_a
\xrightarrow{\,\iota\,}
(\mathcal E_{a,b},\overline\nabla_\beta)
\xrightarrow{\,\pi\,}\mathcal L_b
\longrightarrow0,
\qquad
\iota(\xi)=(\xi,0),\quad \pi(\xi,\eta)=\eta.
\]
An $\mathcal A_q$-linear right inverse of $\pi$ is
$s_0(\eta)=(0,\eta)$. It need not be holomorphic: directly from the
matrix defining $\overline\nabla_\beta$, one obtains
\[
\overline\nabla_\beta s_0(\eta)
-(1\otimes s_0)\overline\nabla_{(b)}\eta
=
\begin{pmatrix}\beta(\eta)\\0\end{pmatrix}.
\]
Thus $\beta$ is precisely the failure of this particular lift of the
quotient to be holomorphic. A different lift may have a different
defect, which is why $\beta$ itself is not the invariant of the extension.

\paragraph{The induced Dolbeault differential.}
For $h\in\operatorname{Hom}_{\mathcal A_q}(\mathcal L_b,\mathcal L_a)$, define
\[
D_{a,b}h
:=\overline\nabla_{(a)}\circ h
-(1\otimes h)\circ\overline\nabla_{(b)}.
\]
The two Leibniz terms cancel, so $D_{a,b}h$ is an $\mathcal A_q$-linear map
from $\mathcal L_b$ to
$\Omega^{0,1}(\mathcal A_q)\otimes_{\mathcal A_q}\mathcal L_a$, of the same type as
$\beta$. In particular, $D_{a,b}h=0$ says exactly that $h$ is a
holomorphic map from $\mathcal L_b$ to $\mathcal L_a$.
The invertibility of the line bimodules and the KLvS twisted flips
identify the two spaces in this differential with
\[
\begin{aligned}
\operatorname{Hom}_{\mathcal A_q}(\mathcal L_b,\mathcal L_a)
&\simeq\mathcal L_{a-b},\\
\operatorname{Hom}_{\mathcal A_q}\bigl(\mathcal L_b,
\Omega^{0,1}(\mathcal A_q)\otimes_{\mathcal A_q}\mathcal L_a\bigr)
&\simeq\Omega^{0,1}(\mathcal A_q)\otimes_{\mathcal A_q}\mathcal L_{a-b}.
\end{aligned}
\]
The degree difference $a-b$ comes from the dual of the source line
module. More concretely, put $k=a-b$. An element
$u\in\mathcal L_k$ represents the left $\mathcal A_q$-linear map
$R_u(\eta)=\eta u$ from $\mathcal L_b$ to $\mathcal L_a$;
every such map is obtained in this way. The order of multiplication
is important.

For a form-valued parameter
$\vartheta\in\Omega^{0,1}({\mathcal A_q})\otimes_{\mathcal A_q}\mathcal L_k$, the
corresponding map is
\[
\beta_\vartheta(\eta)
=(1\otimes m_{b,k})(\Phi_{(b)}\otimes1)
(\eta\otimes\vartheta),
\]
where $m_{b,k}:\mathcal L_b\otimes_{\mathcal A_q}\mathcal L_k
\to\mathcal L_a$ is multiplication and
$\Phi_{(b)}:\mathcal L_b\otimes_{\mathcal A_q}\Omega^{0,1}(\mathcal A_q)
\to\Omega^{0,1}(\mathcal A_q)\otimes_{\mathcal A_q}\mathcal L_b$
is the twisted flip of Section~\ref{subsec:podles-klv}. Compatibility of multiplication
with the induced tensor-product connections gives
\[
D_{a,b}R_u=\beta_{\overline\nabla_{(k)}u}.
\]
This explains the induced-connection notation
$\overline\nabla_{(a-b)}$ used for the Hom line module: under the displayed
identifications, $D_{a,b}$ is its Dolbeault differential. The twisted
flip, rather than an ordinary interchange of factors, is part of this
identification.

\paragraph{Changing the splitting.}
Every $\mathcal A_q$-linear right inverse of $\pi$ has the form
$s_h(\eta)=(h(\eta),\eta)$ for some
$h\in\operatorname{Hom}_{\mathcal A_q}(\mathcal L_b,\mathcal L_a)$.
The corresponding change of the direct-sum presentation is
\begin{equation}
g_h=\begin{pmatrix}1&h\\0&1\end{pmatrix},
\qquad
g_h^{-1}=\begin{pmatrix}1&-h\\0&1\end{pmatrix}.
\label{gautra}
\end{equation}
No invertibility assumption on $h$ is needed. Moreover,
$g_h\iota=\iota$ and $\pi g_h=\pi$: this change acts as the identity
on the specified subbundle and on the specified quotient.
Applying the connection and then $1\otimes g_h^{-1}$ gives the
gauge-transformed connection
\[
\begin{aligned}
(1\otimes g_h^{-1}) \, \overline\nabla_\beta \, g_h
\begin{pmatrix}\xi\\\eta\end{pmatrix}
&=
\begin{pmatrix}
\overline\nabla_{(a)}\xi+\overline\nabla_{(a)}(h(\eta))
+\beta(\eta)-(1\otimes h)\overline\nabla_{(b)}\eta\\
\overline\nabla_{(b)}\eta
\end{pmatrix}\\[3pt]
&=
\begin{pmatrix}
\overline\nabla_{(a)}\xi+(\beta+D_{a,b}h)(\eta)\\
\overline\nabla_{(b)}\eta
\end{pmatrix}.
\end{aligned}
\]
Consequently, changing the splitting changes $\beta$ by the exact
term $D_{a,b}h$. In the line-module description, with $h=R_u$, this
is the rule
$\vartheta\mapsto\vartheta+\overline\nabla_{(a-b)}u$.

Since $\Omega^{0,2}(\mathcal A_q)=0$, every degree-one cochain in this Hom
complex is closed. Its degree-one cohomology is therefore
\[
\frac{\operatorname{Hom}_{\mathcal A_q}\bigl(\mathcal L_b,
\Omega^{0,1}(\mathcal A_q)\otimes_{\mathcal A_q}\mathcal L_a\bigr)}
{D_{a,b}\bigl(\operatorname{Hom}_{\mathcal A_q}
(\mathcal L_b,\mathcal L_a)\bigr)}
\simeq
H^{0,1}(\mathcal L_{a-b})
=
\frac{\Omega^{0,1}(\mathcal A_q)\otimes_{\mathcal A_q}\mathcal L_{a-b}}
{\overline\nabla_{(a-b)}(\mathcal L_{a-b})}.
\]
We denote the resulting class by $[\beta]$, using the above
identification of $\beta$ with a line-module-valued form.

Here equivalence of extensions means a holomorphic isomorphism of
the middle terms that commutes with the given inclusion and
projection and induces the identity on both end terms. Every such
module isomorphism has the form $g_h$ in \eqref{gautra}.

Thus two upper-triangular
connections define equivalent extensions precisely when their
off-diagonal entries determine the same cohomology class. Conversely,
since $\mathcal L_b$ is projective, any extension with these end
terms splits as an $\mathcal A_q$-module sequence; choosing such a splitting
puts its connection in the upper-triangular form above. Hence the
cohomology group parametrizes extensions with identified end terms,
not merely a chosen family of matrix presentations.

\paragraph{Splitting the extension versus decomposing the middle term.}
The lift $s_h$ is holomorphic exactly when
$\overline\nabla_\beta \, s_h=(1\otimes s_h)\, \overline\nabla_{(b)}$, or,
equivalently, when
$\beta+D_{a,b}h=0$. It follows that
\[
[\beta]=0
\quad\Longleftrightarrow\quad
\text{the displayed extension splits holomorphically}.
\]
In that case, $\iota(\mathcal L_a)$ and $s_h(\mathcal L_b)$ are
complementary holomorphic submodules, and $g_h$ transforms the
connection to the diagonal one. Notice that a nonzero off-diagonal
map can still represent the zero class: it is the class, rather
than the chosen representative, that obstructs a holomorphic splitting.

A nonzero class, however, does not imply that the middle holomorphic
bundle is indecomposable. It rules out a splitting compatible with
the specified inclusion of $\mathcal L_a$ and projection onto
$\mathcal L_b$, but an arbitrary holomorphic isomorphism need not
preserve these maps. The same middle term may admit a decomposition
as $\mathcal L_r\oplus\mathcal L_s$ with different summands.
The quantum Euler sequence in Section~\ref{subsec:explicit-quantum-Euler} illustrates this:
its extension class is nonzero although its middle term is
$\mathcal A_q^2=\mathcal L_0\oplus\mathcal L_0$ with the standard holomorphic
structure. The extension class therefore records the bundle
together with its specified subbundle and quotient, not the
isomorphism class of the middle bundle alone.

\subsection{The quantum Euler sequence}
\label{subsec:explicit-quantum-Euler}

The two-term Dolbeault complex of $\mathcal L_k$ is
\[
  0\longrightarrow\mathcal L_k
  \xrightarrow{\,\overline\nabla_{(k)}\,}
  \Omega^{0,1}(\mathcal A_q)\otimes_{\mathcal A_q}\mathcal L_k
  \longrightarrow0,
  \qquad
  \Omega^{0,1}(\mathcal A_q)\otimes_{\mathcal A_q}\mathcal L_k
  \simeq\omega_- \, \mathcal L_{k-2}.
\]
We write $H^{0,1}$ for its degree-one cohomology.
The cohomology dimensions are
\begin{equation}
  \dim H^0(\mathcal L_k)
  =\begin{cases}1-k,&k\leq0,\\0,&k>0,\end{cases}
  \qquad
  \dim H^{0,1}(\mathcal L_k)
  =\begin{cases}0,&k\leq1,\\k-1,&k\geq2.\end{cases}
  \label{eq:line-cohomology-kronecker}
\end{equation}
For $k\geq2$, the cokernel of $\overline\nabla_{(k)}$ is the lowest irreducible summand missed by
it. As a $\mathcal U_q(\mathfrak{su}(2))$-module,
\[
  H^{0,1}(\mathcal L_k)\simeq V_{\frac{k}{2}-1}.
\]
Representatives are
\[
  \bigl[\omega_-(c^*)^\mu(a^*)^{k-2-\mu}\bigr],
  \qquad \mu=0,\ldots,k-2,
\]
up to nonzero powers of $q$ arising from ordering and normalization.
These formulas follow from the Peter--Weyl decomposition of
$\mathcal A(\mathrm{SU}_q(2))$ and the explicit action of $X_-$ in
\eqref{eq:standard-KLvS-delbar}; see \cite{KLvS,KhalkhaliMoatadelro}. In particular, nonzero
extensions of $\mathcal L_b$ by $\mathcal L_a$ occur for $a-b\geq2$,
and their parameter space has dimension $a-b-1$.

The smallest nontrivial extension space occurs for $a=1$ and $b=-1$,
when $H^{0,1}(\mathcal L_2)\simeq\mathbb C$.
A nonzero class is represented, up to normalization, by
$\vartheta_0=\omega_-$. The corresponding map
\[
  \beta_0:\mathcal L_{-1}\longrightarrow
  \Omega^{0,1}(\mathcal A_q)\otimes_{\mathcal A_q}\mathcal L_1
\]
defines the upper-triangular connection on
$\mathcal L_1\oplus\mathcal L_{-1}$.
This extension has an explicit realization with trivial middle term:
\begin{equation}
  0\longrightarrow\mathcal L_1
  \xrightarrow{\,j\,}\mathcal A_q^2
  \xrightarrow{\,\pi\,}\mathcal L_{-1}\longrightarrow0,
  \label{eq:quantum-Euler-Kronecker}
\end{equation}
where
\[
  \pi(f,g)=fa+gc,
  \qquad j(\eta)=(-q\eta c,\eta a).
\]
Since $ac=qca$ one has $\pi(j(\eta))=-q\eta ca+\eta ac=0$.
To verify exactness, use the fundamental unitary
\[
  U=\begin{pmatrix}a&-qc^*\\c&a^*\end{pmatrix}.
\]
Right multiplication of row vectors by $U$ gives a module isomorphism
\[
  \mathcal A_q^2\longrightarrow\mathcal L_{-1}\oplus\mathcal L_1,
  \qquad (f,g)\longmapsto(fa+gc,-qfc^*+ga^*).
\]
It sends $j(\eta)$ to $(0,\eta)$, since
\[
  -q(-q\eta c)c^*+(\eta a)a^*
  =\eta(q^2cc^*+aa^*)=\eta.
\]
Hence $j(\mathcal L_1)=\ker\pi$.
The maps $j$ and $\pi$ are holomorphic: $a,c$ are holomorphic sections of
$\mathcal L_{-1}$, and multiplication is compatible with the standard
tensor-product connection. Thus \eqref{eq:quantum-Euler-Kronecker} is the
quantum counterpart of the classical Euler sequence in \eqref{eq-Euler-Kronecker}.
It does not split holomorphically, because the components of a splitting
would lie in
\[
  \operatorname{Hom}_{\mathrm{hol}}(\mathcal L_{-1},\mathcal A_q)
  \simeq H^0(\mathcal L_1)=0.
\]
Its middle term is nevertheless
$\mathcal A_q^2=\mathcal L_0\oplus\mathcal L_0$.
The extension class therefore detects the chosen filtration, not
indecomposability of the bundle.

\subsection{A three-parameter rank-two family}
\label{subsec:three-parameter-family-construction}

For $a=2$ and $b=-2$, \eqref{eq:line-cohomology-kronecker} gives the
three-dimensional extension space
$H^{0,1}(\mathcal L_4)\simeq\mathbb C^3$. Given
$\lambda=(\lambda_0,\lambda_1,\lambda_2)\in\mathbb C^3$, choose the
representative
\begin{equation}
  \vartheta_\lambda
  =\omega_-\bigl(\lambda_0(a^*)^2+
       \lambda_1c^*a^*+\lambda_2(c^*)^2\bigr)
  \label{eq:theta-lambda-kronecker}
\end{equation}
and let
\[
  \beta_\lambda:\mathcal L_{-2}\longrightarrow
  \Omega^{0,1}(\mathcal A_q)\otimes_{\mathcal A_q}\mathcal L_2
\]
be the associated $\mathcal A_q$-linear map. Set
\begin{equation}
  \mathcal E_\lambda
  =\left(\mathcal L_2\oplus\mathcal L_{-2},
    \begin{pmatrix}\overline\nabla_{(2)}&\beta_\lambda\\
                    0&\overline\nabla_{(-2)}\end{pmatrix}\right).
  \label{eq:E-lambda-kronecker}
\end{equation}
It fits into the holomorphic exact sequence
\begin{equation}
  0\longrightarrow\mathcal L_2\longrightarrow\mathcal E_\lambda
  \longrightarrow\mathcal L_{-2}\longrightarrow0.
  \label{eq:E-lambda-extension-kronecker}
\end{equation}
For $\lambda\neq0$, this sequence is non-split. Modulo scalar
automorphisms of the fixed subbundle and quotient, these nonzero extensions
are parametrized by $\mathbb P H^{0,1}(\mathcal L_4)\simeq\mathbb P^2$.
This parameter space does not classify the middle terms up to arbitrary
holomorphic isomorphism. A standard splitting
$\mathcal E_\lambda\simeq\mathcal L_r\oplus\mathcal L_s$ necessarily
satisfies $r+s=0$ by the winding number. Section~\ref{sec:rank-two-kronecker}
determines the possible pairs and their dependence on the extension class.

\subsection{The classical limit}

At $q=1$, the family becomes
\[
  0\longrightarrow\mathcal O(-2)\longrightarrow E_\lambda
  \longrightarrow\mathcal O(2)\longrightarrow0,
  \qquad [\lambda]\in\mathbb P H^1(\mathcal O(-4))\simeq\mathbb P^2.
\]
The middle term splits by the classical Birkhoff--Grothendieck theorem,
but the splitting type depends on the class. Twisting by $\mathcal O(-1)$
gives
\[
  0\longrightarrow\mathcal O(-3)\longrightarrow E_\lambda(-1)
  \longrightarrow\mathcal O(1)\longrightarrow0.
\]
The connecting map
$\delta_\lambda:H^0(\mathcal O(1))\to H^1(\mathcal O(-3))$
is a map between two-dimensional spaces. In suitable bases it has matrix
\[
  C(\lambda)=\begin{pmatrix}\lambda_0&\lambda_1\\
                            \lambda_1&\lambda_2\end{pmatrix}.
\]
Since $H^0(E_\lambda(-1))\simeq\ker C(\lambda)$, the splitting types are
\[
  E_\lambda\simeq
  \begin{cases}
    \mathcal O(-2)\oplus\mathcal O(2),&\lambda=0,\\
    \mathcal O(-1)\oplus\mathcal O(1),&\lambda\neq0,\ \det C(\lambda)=0,\\
    \mathcal O\oplus\mathcal O,& \det C(\lambda)\neq0.
  \end{cases}
\]
Up to normalization of the middle coordinate, the conic
$\lambda_0\lambda_2-\lambda_1^2=0$ separates the two strata of nonzero
extension classes. Thus a projective plane of filtered extensions has only
two isomorphism classes of middle terms.

\subsection{The quantum connecting map}

\label{subsec:quantum-connecting-map-test}

Tensoring \eqref{eq:E-lambda-extension-kronecker} on the left by
$\mathcal L_1$, with its KLvS bimodule connection, gives
\[
  0\longrightarrow\mathcal L_3
  \longrightarrow\mathcal L_1\otimes_{\mathcal A_q}\mathcal E_\lambda
  \longrightarrow\mathcal L_{-1}\longrightarrow0.
\]
The connecting map
\begin{equation}
  \delta_{\lambda,q}:H^0(\mathcal L_{-1})
               \longrightarrow H^{0,1}(\mathcal L_3)
               \label{qcm}
\end{equation}
acts between two-dimensional spaces, with representatives
\[
  H^0(\mathcal L_{-1})=\operatorname{span}_{\mathbb C}\{a,c\},
  \qquad
  H^{0,1}(\mathcal L_3)
    =\operatorname{span}_{\mathbb C}
        \{[\omega_-a^*],[\omega_-c^*]\}.
\]
It is the cup product, or contraction, of a holomorphic section with
$[\vartheta_\lambda]$, determined by the corresponding
$q$-Clebsch--Gordan map. In compatible weight bases its matrix is linear in
$\lambda$:
\[
  C_q(\lambda)=\lambda_0C_{0,q}+\lambda_1C_{1,q}+\lambda_2C_{2,q},
  \qquad
  \lim_{q\to1}C_q(\lambda)
       =\begin{pmatrix}\lambda_0&\lambda_1\\\lambda_1&\lambda_2\end{pmatrix}.
\]
The numerical matrices $C_{j,q}$ depend on the KLvS tensor multiplication rules and
the chosen normalization. Weight considerations suggest, up to row and column
rescalings and reversals, the pattern
\[
  C_q(\lambda)=
  \begin{pmatrix}\alpha_q\lambda_0&\beta_q\lambda_1\\
                  \gamma_q\lambda_1&\delta_q\lambda_2\end{pmatrix}
\]
with nonzero constants $\alpha_q,\beta_q,\gamma_q,\delta_q$.
This pattern is not used as an established formula below: determining its
coefficients requires an explicit calculation. The determinant
$D_q(\lambda):=\det C_q(\lambda)$ has classical limit
$\lambda_0\lambda_2-\lambda_1^2$ in compatible bases.
Section~\ref{subsec:quantum-connecting-conic} obtains the splitting
stratification without assuming the proposed matrix pattern.

More generally, left tensoring by $\mathcal L_t$ gives
\[
  \delta_{\lambda,t}:H^0(\mathcal L_{-2+t})
                    \longrightarrow H^{0,1}(\mathcal L_{2+t}),
\]
and the long exact sequence yields
\[
  \dim H^0(\mathcal L_t\otimes_{\mathcal A_q}\mathcal E_\lambda)
  =\dim H^0(\mathcal L_{2+t})+\dim\ker\delta_{\lambda,t}.
\]
If $\mathcal E_\lambda\simeq\mathcal L_r\oplus\mathcal L_s$ with standard
structures, these dimensions must equal
$\dim H^0(\mathcal L_{r+t})+\dim H^0(\mathcal L_{s+t})$ for every $t$.
They therefore provide cohomological invariants of the splitting type.

The family $\mathcal E_\lambda$ provides a tractable comparison with the
classical theory: its extension space has dimension three, and the first
connecting map relevant to its splitting type is a $2\times2$ matrix.
Section~\ref{sec:rank-two-kronecker} proves decomposability through the
Kronecker model and
identifies the strata through $\det C_q(\lambda)$. The precise coefficients
of this determinant remain a separate computational question.

\subsection{Higher-rank upper-triangular examples}
\label{subsec:higher-rank-construction}

For $\mathcal E=\bigoplus_{i=1}^r\mathcal L_{n_i}$, set
\[
  \overline\nabla_{\mathcal E}
  =\operatorname{diag}(\overline\nabla_{(n_1)},\ldots,\overline\nabla_{(n_r)})+B,
  \qquad B=(\beta_{ij}),
\]
where
\[
  \beta_{ij}\in\operatorname{Hom}_{\mathcal A_q}
  \bigl(\mathcal L_{n_j},\Omega^{0,1}(\mathcal A_q)
                   \otimes_{\mathcal A_q}\mathcal L_{n_i}\bigr)
  \simeq\Omega^{0,1}(\mathcal A_q)
                   \otimes_{\mathcal A_q}\mathcal L_{n_i-n_j}.
\]
Flatness is automatic. When $B$ is strictly upper triangular, the
connection preserves a filtration with standard line quotients.
For example, on $\mathcal L_2\oplus\mathcal L_0\oplus\mathcal L_{-2}$ take
\begin{equation}
  \overline\nabla=
  \begin{pmatrix}
    \overline\nabla_{(2)}&\beta_{12}&\beta_{13}\\
    0&\overline\nabla_{(0)}&\beta_{23}\\
    0&0&\overline\nabla_{(-2)}
  \end{pmatrix},
  \label{eq:rank-three-upper-triangular-recalled}
\end{equation}
with
\begin{equation}
  [\beta_{12}],[\beta_{23}]\in H^{0,1}(\mathcal L_2)\simeq\mathbb C,
  \qquad [\beta_{13}]\in H^{0,1}(\mathcal L_4)\simeq\mathbb C^3.
  \label{eq:rank-three-parameters}
\end{equation}
The classification of its middle terms under general holomorphic
isomorphisms is given in Section~\ref{subsec:rank-three-family-splits}.

\subsection{The classical Clifford construction on the two-sphere}
\label{subsec:classical-Clifford-projector}

Let $\Gamma_1,\Gamma_2,\Gamma_3\in M_N(\mathbb C)$ be Hermitian matrices
satisfying $\Gamma_i\Gamma_j+\Gamma_j\Gamma_i=2\delta_{ij}I_N$.
For $x\in S^2$, put
\[
  Q_\Gamma(x)=x_1\Gamma_1+x_2\Gamma_2+x_3\Gamma_3,
  \qquad p_\Gamma(x)=\tfrac12(I_N+Q_\Gamma(x)).
\]
Since $Q_\Gamma(x)^2=I_N$, $p_\Gamma$ is a smooth projection and defines
$E_\Gamma=p_\Gamma C^\infty(S^2)^N$. Each irreducible $2\times2$ block
contributes rank one: for $N=2r$, the bundle has rank $r$; in particular,
rank $2^m$ requires matrices of size $2^{m+1}$.
The Grassmannian Dolbeault operator
$\overline\nabla_\Gamma=p_\Gamma\circ\bar\partial$ has curvature
\[
  \overline\nabla_\Gamma^{\,2}
  =p_\Gamma(\bar\partial p_\Gamma)\wedge(\bar\partial p_\Gamma)=0,
\]
because $(0,2)$-forms vanish on the Riemann sphere.

The algebra decomposition
\begin{equation}
  \operatorname{Cl}_3(\mathbb C)
  \simeq M_2(\mathbb C)\oplus M_2(\mathbb C)
  \label{eq:Cl3-decomposition-recalled}
\end{equation}
shows that every such representation is a direct sum of the two irreducible
Pauli representations. A constant unitary change of basis therefore reduces
$p_\Gamma$ to Pauli blocks and gives
\[
  E_\Gamma\simeq
  \mathcal O(-1)^{\oplus r_+}\oplus\mathcal O(1)^{\oplus r_-},
\]
where $r_\pm$ are the two multiplicities. Reversing the orientation convention
interchanges the signs. In particular, $\Gamma_i=\sigma_i\otimes I_r$
gives $r$ copies of one Hopf line bundle. The projector thus presents a
higher-rank bundle whose splitting is already determined by the Clifford
representation.

\subsection{Combining projectors with extension classes}
\label{subsec:projectors-and-extensions}

On the quantum projective line, with
$\mathcal A_q=\mathcal A(\mathbb{CP}_q^1)$,
set $v=(a,c)^{\mathrm t}$. Since $v^*v=a^*a+c^*c=1$, the matrix
\[
  p_q=vv^*
  =\begin{pmatrix}aa^*&ac^*\\ca^*&cc^*\end{pmatrix}
  =\begin{pmatrix}1-q^2B_0&B_-\\B_+&B_0\end{pmatrix}
\]
is a projection (the Bott projection), with $B_0,B_\pm$ as in \eqref{eq:podles-generators}.
Equivalently, $Q_q=2p_q-I_2$ is a self-adjoint involution.
In the left-module convention, $p_q$ acts on row vectors by right
multiplication. Its image is the left module $\mathcal A_q^{1\times2}p_q$,
and the mutually inverse maps
\[
 \mathcal A_q^{1\times2}p_q\longrightarrow\mathcal L_{-1},
 \quad r\longmapsto rv,
 \qquad
 \mathcal L_{-1}\longrightarrow\mathcal A_q^{1\times2}p_q,
 \quad \eta\longmapsto\eta v^*,
\]
identify it with $\mathcal L_{-1}$. Indeed, $v^*v=1$ and $rp_q=r$ on
the image. The complementary projection represents $\mathcal L_1$, as
follows from the second column of the fundamental unitary $U$ in
Section~\ref{subsec:explicit-quantum-Euler}. We equip these images with the
KLvS holomorphic structures transported through these identifications;
see also \cite{MajidSphere,KLvS}.

Amplifying by $I_r$ gives $r$ copies of the same line module. More generally,
a bundle associated to
$\mathcal A_q\subset\mathcal A(\mathrm{SU}_q(2))$
through a finite-dimensional representation of $\mathrm U(1)$ decomposes
according to its characters. For the canonical covariant holomorphic
structures, this is a direct sum of standard line objects. These
constructions therefore describe projective modules but do not themselves
produce indecomposable higher-rank holomorphic bundles.

Let $p_a,p_b$ represent $\mathcal L_a,\mathcal L_b$ and set
$P=p_a\oplus p_b$. Under
$\operatorname{im}P\simeq\mathcal L_a\oplus\mathcal L_b$, transport the
diagonal KLvS connection and add an off-diagonal term:
\[
  \overline\nabla_{P, \, \beta}
  =\begin{pmatrix}\overline\nabla_{(a)}&\beta\\0&\overline\nabla_{(b)}\end{pmatrix}.
\]
For $a=2$, $b=-2$, and $\beta=\beta_\lambda$, this realizes
\eqref{eq:E-lambda-kronecker} in projector form. The projector records the
underlying module, whereas the extension class specifies its holomorphic
structure. No new flatness condition occurs.

\section{The rank-two extension family and the Kronecker model}
\label{sec:rank-two-kronecker}

We determine the middle terms of the family
\eqref{eq:E-lambda-kronecker}. Every member is a direct sum of standard KLvS
line objects, and the two nonzero splitting types are separated by the
determinant of a quantum connecting map. The proof uses the Dolbeault
differential graded category, the quantum Euler sequence, quantum Serre
duality and the representation theory of the two-arrow Kronecker quiver.

\subsection{Splitting of the three-parameter family}
\label{subsec:three-parameter-family-recalled}

Let $\mathcal E_\lambda$ be the polynomial holomorphic bundle defined in
\eqref{eq:E-lambda-kronecker}. Its extension class is specified by
\eqref{eq:theta-lambda-kronecker}, with the ordering and normalization
conventions fixed there.

\begin{theorem}[Splitting of the three-parameter family]
\label{thm:three-parameter-splitting}
For every $\lambda\in\mathbb C^3$, the holomorphic bundle
$\mathcal E_\lambda$ is a direct sum of standard KLvS line bundles.  More
precisely,
\begin{equation}
  \mathcal E_0\simeq\mathcal L_2\oplus\mathcal L_{-2},
  \label{eq:E-zero-splitting}
\end{equation}
and, for $\lambda\neq0$,
\begin{equation}
  \mathcal E_\lambda
  \simeq
  \mathcal L_0^{\oplus2}
  \qquad\text{or}\qquad
  \mathcal L_1\oplus\mathcal L_{-1}.
  \label{eq:E-nonzero-two-types}
\end{equation}
In particular, the family \eqref{eq:E-lambda-kronecker} contains no
indecomposable rank-two holomorphic bundle.
\end{theorem}

This section is devoted to the proof and to the
identification of the two nonzero splitting strata.

\subsection{The Dolbeault differential graded category}
\label{subsec:dolbeault-dg-category}

Let $\mathscr D_{\bar\partial}(\mathcal A_q)$ denote the differential graded
category whose objects are finitely generated projective left
$\mathcal A_q$-modules equipped with flat $\bar\partial$-connections.  For two
objects $(\mathcal E,\overline\nabla_{\mathcal E})$ and
$(\mathcal F,\overline\nabla_{\mathcal F})$, put
\begin{equation}
  \underline{\operatorname{Hom}}^{\,j}_{\bar\partial}
  (\mathcal E,\mathcal F)
  :=
  \operatorname{Hom}_{\mathcal A_q}
  \bigl(
     \mathcal E,
     \Omega^{0,j}(\mathcal A_q)
       \otimes_{\mathcal A_q}\mathcal F
  \bigr),
  \qquad j=0,1,
  \label{eq:dolbeault-hom-complex}
\end{equation}
and set this space equal to zero in all other degrees.  If
$T\in\underline{\operatorname{Hom}}^{\,j}_{\bar\partial}
(\mathcal E,\mathcal F)$, define
\begin{equation}
  D T
  =
  \overline\nabla_{\mathcal F}\circ T
  -(-1)^j(1\wedge T)\circ\overline\nabla_{\mathcal E},
  \label{eq:dolbeault-hom-differential}
\end{equation}
where $\overline\nabla_{\mathcal F}$ is extended to
$\Omega^{0,\bullet}\otimes\mathcal F$ by the graded Leibniz rule.  Flatness
implies $D^2=0$. Composition uses composition of module maps and the
product of forms.

The degree-zero cohomology is precisely the space of holomorphic maps:
\begin{equation}
  H^0\underline{\operatorname{Hom}}_{\bar\partial}^{\,\bullet}
  (\mathcal E,\mathcal F)
  =
  \operatorname{Hom}_{\mathrm{hol}}(\mathcal E,\mathcal F).
  \label{eq:H0-holomorphic-maps}
\end{equation}
For the line objects $\mathcal L_n$, the invertibility of the line bimodules and
the KLvS twisted flips give
\begin{equation}
  H^j\underline{\operatorname{Hom}}_{\bar\partial}^{\,\bullet}
  (\mathcal L_m,\mathcal L_n)
  \simeq
  H^{0,j}(\mathcal L_{n-m}),
  \qquad j=0,1.
  \label{eq:line-Ext-Dolbeault}
\end{equation}

Let $\operatorname{Perf}_{\bar\partial}(\mathcal A_q)$ denote the
triangulated homotopy category of the idempotent-complete pretriangulated
hull of $\mathscr D_{\bar\partial}(\mathcal A_q)$. Thus we first adjoin
shifts and cones, pass to the degree-zero homotopy category, and split
idempotents. The original KLvS holomorphic bundles embed fully faithfully
in this category: their morphism complexes have no negative-degree terms,
so their degree-zero cohomology consists of actual holomorphic maps.
Since the underlying modules in
\eqref{eq:E-lambda-extension-kronecker} are projective, that exact sequence
is split as a sequence of $\mathcal A_q$-modules and determines a
distinguished triangle in
$\operatorname{Perf}_{\bar\partial}(\mathcal A_q)$.

\subsection{A quantum Beilinson object and the Kronecker algebra}
\label{subsec:quantum-beilinson-kronecker}

Consider the object
\begin{equation}
  \mathcal T:=\mathcal L_0\oplus\mathcal L_1.
  \label{eq:quantum-beilinson-object}
\end{equation}
With the cohomology dimensions  given in
\eqref{eq:line-cohomology-kronecker},
it follows from \eqref{eq:line-Ext-Dolbeault} that
\begin{equation}
  H^j
  \underline{\operatorname{End}}_{\bar\partial}^{\,\bullet}
  (\mathcal T)=0,
  \qquad j\neq0.
  \label{eq:T-no-higher-self-Ext}
\end{equation}
Moreover,
\begin{equation}
  B
  :={\operatorname{End}}_{\mathrm{hol}}(\mathcal T)
  \simeq
  \begin{pmatrix}
    \mathbb C&H^0(\mathcal L_{-1})\\
    H^0(\mathcal L_1)&\mathbb C
  \end{pmatrix}
  \simeq
  \begin{pmatrix}
    \mathbb C&\mathbb C^2\\
    0&\mathbb C
  \end{pmatrix}.
  \label{eq:T-endomorphism-algebra}
\end{equation}
Precomposition makes $\operatorname{Hom}(\mathcal T,-)$ a right
$B$-module, or equivalently a left module over the opposite algebra
$K=B^{\mathrm{op}}$. With this convention, $K$ is the path algebra of the
two-arrow Kronecker quiver
\begin{equation}
  Q:\qquad
  0\ \substack{\xrightarrow{\ \alpha\ }\\[-1mm]
                 \xrightarrow[\ \beta\ ]{}}\ 1,
  \qquad
  K:=B^{\mathrm{op}}\simeq\mathbb C Q.
  \label{eq:Kronecker-quiver}
\end{equation}
To justify the reduction to this ordinary algebra, put
$\mathscr B^\bullet=
\underline{\operatorname{End}}_{\bar\partial}^{\,\bullet}(\mathcal T)$.
This dg algebra is concentrated in nonnegative degrees, and
$B=\ker(D:\mathscr B^0\to\mathscr B^1)$ is a subalgebra of degree-zero
cycles. The inclusion $B\hookrightarrow\mathscr B^\bullet$ is a
quasi-isomorphism: it induces the identity on $H^0$, and
\eqref{eq:T-no-higher-self-Ext} gives the vanishing in all other degrees.
The dg Morita theorem therefore yields an equivalence
\begin{equation}
  \mathbf R\!\operatorname{Hom}(\mathcal T,-):
  \operatorname{thick}(\mathcal T)
  \xrightarrow{\ \simeq\ }
  \operatorname{Perf}(K).
  \label{eq:dg-Morita-Kronecker}
\end{equation}
Here $\operatorname{thick}(\mathcal T)$ is the smallest full triangulated
subcategory of $\operatorname{Perf}_{\bar\partial}(\mathcal A_q)$ that
contains $\mathcal T$ and is closed under direct summands, and
$\operatorname{Perf}(K)$ is the category of perfect complexes of left
$K$-modules.
Think of $\mathcal T$ as a basic building block, and $B$ as the algebra of its endomorphisms. The functor $\operatorname{RHom}(\mathcal T,-)$ sends $\mathcal T$ to the regular right $B$-module $B$, equivalently the regular left $K$-module, preserving finite sums, shifts, cones and direct summands. These operations generate $\operatorname{thick}(\mathcal T)$ on one side and $\operatorname{Perf}(K)$ on the other. The vanishing of self-extensions in nonzero degrees ensures that the ordinary algebra $B$ retains all the required information. Conversely, one reconstructs an object by replacing copies of $B$ with $\mathcal T$, matrix entries with the corresponding endomorphisms, and projective summands with the corresponding direct summands; see \cite{KellerDG}. The two arrows act by precomposition
with the holomorphic maps $R_a,R_c:\mathcal L_1\to\mathcal L_0$.
The term \emph{Beilinson object} refers here to this two-line object and
its Kronecker endomorphism algebra; generation of the entire KLvS category
is not assumed.

The quantum Euler sequence \eqref{eq:quantum-Euler-Kronecker}
provides the generation argument.
Tensoring on the left by the line bimodule $\mathcal L_n$, with its KLvS
bimodule connection, gives
\begin{equation}
  0\longrightarrow\mathcal L_{n+1}
  \longrightarrow\mathcal L_n^{\oplus2}
  \longrightarrow\mathcal L_{n-1}
  \longrightarrow0.
  \label{eq:shifted-quantum-Euler}
\end{equation}
Induction in both directions shows that every standard line object belongs
to $\operatorname{thick}(\mathcal T)$.  Hence
\begin{equation}
  \operatorname{thick}(\mathcal T)
  =
  \operatorname{thick}
  \{\mathcal L_n:n\in\mathbb Z\}.
  \label{eq:standard-line-generated-category}
\end{equation}
In particular, the extension object $\mathcal E_\lambda$ belongs to this
subcategory.

\subsection{The Kronecker representation attached to
\texorpdfstring{$\mathcal E_\lambda$}{E(lambda)}}
\label{subsec:E-lambda-Kronecker-representation}

Assume for the moment that $\lambda\neq0$.  Applying Dolbeault cohomology to
\eqref{eq:E-lambda-extension-kronecker} gives
\begin{equation}
  0\longrightarrow H^0(\mathcal E_\lambda)
  \longrightarrow H^0(\mathcal L_{-2})
  \xrightarrow{\ \delta_\lambda\ }
  H^{0,1}(\mathcal L_2)
  \longrightarrow H^{0,1}(\mathcal E_\lambda)
  \longrightarrow0.
  \label{eq:E-lambda-long-exact-L0}
\end{equation}
Here
\[
  \dim H^0(\mathcal L_{-2})=3,
  \qquad
  \dim H^{0,1}(\mathcal L_2)=1.
\]
The connecting map is cup product with the extension class:
\begin{equation}
  \delta_\lambda(s)
  =s\smile[\vartheta_\lambda].
  \label{eq:delta-lambda-cup-product}
\end{equation}
Quantum Serre duality identifies the pairing
\begin{equation}
  H^0(\mathcal L_{-2})
  \otimes H^{0,1}(\mathcal L_4)
  \longrightarrow H^{0,1}(\mathcal L_2)
  \simeq\mathbb C
  \label{eq:quantum-Serre-pairing}
\end{equation}
with a perfect pairing; see \cite{KhalkhaliMoatadelro}.  We can check
nondegeneracy directly in the factor order used here.  Let $h$ be the
normalized Haar state of $\mathcal A(\mathrm{SU}_q(2))$.  The functional
\[
  [\omega_-f]\longmapsto h(f),\qquad f\in\mathcal A_q,
\]
identifies $H^{0,1}(\mathcal L_2)$ with $\mathbb C$: Haar invariance
annihilates the image of $\overline\nabla_{(2)}$, and $h(1)=1$.
Represent $[\vartheta_\lambda]$ by $\omega_-u$, where
\[
  u=\lambda_0(a^*)^2+\lambda_1c^*a^*+\lambda_2(c^*)^2.
\]
The KLvS twisted flip is induced by multiplication in the total calculus
\cite[Lemma~3.6]{KLvS}.  Thus the coefficient order in
$\beta_{\vartheta_\lambda}(s)$ is $su$, for
$s\in H^0(\mathcal L_{-2})$.  More explicitly,
$\omega_-s=q^{-2}s\omega_-$ gives
$s\omega_-u=q^2\omega_-su$ in the multiplication identifications.
Consequently the pairing is a fixed nonzero scalar multiple of $h(su)$;
rescaling the one-dimensional target changes only that scalar.
If $u\neq0$, then $s=u^*$ is a holomorphic section of
$\mathcal L_{-2}$ and $h(u^*u)>0$ by faithfulness of the Haar state.
The pairing is therefore nondegenerate in its second factor, and the two
factors both have dimension three, so it is perfect.  In particular,
$[\vartheta_\lambda]\neq0$ implies that $\delta_\lambda$ is a nonzero
functional, hence is surjective.  It follows that
\begin{equation}
  \dim H^0(\mathcal E_\lambda)=2,
  \qquad
  H^{0,1}(\mathcal E_\lambda)=0.
  \label{eq:E-lambda-H0-H1}
\end{equation}
Next apply
$\mathbf R\!\operatorname{Hom}(\mathcal L_1,-)$ to
\eqref{eq:E-lambda-extension-kronecker}.  Formula
\eqref{eq:line-Ext-Dolbeault} gives
\begin{align*}
  \operatorname{Hom}_{\mathrm{hol}}(\mathcal L_1,\mathcal L_2)
  &\simeq H^0(\mathcal L_1)=0,\\
  \operatorname{Hom}_{\mathrm{hol}}(\mathcal L_1,\mathcal L_{-2})
  &\simeq H^0(\mathcal L_{-3})\simeq\mathbb C^4,\\
  \operatorname{Ext}^1_{\bar\partial}(\mathcal L_1,\mathcal L_2)
  &\simeq H^{0,1}(\mathcal L_1)=0,\\
  \operatorname{Ext}^1_{\bar\partial}(\mathcal L_1,\mathcal L_{-2})
  &\simeq H^{0,1}(\mathcal L_{-3})=0.
\end{align*}
Therefore
\begin{equation}
  \dim\operatorname{Hom}_{\mathrm{hol}}
  (\mathcal L_1,\mathcal E_\lambda)=4,
  \qquad
  \operatorname{Ext}^1_{\bar\partial}
  (\mathcal L_1,\mathcal E_\lambda)=0.
  \label{eq:Hom-L1-E-lambda}
\end{equation}

Combining \eqref{eq:E-lambda-H0-H1} and
\eqref{eq:Hom-L1-E-lambda}, one finds that
\begin{equation}
  M_\lambda
  :=\mathbf R\!\operatorname{Hom}(\mathcal T,\mathcal E_\lambda)
  \label{eq:M-lambda-definition}
\end{equation}
is concentrated in cohomological degree zero.  It is therefore an ordinary
finite-dimensional representation of the Kronecker quiver
\eqref{eq:Kronecker-quiver}, with dimension vector
\begin{equation}
  \underline{\dim}\,M_\lambda=(2,4).
  \label{eq:M-lambda-dimension-vector}
\end{equation}
Thus $M_\lambda$ may be written as a pair of linear maps
\begin{equation}
  \mathbb C^2
  \ \substack{\xrightarrow{\ A_\lambda\ }\\[-1mm]
                \xrightarrow[\ B_\lambda\ ]{}}
  \ \mathbb C^4.
  \label{eq:M-lambda-two-arrows}
\end{equation}

\paragraph{Decomposability of the dimension vector \texorpdfstring{$(2,4)$}{(2,4)}.}

We first record an elementary fact which does not require the classification
of all Kronecker representations.

\begin{lemma}
\label{lem:Kronecker-2-4}
Every representation of the two-arrow Kronecker quiver having
dimension vector $(2,4)$ is decomposable.
\end{lemma}

\begin{proof}
Let
\[
  A,B:V\longrightarrow W,
  \qquad
  \dim V=2,
  \quad
  \dim W=4,
\]
be the two arrows, and form
\begin{equation}
  \Phi:V\oplus V\longrightarrow W,
  \qquad
  \Phi(v_1,v_2)=A(v_1)+B(v_2).
  \label{eq:combined-Kronecker-map}
\end{equation}
If $\Phi$ is not surjective, let $W'=\operatorname{im}\Phi$ and choose a
nonzero complement $W_0$ with $W=W'\oplus W_0$.  Since the images of both
$A$ and $B$ lie in $W'$, one has
\[
  (V,W;A,B)
  \simeq
  (V,W';A,B)\oplus(0,W_0;0,0).
\]

Suppose now that $\Phi$ is surjective.  Since both its domain and codomain
have dimension four, it is an isomorphism.  Choose a decomposition
$V=V_1\oplus V_2$ into one-dimensional subspaces and put
\[
  W_i:=\Phi(V_i\oplus V_i)
  =A(V_i)+B(V_i),
  \qquad i=1,2.
\]
The injectivity of $\Phi$ gives $W=W_1\oplus W_2$, and both arrows carry
$V_i$ into $W_i$.  Hence
\[
  (V,W;A,B)
  \simeq
  (V_1,W_1;A|_{V_1},B|_{V_1})
  \oplus
  (V_2,W_2;A|_{V_2},B|_{V_2}).
\]
Thus the representation is decomposable in either case.
\end{proof}

Applying Lemma~\ref{lem:Kronecker-2-4} to $M_\lambda$ gives a nontrivial idempotent
\[
  e\in\operatorname{End}_{K}(M_\lambda).
\]
The fully faithful equivalence \eqref{eq:dg-Morita-Kronecker} transports it
to an idempotent class
\[
  [p]\in
  H^0\underline{\operatorname{End}}_{\bar\partial}^{\,\bullet}
  (\mathcal E_\lambda).
\]
The endomorphism complex of the original bundle $\mathcal E_\lambda$ has no
negative-degree term.  Hence there are no degree-zero boundaries, and the
identity $[p]^2=[p]$ is represented by an actual holomorphic
$\mathcal A_q$-linear endomorphism satisfying
\begin{equation}
  p^2=p,
  \qquad
  (1\otimes p)\overline\nabla_\lambda
  =\overline\nabla_\lambda p.
  \label{eq:actual-holomorphic-idempotent}
\end{equation}
Therefore
\begin{equation}
  \mathcal E_\lambda
  =\operatorname{im}p\oplus\ker p
  \label{eq:actual-holomorphic-splitting}
\end{equation}
inside the KLvS holomorphic category.  Both summands are finitely generated
projective modules, and their induced $\bar\partial$-connections are flat.
Thus $\mathcal E_\lambda$ is holomorphically decomposable for
all $\lambda\neq0$.

\subsection{Determination of the splitting type}
\label{subsec:determination-splitting-type}

We now identify the two summands in \eqref{eq:actual-holomorphic-splitting}.  Recall the classification of
finite-dimensional indecomposable representations of the two-arrow
Kronecker quiver.  Their dimension vectors are
\begin{equation}
  (n,n+1),
  \qquad
  (n,n),
  \qquad
  (n+1,n),
  \label{eq:Kronecker-indec-dim-vectors}
\end{equation}
where $n\geq0$ in the first and third families and $n\geq1$ in the
middle family. These are respectively the preprojective, regular and
preinjective families; see
\cite{RingelTame}.

For $n\geq0$, formula \eqref{eq:line-cohomology-kronecker} shows that the
image of $\mathcal L_{1-n}$ is concentrated in degree zero and has dimension
vector $(n,n+1)$. It is indecomposable because
$\operatorname{End}_{\mathrm{hol}}(\mathcal L_{1-n})=\mathbb C$. Hence
\begin{equation}
  \mathbf R\!\operatorname{Hom}
  (\mathcal T,\mathcal L_{1-n})
  \simeq P_n,
  \qquad
  \underline{\dim}\,P_n=(n,n+1),
  \label{eq:line-to-preprojective}
\end{equation}
where $P_n$ is the unique indecomposable preprojective representation of
that dimension vector.  The cases $n=0,1,2$ read
\begin{equation}
  \mathcal L_1\longleftrightarrow(0,1),
  \qquad
  \mathcal L_0\longleftrightarrow(1,2),
  \qquad
  \mathcal L_{-1}\longleftrightarrow(2,3).
  \label{eq:first-three-line-images}
\end{equation}

There is also a useful rank interpretation.  If a finite-dimensional
left module over the Kronecker algebra $K$ of
\eqref{eq:Kronecker-quiver} has dimension vector $(a,b)$, then the virtual
generic rank of its
inverse image under \eqref{eq:dg-Morita-Kronecker} is
\begin{equation}
  \operatorname{rk}(a,b)=b-a.
  \label{eq:rank-dimension-vector}
\end{equation}
Indeed, both sides are additive on $K_0(K)$.  The simple module at vertex
$1$ is the image of $\mathcal L_1$ and has rank $1$, while the simple module
at vertex $0$ is the image of $\mathcal L_2[1]$ and therefore has virtual
rank $-1$.

For completeness, the generic rank used here is well defined in the
polynomial category. Set $t=B_0$, $x=B_-$ and $y=B_+$. The relations
\eqref{eq:podles-relations-one}--\eqref{eq:podles-relations-two} identify
$\mathcal A_q$ with the generalized Weyl algebra over $\mathbb C[t]$
determined by the automorphism $\sigma$ of $\mathbb C[t]$
and the polynomial $f_0$, where
\[
 \sigma(t)=q^2t,\qquad f_0(t)=t(1-t),\qquad
 yx=f_0(t),\quad xy=\sigma(f_0(t)).
\]
The structure theorem for generalized Weyl algebras shows that
this is a left and right noetherian domain, hence an Ore domain
\cite{BavulaGWA}; see also \cite[Section~5.8]{KaygunGWA} for this
presentation of the Podle\'s sphere. Explicitly, it embeds in the skew
Laurent ring $\mathbb C(t)[z,z^{-1};\sigma]$ by
\[
 t\longmapsto t,\qquad x\longmapsto z,\qquad
 y\longmapsto f_0(t)z^{-1}.
\]
Let $\mathcal Q_q$ be its division ring of fractions. For a finitely
generated projective left module $P$, set
\[
 \operatorname{rk}P
 :=\dim_{\mathcal Q_q}\bigl(\mathcal Q_q\otimes_{\mathcal A_q}P\bigr).
\]
A nonzero such module has positive integral rank: it embeds as a direct
summand of a finite free module, and localization is injective on that
module. Applied to an actual holomorphic projective summand of
$\mathcal E_\lambda$, this excludes the regular dimension vectors
$(n,n)$, of rank zero, and the preinjective vectors $(n+1,n)$, of virtual
rank $-1$. Only the preprojective vectors $(n,n+1)$ remain.

Since $\mathcal E_\lambda$ has rank two, its indecomposable decomposition
contains exactly two positive-rank summands.  Write their dimension vectors
as
\[
  (n,n+1)
  \qquad\text{and}\qquad
  (m,m+1).
\]
Their sum is $(2,4)$, so
\[
  n+m=2.
\]
Up to order, there are only two possibilities:
\[
  (n,m)=(1,1)
  \qquad\text{or}\qquad
  (n,m)=(0,2).
\]
Using \eqref{eq:line-to-preprojective}, these give respectively
\begin{equation}
  \mathcal E_\lambda\simeq\mathcal L_0^{\oplus2}
  \qquad\text{or}\qquad
  \mathcal E_\lambda\simeq\mathcal L_1\oplus\mathcal L_{-1}.
  \label{eq:identified-two-splittings}
\end{equation}
This completes the proof of Theorem~\ref{thm:three-parameter-splitting}.

\subsection{The connecting map and the splitting conic}
\label{subsec:quantum-connecting-conic}

The two nonzero splitting types are detected by the connecting map obtained
by applying
$\mathbf R\!\operatorname{Hom}(\mathcal L_{-1},-)$ to
\eqref{eq:E-lambda-extension-kronecker}.  It is
\begin{equation}
  \delta_{\lambda,q} :
  H^0(\mathcal L_{-1})
  \longrightarrow
  H^{0,1}(\mathcal L_3).
  \label{eq:quantum-connecting-map-final}
\end{equation}
This is the map \eqref{qcm}, with associated matrix $C_q(\lambda)$.
Both spaces are two-dimensional, and $C_q(\lambda)$ depends linearly on the
extension class $\lambda$.  The long exact sequence gives
\begin{equation}
  \ker C_q(\lambda)
  \simeq
  \operatorname{Hom}_{\mathrm{hol}}
  (\mathcal L_{-1},\mathcal E_\lambda).
  \label{eq:kernel-Cq-Hom}
\end{equation}
The three splitting types now give
\begin{equation}
  \dim\ker C_q(\lambda)
  =
  \begin{cases}
    2,&\lambda=0,\\
    1,&\mathcal E_\lambda\simeq
       \mathcal L_1\oplus\mathcal L_{-1},\\
    0,&\mathcal E_\lambda\simeq
       \mathcal L_0^{\oplus2}.
  \end{cases}
  \label{eq:kernel-Cq-three-cases}
\end{equation}
Consequently,
\begin{equation}
  \mathcal E_\lambda
  \simeq
  \begin{cases}
    \mathcal L_2\oplus\mathcal L_{-2},
      &\lambda=0,\\[1mm]
    \mathcal L_1\oplus\mathcal L_{-1},
      &\lambda\neq0\ \text{and}\ \det C_q(\lambda)=0,\\[1mm]
    \mathcal L_0^{\oplus2},
      &\det C_q(\lambda)\neq0.
  \end{cases}
  \label{eq:final-splitting-stratification}
\end{equation}
In particular,
\begin{equation}
  C_q(\lambda)=0
  \quad\Longleftrightarrow\quad
  \lambda=0.
  \label{eq:Cq-injective}
\end{equation}
Thus the linear map
\[
  H^{0,1}(\mathcal L_4)
  \longrightarrow
  \operatorname{Hom}_{\mathbb C}
  \bigl(H^0(\mathcal L_{-1}),H^{0,1}(\mathcal L_3)\bigr),
  \qquad
  \lambda\longmapsto C_q(\lambda),
\]
is injective.  Its image is a three-dimensional linear subspace of
$M_2(\mathbb C)$.  Such a subspace cannot consist entirely of singular
matrices: projectively, the singular matrices form a smooth quadric surface
in $\mathbb P^3$, which contains no projective plane.  It follows that
$\det C_q(\lambda)$ is a nonzero homogeneous quadratic polynomial.
Therefore
\begin{equation}
  \mathcal C_q
  :=
  \bigl\{[\lambda]\in\mathbb P^2:
          \det C_q(\lambda)=0\bigr\}
  \label{eq:q-splitting-conic}
\end{equation}
is a projective conic, possibly reducible.  The splitting stratification is
\begin{equation}
  [\lambda]\notin\mathcal C_q
  \quad\Longleftrightarrow\quad
  \mathcal E_\lambda\simeq\mathcal L_0^{\oplus2},
  \label{eq:off-conic-trivial}
\end{equation}
whereas
\begin{equation}
  [\lambda]\in\mathcal C_q
  \quad\Longleftrightarrow\quad
  \mathcal E_\lambda\simeq
  \mathcal L_1\oplus\mathcal L_{-1}.
  \label{eq:on-conic-balanced}
\end{equation}
The explicit $q$-Clebsch--Gordan calculation remains useful: it determines
the precise equation of $\mathcal C_q$ and decides whether the conic is
smooth or reducible.  It is no longer needed, however, to decide whether the
middle term is decomposable.

The splitting in Theorem~\ref{thm:three-parameter-splitting} is obtained
by polynomial holomorphic maps. It therefore persists under scalar
extension to a smooth Fr\'echet completion in which these maps and the
KLvS connections extend. This concerns the family
\eqref{eq:E-lambda-kronecker}, whose objects belong to
$\operatorname{thick}(\mathcal T)$; it does not establish generation of
all KLvS holomorphic bundles. We next extend the argument to finite
iterated extensions of standard line objects.

\section{Decomposition of line-filtered bundles}
\label{sec:all-line-filtered-examples-split}

The Kronecker argument extends from the three-parameter family to all
polynomial holomorphic bundles obtained by finite iterated extensions
of the line objects $\mathcal L_n$. Their middle terms split into
direct sums of these objects, although the prescribed extension
sequences need not split.

\subsection{A splitting theorem for line-filtered bundles}
\label{subsec:line-filtered-splitting-theorem}

We use the Dolbeault dg category and the Beilinson object
$\mathcal T=\mathcal L_0\oplus\mathcal L_1$ from
Sections~\ref{subsec:dolbeault-dg-category}
and~\ref{subsec:quantum-beilinson-kronecker}. In particular,
\eqref{eq:dg-Morita-Kronecker} identifies
$\operatorname{thick}(\mathcal T)$ with $\operatorname{Perf}(K)$, where
$K=\operatorname{End}_{\mathrm{hol}}(\mathcal T)^{\mathrm{op}}$ is the
two-arrow Kronecker algebra. By \eqref{eq:standard-line-generated-category},
this subcategory contains all standard line objects and their finite
extensions. The Hom cohomology and line-cohomology formulas are
\eqref{eq:line-Ext-Dolbeault} and \eqref{eq:line-cohomology-kronecker}.

Call a polynomial KLvS holomorphic bundle \(\mathcal E\)
\emph{line-filtered} if it has a finite filtration by holomorphic subbundles
\begin{equation}
  0=\mathcal E_0\subset\mathcal E_1\subset\cdots
  \subset\mathcal E_r=\mathcal E
  \label{eq:line-filtration}
\end{equation}
such that
\begin{equation}
  \mathcal E_i/\mathcal E_{i-1}
  \simeq\mathcal L_{n_i}
  \label{eq:line-filtration-quotients}
\end{equation}
with the standard KLvS holomorphic structure.  Since the quotients are
projective as \(\mathcal A_q\)-modules, the filtration splits at the level of
underlying modules.  The nontrivial information lies in the
upper-triangular \(\bar\partial\)-connection.

\begin{theorem}[Splitting of line-filtered KLvS bundles]
\label{thm:line-filtered-splitting}
Every polynomial line-filtered KLvS holomorphic bundle is
holomorphically isomorphic to a direct sum of standard line objects:
\begin{equation}
  \mathcal E
  \simeq
  \mathcal L_{m_1}\oplus\cdots\oplus\mathcal L_{m_r}.
  \label{eq:line-filtered-final-splitting}
\end{equation}
The integers \(m_1,\ldots,m_r\) are uniquely determined up to permutation.
\end{theorem}

\begin{proof}
Choose an integer \(t\) so negative that
\begin{equation}
  t+n_i\leq1
  \qquad\text{for every }i.
  \label{eq:negative-twist-condition}
\end{equation}
Left tensoring by \(\mathcal L_t\), with its KLvS bimodule
connection and twisted flip, acts on a holomorphic left-module map
\(f:\mathcal E\to\mathcal F\) by \(1\otimes f\).  The tensor-product
connection formula shows that \(1\otimes f\) is holomorphic; no right
linearity of \(f\) is required.  The invertibility of the line bimodule
makes this functor exact, with inverse given by left tensoring by
\(\mathcal L_{-t}\).  Indeed, multiplication identifies
\(\mathcal L_{-t}\otimes_{\mathcal A_q}\mathcal L_t\) and
\(\mathcal L_t\otimes_{\mathcal A_q}\mathcal L_{-t}\) with
\(\mathcal A_q\), carrying the induced tensor-product connections and
composite twisted flips to those of \(\mathcal A_q\)
\cite[Propositions~3.7--3.8]{KLvS}.  The resulting natural inverse
identifications therefore intertwine the induced connections on every
left module.  Since tensoring sends \(\mathcal L_n\) to \(\mathcal L_{t+n}\), this
restricts to an exact autoequivalence of the line-filtered category.  Put
\[
  \mathcal E(t):=\mathcal L_t\otimes_{\mathcal A_q}\mathcal E.
\]
Its filtration factors are \(\mathcal L_{t+n_i}\).  By
\eqref{eq:line-Ext-Dolbeault} and
\eqref{eq:negative-twist-condition},
\begin{align}
  \operatorname{Ext}^1_{\bar\partial}
  (\mathcal L_0,\mathcal L_{t+n_i})
  &\simeq H^{0,1}(\mathcal L_{t+n_i})=0,
  \notag\\
  \operatorname{Ext}^1_{\bar\partial}
  (\mathcal L_1,\mathcal L_{t+n_i})
  &\simeq H^{0,1}(\mathcal L_{t+n_i-1})=0.
  \label{eq:negative-twist-Ext-vanishing}
\end{align}
The long exact sequences associated with the filtration now give
\begin{equation}
  \operatorname{Ext}^1_{\bar\partial}
  (\mathcal T,\mathcal E(t))=0.
  \label{eq:T-Et-Ext-vanishing}
\end{equation}
There are no higher Dolbeault Ext groups.  Hence
\begin{equation}
  M_{\mathcal E,t}
  :=\mathbf R\!\operatorname{Hom}(\mathcal T,\mathcal E(t))
  \label{eq:filtered-Kronecker-module}
\end{equation}
is concentrated in degree zero and is an ordinary finite-dimensional
representation of the Kronecker quiver.

By Section~\ref{subsec:determination-splitting-type}, the indecomposable
Kronecker dimension vectors are as in \eqref{eq:Kronecker-indec-dim-vectors},
and the preprojective module $P_n$ corresponds to $\mathcal L_{1-n}$
as in \eqref{eq:line-to-preprojective}. The generic-rank function is
\eqref{eq:rank-dimension-vector}. Since $\mathcal A_q$ is an Ore domain,
every nonzero projective direct summand has positive generic rank.

Decompose \(M_{\mathcal E,t}\) into indecomposable Kronecker modules. Then every
projection onto an indecomposable summand gives an idempotent in
\[
  \operatorname{End}_K(M_{\mathcal E,t})
  \simeq
  H^0\underline{\operatorname{End}}_{\bar\partial}^{\,\bullet}
  (\mathcal E(t)).
\]
The Dolbeault endomorphism complex of an ordinary holomorphic bundle has no
negative-degree term, and hence has no degree-zero boundaries.  The
idempotent is therefore represented by an actual holomorphic endomorphism
\begin{equation}
  p^2=p,
  \qquad
  (1\otimes p)\overline\nabla_{\mathcal E(t)}
  =\overline\nabla_{\mathcal E(t)}p.
  \label{eq:filtered-actual-idempotent}
\end{equation}
Its image is a nonzero finitely generated projective holomorphic direct
summand and therefore has positive rank.  Formula
\eqref{eq:rank-dimension-vector} rules out the regular vectors \((n,n)\),
which have rank zero, and the preinjective vectors \((n+1,n)\), which have
rank \(-1\).  Thus every indecomposable summand is preprojective and, by
\eqref{eq:line-to-preprojective}, corresponds to a standard line
object.  Consequently \(\mathcal E(t)\) is a direct sum of standard lines.
Tensoring by \(\mathcal L_{-t}\) proves
\eqref{eq:line-filtered-final-splitting}.

For uniqueness, choose a common sufficiently negative twist for any two
proposed standard-line decompositions. Their images are direct sums of the
preprojective modules $P_n$. The Krull--Schmidt theorem for finite-dimensional
Kronecker modules determines their multiplicities, and
\eqref{eq:line-to-preprojective} recovers the line indices after undoing
the twist. This also agrees with the distinct winding numbers of the
line modules \(\mathcal L_n\).
\end{proof}

The line-filtered category, with the exact sequences inherited from the
KLvS category, is not semisimple. The quantum Euler sequence
\eqref{eq:quantum-Euler-Kronecker} is a non-split sequence entirely within
this category, although its middle term is holomorphically isomorphic to
\(\mathcal L_0\oplus\mathcal L_0\).

\subsection{All rank-two upper-triangular extensions}
\label{subsec:all-rank-two-extensions-split}

Consider an arbitrary upper-triangular extension constructed above:
\begin{equation}
  0\longrightarrow\mathcal L_a
  \longrightarrow\mathcal E_\beta
  \longrightarrow\mathcal L_b
  \longrightarrow0,
  \qquad
  [\beta]\in H^{0,1}(\mathcal L_{a-b}).
  \label{eq:general-rank-two-filtered-extension}
\end{equation}
Theorem~\ref{thm:line-filtered-splitting} implies that \(\mathcal E_\beta\) is a
direct sum of two standard lines.  The possible splitting types can be
identified explicitly.

\begin{corollary}[Rank-two balancing]
\label{cor:rank-two-balancing}
Put \(d=a-b\).  If \(d\leq1\) the extension group vanishes and
\[
  \mathcal E_\beta\simeq\mathcal L_a\oplus\mathcal L_b.
\]
If \(d\geq2\), there is a unique integer
\[
  0\leq j\leq\left\lfloor\frac d2\right\rfloor
\]
such that
\begin{equation}
  \mathcal E_\beta
  \simeq
  \mathcal L_{a-j}\oplus\mathcal L_{b+j}.
  \label{eq:rank-two-balancing-formula}
\end{equation}
Moreover, \(j=0\) if and only if \([\beta]=0\).
\end{corollary}

\begin{proof}
Only the nonzero case with \(d\geq2\) requires discussion.  Tensor
\eqref{eq:general-rank-two-filtered-extension} on the left by
\(\mathcal L_{-b}\).  This gives
\begin{equation}
  0\longrightarrow\mathcal L_d
  \longrightarrow\widetilde{\mathcal E}_\beta
  \longrightarrow\mathcal L_0
  \longrightarrow0.
  \label{eq:normalized-rank-two-extension}
\end{equation}
Since \(d>0\), one has \(H^0(\mathcal L_d)=0\).  The beginning of the long
exact sequence is
\begin{equation}
  0\longrightarrow H^0(\widetilde{\mathcal E}_\beta)
  \longrightarrow H^0(\mathcal L_0)
  \xrightarrow{\ \delta_\beta\ }
  H^{0,1}(\mathcal L_d),
  \qquad
  \delta_\beta(1)=[\beta].
  \label{eq:normalized-rank-two-connecting-map}
\end{equation}
Hence \([\beta]\neq0\) implies
\(H^0(\widetilde{\mathcal E}_\beta)=0\).  Write
\[
  \widetilde{\mathcal E}_\beta
  \simeq\mathcal L_u\oplus\mathcal L_v,
  \qquad u\geq v.
\]
The vanishing of global sections forces \(u,v>0\), while the winding number
gives \(u+v=d\).  Therefore
\[
  (u,v)=(d-j,j),
  \qquad
  1\leq j\leq\left\lfloor\frac d2\right\rfloor.
\]
Tensoring back by \(\mathcal L_b\) gives
\eqref{eq:rank-two-balancing-formula}.
\end{proof}

For $a=2$ and $b=-2$, Corollary~\ref{cor:rank-two-balancing} recovers
\eqref{eq:final-splitting-stratification}. The determinant of the connecting
map \eqref{eq:quantum-connecting-map-final} distinguishes the two nonzero
types. Its explicit coefficients determine the equation of the conic;
they are not needed to establish decomposability.

\subsection{The rank-three upper-triangular family}
\label{subsec:rank-three-family-splits}

Consider the rank-three family
\eqref{eq:rank-three-upper-triangular-recalled} with parameters
\eqref{eq:rank-three-parameters}. It is line-filtered and hence splits
into three standard line objects.

Write its splitting type as
\begin{equation}
  \mathcal E
  \simeq
  \mathcal L_{m_1}\oplus\mathcal L_{m_2}\oplus\mathcal L_{m_3},
  \qquad
  m_1\geq m_2\geq m_3.
  \label{eq:rank-three-splitting-indices}
\end{equation}
The winding number gives
\[
  m_1+m_2+m_3=0.
\]
For any integer \(t\), the long exact cohomology sequences of the filtration
give
\begin{equation}
  \dim H^0(\mathcal L_t\otimes\mathcal E)
  \leq
  \sum_{n\in\{2,0,-2\}}
  \dim H^0(\mathcal L_{n+t}).
  \label{eq:rank-three-h0-inequality}
\end{equation}
At \(t=1\), the right-hand side is equal to \(2\), whereas the left-hand
side is
\[
  \sum_{i=1}^3\max\{-m_i,0\}.
\]
Since \(m_1+m_2+m_3=0\), this is also the sum of the positive entries among
the \(m_i\).  It follows that
\[
  m_1\leq2,
  \qquad
  m_1+m_2\leq2.
\]
An elementary enumeration now gives exactly five possible splitting types:
\begin{equation}
  \begin{aligned}
  \mathcal E&\simeq
    \mathcal L_2\oplus\mathcal L_0\oplus\mathcal L_{-2},\\
  \mathcal E&\simeq
    \mathcal L_1^{\oplus2}\oplus\mathcal L_{-2},\\
  \mathcal E&\simeq
    \mathcal L_2\oplus\mathcal L_{-1}^{\oplus2},\\
  \mathcal E&\simeq
    \mathcal L_1\oplus\mathcal L_0\oplus\mathcal L_{-1},\\
  \mathcal E&\simeq
    \mathcal L_0^{\oplus3}.
  \end{aligned}
  \label{eq:five-rank-three-types}
\end{equation}
All five types occur by varying the entries $\beta_{ij}$ in
\eqref{eq:rank-three-upper-triangular-recalled}. If all three entries
vanish, one obtains the first type. If $\beta_{13}=\beta_{23}=0$ and
$[\beta_{12}]\neq0$, the rank-two balancing Corollary~\ref{cor:rank-two-balancing} gives
\[
  \mathcal L_1^{\oplus2}\oplus\mathcal L_{-2}.
\]
If $\beta_{12}=\beta_{13}=0$ and $[\beta_{23}]\neq0$, one obtains
\[
  \mathcal L_2\oplus\mathcal L_{-1}^{\oplus2}.
\]
Finally, if $\beta_{12}=\beta_{23}=0$ and $[\beta_{13}]\neq0$, the object is
\(\mathcal L_0\oplus\mathcal E_{\beta_{13}}\), and the two nonzero
rank-two strata give the final two types in
\eqref{eq:five-rank-three-types}.

The splitting type is already determined by three twisted section dimensions.
Put
\begin{equation}
  h_t(\mathcal E)
  :=\dim H^0(\mathcal L_t\otimes_{\mathcal A_q}\mathcal E).
  \label{eq:rank-three-ht}
\end{equation}
Then
\begin{equation}
\begin{array}{c@{\qquad}ccc}
\text{splitting indices}&h_0&h_1&h_2\\[2pt]
(2,0,-2)&4&2&1\\
(1,1,-2)&3&2&1\\
(2,-1,-1)&4&2&0\\
(1,0,-1)&3&1&0\\
(0,0,0)&3&0&0
\end{array}
  \label{eq:rank-three-section-table}
\end{equation}
Thus these three section dimensions determine the splitting type
completely. They can be computed from the connecting maps in the long
exact sequences of the twisted filtration, including when several
off-diagonal entries are simultaneously nonzero.

\subsection{Scope and the remaining generation problem}
\label{subsec:evidence-and-unfiltered-problem}

The splitting theorem applies to finite extensions of standard KLvS line
objects. The deformation parameter $q$ enters the cup products and the
determinantal equations for their splitting strata. The theorem does not
assert that every KLvS holomorphic bundle admits such a filtration.
Nonstandard holomorphic structures on the canonical line modules
\cite{GravemanEtAl} show why a decomposition into standard line objects
cannot be conjectured for the unrestricted category.

A full matrix $\bar\partial$-connection need not come with a filtration
by standard line objects. Section~\ref{sec:nc-riccati-splitting} studies its
invariant rank-one summands through the holomorphic projector equation
and the associated Riccati equations. Any broader splitting question must
specify whether nonstandard holomorphic line structures are allowed as
summands, as well as the underlying polynomial or completed category.

A separate problem is to identify further explicitly restricted classes
of bundles that admit finite filtrations or suitable resolutions by
standard lines. For a chosen subcategory extending the line-filtered
category, Serre-type vanishing and global-generation statements could
provide tools for such a comparison, but their hypotheses and any
required exactness must be established within that subcategory.
Riemann--Roch and Serre duality for the standard line objects
\cite{KhalkhaliMoatadelro} do not by themselves establish these statements
for arbitrary holomorphic projective modules. We make no unrestricted
standard-line generation assertion.

\medskip

\begin{remark}[The Clifford, Bott, and projector constructions]
\label{subsec:Clifford-Bott-not-counterexamples}
The Clifford and Bott decompositions were described in
Sections~\ref{subsec:classical-Clifford-projector}
and~\ref{subsec:projectors-and-extensions}.
For the projector-plus-extension
construction of Section~\ref{subsec:projectors-and-extensions}, an
upper-triangular perturbation is line-filtered after identifying the
underlying module with a sum of the line modules $\mathcal L_n$. Its decomposition is
therefore covered by Theorem~\ref{thm:line-filtered-splitting}.
\end{remark}

\section{The noncommutative Riccati equation and the splitting problem}
\label{sec:nc-riccati-splitting}

We now study full rank-two connections without a prescribed holomorphic
line filtration. Invariant graph submodules satisfy noncommutative Riccati
equations. We solve them for an explicit constant-matrix family, then
discuss an analytic approach to small perturbations in a chosen completion.
The latter discussion is conditional on estimates specified below.

We retain the left-module convention for the quantum sphere. Let
$
  \mathcal E_{a,b}:=\mathcal L_a\oplus \mathcal L_b
$
and consider a general matrix \(\bar\partial\)-connection
\begin{equation}
  \overline\nabla_B=
  \begin{pmatrix}
    \overline\nabla_{(a)}+\alpha & \beta\\
    \gamma & \overline\nabla_{(b)}+\delta
  \end{pmatrix}.
  \label{eq:full-nondiagonal-dbar}
\end{equation}
Here
\begin{align*}
  \alpha&\in
  \operatorname{Hom}_{\mathcal A_q}
  \bigl(\mathcal L_a,\Omega^{0,1}(\mathcal A_q)\otimes_{\mathcal A_q}\mathcal L_a\bigr),\\
  \beta&\in
  \operatorname{Hom}_{\mathcal A_q}
  \bigl(\mathcal L_b,\Omega^{0,1}(\mathcal A_q)\otimes_{\mathcal A_q}\mathcal L_a\bigr),\\
  \gamma&\in
  \operatorname{Hom}_{\mathcal A_q}
  \bigl(\mathcal L_a,\Omega^{0,1}(\mathcal A_q)\otimes_{\mathcal A_q}\mathcal L_b\bigr),\\
  \delta&\in
  \operatorname{Hom}_{\mathcal A_q}
  \bigl(\mathcal L_b,\Omega^{0,1}(\mathcal A_q)\otimes_{\mathcal A_q}\mathcal L_b\bigr).
\end{align*}
Since \(\Omega^{0,2}(\mathcal A_q)=0\), every such connection is automatically flat.
A holomorphic splitting is described by a holomorphic rank-one idempotent.
When the summands admit graph presentations, their invariance is expressed
by the Riccati equations below. Graph presentations are chart-dependent,
as discussed in Section~\ref{subsec:classical-BG-Riccati}.

\subsection{Invariant graph submodules}
\label{subsec:invariant-graphs}

Let
$
  h\in \operatorname{Hom}_{\mathcal A_q}(\mathcal L_a,\mathcal L_b).
$
Its graph is the left \(\mathcal A_q\)-submodule
\[
  \Gamma_h:=\{(\xi,h\xi):\xi\in \mathcal L_a\}\subset \mathcal L_a\oplus \mathcal L_b.
\]
Using the standard connection on the Hom line module
\(\operatorname{Hom}_{\mathcal A_q}(\mathcal L_a,\mathcal L_b)\simeq \mathcal L_{b-a}\), define
\begin{equation}
  \overline\nabla_{(b-a)}h
  :=\overline\nabla_{(b)}\circ h-(1\otimes h)\circ\overline\nabla_{(a)}.
  \label{eq:hom-connection-h}
\end{equation}
Here and below we suppress $1\otimes -$ when extending a module map to
form-valued sections. The graph \(\Gamma_h\) is preserved
by \(\overline\nabla_B\) if and only if, for every \(\xi\in \mathcal L_a\),
\[
  \overline\nabla_B(\xi,h\xi)
  \in \Omega^{0,1}(\mathcal A_q)\otimes_{\mathcal A_q}\Gamma_h.
\]
Equivalently, the second component of \(\overline\nabla_B(\xi,h\xi)\) must be
obtained by applying \(h\) to the first component. A direct calculation gives
\begin{equation}
  \overline\nabla_{(b-a)}h
  +\gamma+\delta h-h\alpha-h\beta h=0.
  \label{eq:riccati-h}
\end{equation}
This is the first noncommutative Riccati equation.

We call this a \emph{noncommutative Riccati equation} because
it is the invariant-graph equation associated with the block Dolbeault
connection. Its quadratic term arises by the same mechanism that transforms
a two-dimensional linear differential system into a scalar Riccati equation
through the projective coordinate \(h=v/u\). In the present setting the
unknown is a module homomorphism, differentiation is replaced by the
induced Dolbeault connection on the Hom module, and products are ordered
compositions. The holomorphic-idempotent system provides the intrinsic
formulation of the splitting problem; the Riccati equations describe
invariant line submodules in graph coordinates, and complementary
invariant graphs yield a holomorphic splitting.

Similarly, interchanging the roles of $\mathcal L_a$ and $\mathcal L_b$,
for
\[
  k\in \operatorname{Hom}_{\mathcal A_q}(\mathcal L_b,\mathcal L_a),
\]
consider the graph
\[
  \Gamma_k:=\{(k\eta,\eta):\eta\in \mathcal L_b\}\subset \mathcal L_a\oplus \mathcal L_b.
\]
With
\begin{equation}
  \overline\nabla_{(a-b)}k
  :=\overline\nabla_{(a)}\circ k-(1\otimes k)\circ\overline\nabla_{(b)},
  \label{eq:hom-connection-k}
\end{equation}
its invariance is equivalent to
\begin{equation}
  \overline\nabla_{(a-b)}k
  +\beta+\alpha k-k\delta-k\gamma k=0.
  \label{eq:riccati-k}
\end{equation}
The order of composition is essential. In the commutative setting these
formulas reduce to the usual matrix Riccati equations.

\subsection{Transverse solutions and holomorphic splitting}
\label{subsec:transverse-riccati-solutions}

Suppose that \(h\) and \(k\) solve \eqref{eq:riccati-h} and
\eqref{eq:riccati-k}, respectively. The two graph submodules are transverse
precisely when the block matrix
\begin{equation}
  g_{h,k}:=
  \begin{pmatrix}
    1 & k\\
    h & 1
  \end{pmatrix}
  \label{eq:ghk}
\end{equation}
is an automorphism of \(\mathcal L_a\oplus \mathcal L_b\). This is equivalent
to invertibility of either Schur complement, and hence of both:
\begin{equation}
  1-kh\in \operatorname{End}_{\mathcal A_q}(\mathcal L_a),
  \qquad
  1-hk\in \operatorname{End}_{\mathcal A_q}(\mathcal L_b)
  \label{eq:transversality-condition}
\end{equation}
In this case,
\begin{equation}
  g_{h,k}^{-1}=
  \begin{pmatrix}
    (1-kh)^{-1} & -k(1-hk)^{-1}\\[1mm]
    -h(1-kh)^{-1} & (1-hk)^{-1}
  \end{pmatrix}.
  \label{eq:inverse-ghk}
\end{equation}
These inverses must exist in the stated endomorphism algebras, not merely
in a larger completion. The condition does have algebraic examples: if
$h=0$ or $k=0$, both complements are the identity; if $a=b$,
$h=\lambda\operatorname{id}$ and $k=\mu\operatorname{id}$ work whenever
$1-\lambda\mu\neq0$. These examples concern invertibility only; the
Riccati equations must still be satisfied. In general, formal Neumann
series do not establish invertibility in the polynomial algebra. The
small-norm argument used below belongs to an analytic completion.

The two invariant graph submodules give a holomorphic direct-sum decomposition,
and one obtains
\begin{equation}
  g_{h,k}^{-1} \, \overline\nabla_B \, g_{h,k}
  =
  \begin{pmatrix}
    \overline\nabla_{(a)}+\alpha+\beta h & 0\\[1mm]
    0 & \overline\nabla_{(b)}+\delta+\gamma k
  \end{pmatrix}.
  \label{eq:diagonalization-by-riccati}
\end{equation}
Thus two transverse Riccati solutions split the connection into a direct sum of
holomorphic line structures on the underlying line modules \(\mathcal L_a\) and \(\mathcal L_b\).
They need not be the standard KLvS line structures. Equivalence to a standard
structure must be checked separately for each induced connection and the
chosen gauge group.

Equivalently, the first summand is the image of the idempotent
\begin{equation}
  p_{h,k}:=
  g_{h,k}
  \begin{pmatrix}1&0\\0&0\end{pmatrix}
  g_{h,k}^{-1}.
  \label{eq:idempotent-from-graphs}
\end{equation}
Explicitly,
\begin{equation}
  p_{h,k}=
  \begin{pmatrix}
    (1-kh)^{-1} & -k(1-hk)^{-1}\\[1mm]
    h(1-kh)^{-1} & -hk(1-hk)^{-1}
  \end{pmatrix}.
  \label{eq:explicit-projector}
\end{equation}
It satisfies
\begin{equation}
  p_{h,k}^2=p_{h,k}.
  \label{eq:idempotent-relation}
\end{equation}
Moreover, if \(\overline\nabla_{\operatorname{End}}\) denotes the connection induced
by the diagonal background connection
\(\operatorname{diag}(\overline\nabla_{(a)},\overline\nabla_{(b)})\) on
\(\operatorname{End}_{\mathcal A_q}(\mathcal E_{a,b})\), then holomorphicity of the projector is
\begin{equation}
  \overline\nabla^{}_{\operatorname{End}} p_{h,k} +Bp_{h,k}-p_{h,k}B=0,
  \label{eq:holomorphic-projector-equation}
\end{equation}
where
\[
  B=\begin{pmatrix}\alpha&\beta\\ \gamma&\delta\end{pmatrix}.
\]
Conversely, the image of a holomorphic rank-one projector is represented
by a graph over $\mathcal L_a$ whenever its projection onto
$\mathcal L_a$ is a module isomorphism; the analogous statement holds for
$\mathcal L_b$. Complementary summands give two graph solutions when both
lie in the corresponding graph charts. This need not hold globally relative
to one fixed decomposition: a summand of a different underlying line type
cannot be a global graph over $\mathcal L_a$ or $\mathcal L_b$.

Within a graph chart for the two summands, these descriptions agree:
\begin{equation}
  \begin{array}{c}
  \text{holomorphic rank-one splitting of }(\mathcal E_{a,b},\overline\nabla_B)\\[1mm]
  \Longleftrightarrow\\[-1mm]
  \text{holomorphic rank-one projector}\\[1mm]
  \Longleftrightarrow\\[-1mm]
  \text{two transverse invariant graph submodules}\\[1mm]
  \Longleftrightarrow\\[-1mm]
  \text{two transverse solutions of the Riccati equations.}
  \end{array}
  \label{eq:riccati-equivalence-summary}
\end{equation}
For the global problem, the idempotent formulation does not require a fixed
graph chart.

\subsection{Relation with the classical Birkhoff--Grothendieck theorem}
\label{subsec:classical-BG-Riccati}

The Riccati equation already appears in the classical case.
 Let
\(E=L_a\oplus L_b\) be a smooth splitting of a rank-two complex vector bundle
on the classical Riemann sphere and write a Dolbeault operator in block form,
\[
  \bar\partial_E=
  \begin{pmatrix}
    \bar\partial_a+\alpha & \beta\\
    \gamma & \bar\partial_b+\delta
  \end{pmatrix}.
\]
A line subbundle written as the graph of a smooth map \(h:L_a\to L_b\) is
holomorphic exactly when
\begin{equation}
  \bar\partial_{\operatorname{Hom}}h
  +\gamma+\delta h-h\alpha-h\beta h=0,
  \label{eq:classical-riccati-formal}
\end{equation}
which is the ordinary matrix Riccati equation written in Dolbeault form. Thus
the classical Birkhoff--Grothendieck theorem may be viewed as a global
solvability theorem for such equations, provided one allows changes of chart in
the Grassmannian of line subbundles when the chosen graph representation breaks
down.

The KLvS splitting problem asks whether every rank-two holomorphic bundle
similarly admits complementary holomorphic rank-one summands after
replacing the ordinary Dolbeault calculus by the covariant calculus of the
standard Podle\'s sphere. The formal equations have the same shape, but the unknowns
now lie in noncommutative Hom line modules and the products in
\eqref{eq:riccati-h} and \eqref{eq:riccati-k} are ordered compositions.

\subsection{An explicit solvable non-diagonal family}
\label{subsec:explicit-solvable-riccati-family}

There is a simple family where the Riccati equation can be solved
exactly. Take \(a=b=n\), so that
\[
  \mathcal E=\mathcal L_n\oplus \mathcal L_n.
\]
Let
\begin{equation}
  \Theta\in
  \operatorname{Hom}_{\mathcal A_q}
  \bigl(\mathcal L_n,\Omega^{0,1}(\mathcal A_q)\otimes_{\mathcal A_q}\mathcal L_n\bigr)
  \label{eq:theta-endomorphism-valued-form}
\end{equation}
be any \(\mathcal A_q\)-linear \((0,1)\)-form-valued endomorphism of \(\mathcal L_n\). With
complex numbers \(u,v,b,c\),  put
\begin{equation}
  M=\begin{pmatrix}u&b\\ c&v\end{pmatrix}.
  \label{eq:constant-matrix-M}
\end{equation}
Consider the connection
\begin{equation}
  \overline\nabla_M=\overline\nabla_{(n)}I_2+
  \Theta M
  =
  \begin{pmatrix}
    \overline\nabla_{(n)}+u\Theta & b\Theta\\
    c\Theta & \overline\nabla_{(n)}+v\Theta
  \end{pmatrix}.
  \label{eq:constant-coefficient-family}
\end{equation}
Both off-diagonal entries are nonzero when $bc\neq0$ and $\Theta\neq0$.
Under the distinct-eigenvalue assumption below, a constant change of frame
diagonalizes $M$, and the Riccati equations admit constant graph solutions.

If $\Theta=0$, the connection is already diagonal, so assume
$\Theta\neq0$ in the scalar reduction. Let \(h\in\mathbb C\).
Since \(\overline\nabla_{(0)}h=0\), equation
\eqref{eq:riccati-h} becomes
\begin{equation}
  c+(v-u)h-bh^2=0,
  \label{eq:scalar-riccati-example}
\end{equation}
or
\begin{equation}
  bh^2+(u-v)h-c=0.
  \label{eq:scalar-riccati-example-standard}
\end{equation}
Assume \(b\neq0\), \(c\neq0\), and
\begin{equation}
  \Delta:=(u-v)^2+4bc\neq0.
  \label{eq:discriminant-riccati}
\end{equation}
Then the two roots are
\begin{equation}
  h_{\pm}=\frac{v-u\pm\sqrt{\Delta}}{2b}.
  \label{eq:two-h-roots}
\end{equation}
They define two distinct graph submodules and hence two transverse invariant
line summands. With
\begin{equation}
  S=\begin{pmatrix}1&1\\ h_+&h_-\end{pmatrix},
  \label{eq:diagonalizing-matrix-S}
\end{equation}
equation \eqref{eq:scalar-riccati-example-standard} gives
\begin{equation}
  S^{-1}\overline\nabla_M S
  =
  \begin{pmatrix}
    \overline\nabla_{(n)}+\lambda_+\Theta & 0\\
    0 & \overline\nabla_{(n)}+\lambda_-\Theta
  \end{pmatrix},
  \label{eq:explicit-diagonalization}
\end{equation}
where
\begin{equation}
  \lambda_{\pm}=\frac{u+v\pm\sqrt{\Delta}}{2}.
  \label{eq:eigenvalues-lambda}
\end{equation}
Thus
\begin{equation}
  (\mathcal L_n^{\oplus2},\overline\nabla_{(n)}I_2+\Theta M)
  \simeq
  (\mathcal L_n,\overline\nabla_{(n)}+\lambda_+\Theta)
  \oplus
  (\mathcal L_n,\overline\nabla_{(n)}+\lambda_-\Theta).
  \label{eq:explicit-family-splits}
\end{equation}
For $\Theta\neq0$, this gives a splitting with both off-diagonal entries
nonzero in the original presentation. The resulting line structures may be
nonstandard; their equivalence to standard KLvS structures must be checked
separately for these particular perturbations and the chosen gauge group.

Equivalently, the two holomorphic projections are the spectral projections
\begin{equation}
  p_{\pm}=\frac{M-\lambda_{\mp}I_2}{\lambda_{\pm}-\lambda_{\mp}}.
  \label{eq:spectral-projections}
\end{equation}
They are constant matrices. Hence \(\overline\nabla^{}_{\operatorname{End}} \, p_{\pm}=0\),
and since \([M,p_{\pm}]=0\), they satisfy the holomorphic projector equation
\eqref{eq:holomorphic-projector-equation}.

In the particularly transparent subfamily
\[
  u=v=0,
  \qquad bc\neq0,
\]
choose \(\mu\in\mathbb C\) with \(\mu^2=bc\). Then
\[
  \lambda_{\pm}=\pm\mu,
  \qquad
  h_{\pm}=\pm\frac{\mu}{b},
\]
and
\begin{equation}
  \overline\nabla_{(n)}I_2+
  \begin{pmatrix}0&b\Theta\\ c\Theta&0\end{pmatrix}
  \simeq
  (\overline\nabla_{(n)}+\mu\Theta)
  \oplus
  (\overline\nabla_{(n)}-\mu\Theta).
  \label{eq:symmetric-offdiagonal-example}
\end{equation}
This example reduces the graph equation to a scalar quadratic while
retaining arbitrary form-valued data in $\Theta$.

\subsection{Analytic perspectives and the remaining global problem}
\label{subsec:local-analytic-triangularization}

The following discussion describes a possible analytic approach rather
than a theorem for an unspecified completion. Here \emph{local} refers to
a neighbourhood of the standard diagonal connection in the space of
connections, not to a coordinate neighbourhood on the quantum sphere.
Fix a smooth Fr\'echet algebra $\mathcal A_q^\infty$ carrying the KLvS
calculus, together with Banach or Sobolev completions of the line modules
and their form-valued Hom spaces. At a fixed regularity index $s$, assume
that
\[
 \overline\nabla_{(m)}:
 W^{s+1}(\mathcal L_m)\longrightarrow
 W^s\bigl(\Omega^{0,1}(\mathcal A_q^\infty)
       \otimes_{\mathcal A_q^\infty}\mathcal L_m\bigr)
\]
is bounded and has a bounded right inverse with image in a fixed
complement of its kernel whenever the relevant degree-one cohomology
vanishes. Also assume continuous multiplication and composition estimates
in these norms. For the analytic argument to use the polynomial
cohomology formula
\begin{equation}
 \dim H^{0,1}(\mathcal L_m)=\max\{m-1,0\},
 \label{eq:h01-line-formula-riccati}
\end{equation}
one must additionally establish its persistence in the chosen completion.
None of these estimates follows merely from algebraic vanishing of the
cokernel. A smooth conclusion would further require regularity or
bootstrapping from the completed solution space.

\paragraph{Conditional triangularization.}
Assume first that \(a\geq b\). Then \(b-a\leq0\), and hence
\[
  H^{0,1}(\mathcal L_{b-a})=0.
\]
Let
\[
  G_{b-a}:
  \Omega^{0,1}(\mathcal A_q^\infty)\otimes_{\mathcal A_q^\infty}\mathcal L_{b-a}
  \longrightarrow \mathcal L_{b-a}
\]
be a right inverse for \(\overline\nabla_{(b-a)}\), bounded between the
chosen completed spaces and taking values in the fixed complement of the
holomorphic kernel.
Equation \eqref{eq:riccati-h} may be rewritten as the fixed-point equation
\begin{equation}
  h=
  G_{b-a}\bigl(-\gamma-\delta h+h\alpha+h\beta h\bigr).
  \label{eq:fixed-point-h}
\end{equation}
Under the assumed multiplication and right-inverse bounds, the right-hand
side maps a sufficiently small closed ball into itself and is a contraction
when \(\alpha,\beta,\gamma,\delta\) are sufficiently small. Hence the graph equation
has a unique small solution $h$ in the chosen complement. The unipotent gauge transformation
\begin{equation}
  g_h=\begin{pmatrix}1&0\\ h&1\end{pmatrix}
  \label{eq:unipotent-gauge-h}
\end{equation}
then puts the connection into upper-triangular form:
\begin{equation}
  g_h^{-1}\overline\nabla_B g_h
  =
  \begin{pmatrix}
    \overline\nabla_{(a)}+\alpha+\beta h & \beta\\
    0 & \overline\nabla_{(b)}+\delta-h\beta
  \end{pmatrix}.
  \label{eq:local-triangularization-h}
\end{equation}
Similarly, if \(b\geq a\), then \(a-b\leq0\), and the same argument applied to
\eqref{eq:riccati-k} gives a small solution \(k\) and a lower-triangular form.

Thus the stated analytic hypotheses would imply triangularization of
sufficiently small full matrix perturbations in the completed category.
Establishing those hypotheses for a concrete KLvS completion is a separate
analytic problem.

\paragraph{Adjacent degrees and conditional splitting.}
\label{subsec:local-analytic-splitting-adjacent}

If
\begin{equation}
  |a-b|\leq1,
  \label{eq:adjacent-degrees-condition}
\end{equation}
then both
\[
  a-b\leq1,
  \qquad
  b-a\leq1.
\]
Therefore
\[
  H^{0,1}(\mathcal L_{a-b})=H^{0,1}(\mathcal L_{b-a})=0.
\]
Under the same analytic hypotheses, the two
Riccati equations \eqref{eq:riccati-h} and \eqref{eq:riccati-k} have small
solutions for every sufficiently small perturbation \(\alpha,\beta,\gamma,\delta\).
If the chosen norms make the endomorphism spaces Banach algebras,
smallness of \(h\) and \(k\) makes \(1-kh\) and \(1-hk\) invertible there
by Neumann series. Hence the two graph submodules are transverse
and the connection is diagonalized by \(g_{h,k}\). We obtain
\begin{equation}
  (\mathcal E_{a,b},\overline\nabla_B)
  \simeq
  (\mathcal L_a,\overline\nabla_{(a)}+\alpha+\beta h)
  \oplus
  (\mathcal L_b,\overline\nabla_{(b)}+\delta+\gamma k).
  \label{eq:local-line-splitting}
\end{equation}
This is the conditional local analytic splitting statement for adjacent
degrees.

If one also assumes local line-rigidity in the same regularity class, then the
line structures appearing in \eqref{eq:local-line-splitting} are gauge equivalent
to the standard KLvS structures. Formally, the relevant infinitesimal statement
is controlled by
\[
  H^{0,1}\bigl(\operatorname{End}_{\mathcal A_q}(\mathcal L_n)\bigr)
  \simeq H^{0,1}(\mathcal L_0)=0.
\]
With the additional local line-rigidity assumption, solving the
corresponding nonlinear gauge equation would give the standard splitting
\begin{equation}
  (\mathcal E_{a,b},\overline\nabla_B)
  \simeq
  (\mathcal L_a,\overline\nabla_{(a)})
  \oplus
  (\mathcal L_b,\overline\nabla_{(b)}).
  \label{eq:local-standard-splitting}
\end{equation}
However, this last standard-form conclusion uses local line-rigidity in the
chosen analytic completion and must be kept separate from the existence of
the holomorphic rank-one splitting. It makes no assertion about polynomial
gauge transformations: Proposition~\ref{prop:polynomial-line-nonrigidity}
shows that polynomial line-rigidity fails even for arbitrarily small nonzero
perturbations of the trivial line connection.

\paragraph{The global existence problem.}
\label{subsec:remaining-global-riccati-problem}

The established polynomial results cover line-filtered bundles and the
constant-matrix family under \eqref{eq:discriminant-riccati}. The analytic
discussion gives a conditional route for small perturbations, subject to
the completion and regularity hypotheses above. They do not resolve
the unrestricted global problem for arbitrary $\alpha,\beta,\gamma,\delta$.
One must find complementary invariant rank-one summands, allowing changes
of graph chart and of splitting type.

The global idempotent formulation asks whether every holomorphic
connection of the form \eqref{eq:full-nondiagonal-dbar} admits a nontrivial
holomorphic rank-one projector. In formula form, one seeks an idempotent
\begin{equation}
  p\in \operatorname{End}_{\mathcal A_q}(\mathcal L_a\oplus \mathcal L_b),
  \qquad p^2=p,
  \label{eq:global-idempotent-problem}
\end{equation}
of generic rank one satisfying
\begin{equation}
  \overline\nabla^{}_{\operatorname{End}}p+Bp-pB=0.
  \label{eq:global-projector-equation}
\end{equation}
A solution automatically gives the complementary holomorphic projector
$1-p$, and hence a splitting into rank-one projective summands. Identifying
their underlying modules with canonical $\mathcal L_n$, and their
holomorphic structures with the standard KLvS structures, are additional
classification questions.

The projector equation is global, while the Riccati equations describe
invariance in graph charts. A single invariant graph supplies a holomorphic
line submodule; a direct-sum decomposition requires a complementary
invariant summand. The higher-rank analogue is a Grassmann--Riccati equation
for invariant submodules. These formulations separate the explicit splitting
results above from the global existence problem for general KLvS bundles.

\section{Sheaves, associated bundles, and the algebraic quantum projective line}
\label{sec:sheaf-approach}

The purpose of this section is to formulate the algebraic, sheaf-theoretic side of
the quantum Birkhoff--Grothendieck problem.  The point of departure is the
sheaf approach to quantum principal bundles developed in 
\cite{AschieriFioresiLatini,AschieriFioresiLatiniWeber} inspired by \cite{Pflaum}.  We shall also compare
this approach with the quotient-category construction of Artin and Zhang and
with twisted homogeneous coordinate rings in the sense of Artin and Van den
Bergh \cite{ArtinZhang,ArtinVdB,ZhangTwist}.

Several categories occur naturally and must be kept separate.  The first is the
category of quantum principal bundles, encoded by sheaves of comodule algebras.
The second is the category of vector bundles associated to finite-dimensional
comodules of the structure Hopf algebra.  The third is the larger category of
all finite locally free modules over the base structure sheaf.  Finally, to make
contact with Section~\ref{sec:klv-module-picture}, one must enrich the sheaf
picture by a differential calculus and a compatible Dolbeault operator.  Passing
from one of these categories to another can forget substantial information.
In particular, the classification of principal bundles or of equivariant
associated bundles is finer than the classification of the underlying vector
bundles.

For the quantum projective line an important simplification occurs:
the two affine coordinate algebras of the base, and also their overlap algebra,
are commutative polynomial or Laurent polynomial algebras.  Thus, on the topology generated by the standard two-chart cover of \cite{AschieriFioresiLatini} the bare structure sheaf of the quantum projective
line is the restriction of the ordinary structure sheaf of $\mathbb P^1$,
although its homogeneous coordinate
ring, the total principal-bundle sheaf, and the natural differential calculi are
noncommutative.  It follows that the Birkhoff--Grothendieck theorem already holds
for finite locally free modules over this bare structure sheaf.  This statement
is useful, but it does not settle weak holomorphic splitting in the KLvS
category: passage to the bare structure sheaf forgets the differential and
regularity data defining the holomorphic problem, as well as any separately
retained equivariance.

\subsection{Quantum ringed spaces and sheaves of principal comodule algebras}
\label{subsec:quantum-ringed-principal}

Let $\Bbbk$ be a field and let $H$ be a Hopf algebra over $\Bbbk$, with coproduct
$\Delta$, counit $\varepsilon$, and bijective antipode $S$.  A \emph{quantum
ringed space} is a pair $(M,\mathcal O_M)$ consisting of a topological space $M$
and a sheaf $\mathcal O_M$ of not necessarily commutative unital
$\Bbbk$-algebras.  This terminology is deliberately flexible: the underlying
space records the chosen localization topology, whereas the sheaf records the
possibly noncommutative local coordinate algebras.

A sheaf $\mathcal F$ of right $H$-comodule algebras assigns to every open set
$U\subseteq M$ a right $H$-comodule algebra
\[
  \delta_U:\mathcal F(U)\longrightarrow \mathcal F(U)\otimes H,
\]
with all restriction maps $H$-colinear algebra homomorphisms.  Its coinvariant
subsheaf is
\[
  \mathcal F^{\operatorname{co}H}(U)
  :=\{f\in\mathcal F(U):\delta_U(f)=f\otimes 1\}.
\]
The basic requirement is that the base sheaf be recovered as the coinvariant
subsheaf,
\begin{equation}
  \mathcal O_M=\mathcal F^{\operatorname{co}H}.
  \label{eq:base-as-coinvariants-sheaf}
\end{equation}

For an affine extension $B=A^{\operatorname{co}H}\subseteq A$, the canonical map
is
\begin{equation}
  \operatorname{can}:A\otimes_B A\longrightarrow A\otimes H,
  \qquad
  a\otimes_B a'\longmapsto aa'_{(0)}\otimes a'_{(1)}.
  \label{eq:canonical-map-affine}
\end{equation}
The extension is Hopf--Galois when this map is bijective.  In 
\cite{AschieriFioresiLatini} a principal comodule algebra also satisfies an
equivariant projectivity condition; this is the algebraic substitute for the
local freeness of an ordinary principal action.

\begin{definition}[Sheaf-theoretic quantum principal bundle]
\label{def:sheaf-quantum-principal}
Let $(M,\mathcal O_M)$ be a quantum ringed space and let $H$ be a Hopf algebra.
A sheaf $\mathcal F$ of right $H$-comodule algebras is an $H$-principal bundle
on $(M,\mathcal O_M)$ if
\[
  \mathcal F^{\operatorname{co}H}=\mathcal O_M
\]
and there exists an open cover $\{U_i\}$ of $M$ such that every extension
\[
  \mathcal O_M(U_i)
  =\mathcal F(U_i)^{\operatorname{co}H}
  \subseteq \mathcal F(U_i)
\]
is a principal $H$-comodule algebra.

It is \emph{locally cleft} if one may choose the cover so that every local
extension is cleft, equivalently, if there are convolution-invertible right
$H$-colinear maps
\[
  j_i:H\longrightarrow\mathcal F(U_i),
  \qquad j_i(1)=1.
\]
It is \emph{locally trivial} if the $j_i$ can be chosen as algebra maps.
\end{definition}

A cleaving map gives an isomorphism of left
$\mathcal O_M(U_i)$-modules and right $H$-comodules
\begin{equation}
  \mathcal O_M(U_i)\otimes H
  \xrightarrow{\ \simeq\ }
  \mathcal F(U_i),
  \qquad b\otimes h\longmapsto b\,j_i(h).
  \label{eq:local-cleft-trivialization}
\end{equation}
The algebra structure on the left-hand side may be a crossed-product structure;
only in the locally trivial case is it the appropriate smash or tensor product
algebra without a nontrivial cocycle.  This distinction has no classical analogue
for an ordinary trivial principal bundle, but it is fundamental in the quantum
setting.

The sheaf formulation is forced upon us by projective geometry.  If $G/P$ is a
projective homogeneous variety, the algebra of global regular functions on the
base is usually only $\Bbbk$.  Consequently, the global coinvariant inclusion
$\Bbbk\subseteq\mathcal O(G)$ does not retain the geometry of $G/P$ and is not the
correct affine Hopf--Galois model.  One instead localizes the total algebra and
the homogeneous coordinate ring on a projective cover and imposes principality
locally.  In this way the base is recovered from projective localizations rather
than from global coinvariants.

\subsection{Associated sheaves and the loss of equivariant information}
\label{subsec:associated-sheaves}

Let $V$ be a finite-dimensional left $H$-comodule, with coaction
\[
  \lambda_V(v)=v_{(-1)}\otimes v_{(0)}.
\]
For a right $H$-comodule algebra $A$, the usual associated module is the cotensor
product
\[
  A\square_H V
  =\ker\bigl(
       \delta_A\otimes\operatorname{id}_V
       -\operatorname{id}_A\otimes\lambda_V
     \bigr)
  \subseteq A\otimes V.
\]
The same construction may be performed open set by open set.

\begin{definition}[Associated sheaf]
\label{def:associated-sheaf}
Let $\mathcal F$ be an $H$-principal bundle on $(M,\mathcal O_M)$ and let $V$ be a
finite-dimensional left $H$-comodule.  The associated sheaf is
\begin{equation}
  \mathcal E_V(U):=\mathcal F(U)\square_H V.
  \label{eq:associated-sheaf-cotensor}
\end{equation}
It is naturally a left $\mathcal O_M(U)$-module by multiplication on the first
factor.
\end{definition}

Because $V$ is finite-dimensional, $\mathcal F\otimes V$ is a sheaf;
the cotensor condition can be checked locally since restriction maps are
$H$-colinear. Thus \eqref{eq:associated-sheaf-cotensor} defines a sheaf.
If $\mathcal F$ is
\emph{locally cleft}, choose a cleft cover. On each chart of this cover,
\eqref{eq:local-cleft-trivialization} gives
\begin{equation}
  \mathcal E_V(U_i)
  \simeq
  \bigl(\mathcal O_M(U_i)\otimes H\bigr)\square_H V
  \simeq
  \mathcal O_M(U_i)\otimes V.
  \label{eq:associated-local-free}
\end{equation}
Hence associated sheaves of a locally cleft principal sheaf are finite
locally free. This applies to the locally trivial quantum projective-line
examples considered below. For a merely locally principal sheaf, the affine
theory gives finitely generated projective associated modules on the principal
charts \cite{BrzezinskiHajac}. That projectivity alone does not imply freeness
on those charts, or local freeness on an arbitrary topology in
Definition~\ref{def:sheaf-quantum-principal}; a suitable refinement or an
additional local hypothesis is needed.

The converse is false at this level of generality.  A locally free
$\mathcal O_M$-module need not come equipped with a lift to the total sheaf
$\mathcal F$, nor need it arise from a finite-dimensional $H$-comodule.  Thus the
essential image of
\[
  V\longmapsto \mathcal E_V
\]
is the collection of underlying associated bundles. An enriched associated
category additionally retains the comodule or equivariant data; it need not
be a full subcategory of the category of all vector bundles.  Even within this image, forgetting the $H$-comodule that produced the
bundle can identify objects that were inequivalent before forgetting the
principal-bundle data.

The classical bundle $\mathrm{SL}_2\to\mathbb P^1$ already illustrates this
loss.  A two-dimensional representation of the Borel subgroup may be an
extension of two characters.  The associated sequence of line bundles can be
non-split as a homogeneous or equivariant sequence even when its middle term is
a split ordinary vector bundle.  The Euler sequence
\begin{equation}
  0\longrightarrow\mathcal O(-1)
  \longrightarrow\mathcal O^{\oplus 2}
  \longrightarrow\mathcal O(1)
  \longrightarrow 0
  \label{eq:euler-sequence-associated}
\end{equation}
is the basic example.  Thus a Birkhoff--Grothendieck theorem classifies the
underlying middle bundle; it does not recover the representation, the extension
class in an equivariant category, or a reduction of the frame bundle to a
smaller structure group.  The same warning applies, with greater force, to the
quantum principal-bundle setting.

\subsection{Projective Ore localization and quantum projective spaces}
\label{subsec:projective-ore-localization}

Let
\[
  A=\bigoplus_{n\geq 0}A_n
\]
be a connected graded algebra generated in degree one.  Suppose that homogeneous
elements $x_0,\ldots,x_N\in A_1$ determine Ore localizations.  The projective
local algebra on the chart $U_i$ is the degree-zero part
\begin{equation}
  \mathcal O_A(U_i)
  :=\bigl(A[x_i^{-1}]\bigr)_0,
  \label{eq:projective-localization-degree-zero}
\end{equation}
and similarly on multiple intersections.  Under the appropriate compatibility
conditions these algebras, together with their restriction maps, define a sheaf
on the standard projective cover.  This is the noncommutative analogue of the
usual construction of $\operatorname{Proj}A$ from degree-zero localizations.

For the  quantum projective space one begins with
\begin{equation}
  A_q^{(N)}
  =\Bbbk\langle x_0,\ldots,x_N\rangle
   \big/
   \bigl(x_i x_j-qx_jx_i:\ i<j\bigr),
  \label{eq:q-projective-coordinate-ring}
\end{equation}
where changing $q$ to $q^{-1}$ gives the convention used in parts of the
quantum-group literature.  The generators are normal and their powers form Ore
sets.  The chart algebras in
\eqref{eq:projective-localization-degree-zero} are quantum affine spaces of
dimension $N$.

At the same time in \cite{AschieriFioresiLatini}, a total sheaf is constructed by localizing
$\mathcal O_q(\mathrm{SL}_{N+1})$ along the corresponding quantum minors.  Its
coinvariant subsheaf is the projectively localized sheaf associated with
\eqref{eq:q-projective-coordinate-ring}, and the resulting total sheaf is locally
principal, in the standard examples even locally trivial, for the quantum
parabolic Hopf algebra $\mathcal O_q(P)$.  The construction therefore contains
three related but distinct pieces of information:
\[
  \text{the homogeneous ring }A_q^{(N)},
  \qquad
  \text{the base sheaf }\mathcal O_{q,\mathbb P^N},
  \qquad
  \text{the principal sheaf }\mathcal F_q.
\]
They should not be identified.

\subsection{The special feature of the quantum projective line}
\label{subsec:quantum-line-hidden-commutativity}

Set
\begin{equation}
  A_q
  :=\Bbbk\langle x_0,x_1\rangle/(x_0x_1-qx_1x_0),
  \qquad q\in\Bbbk^\times.
  \label{eq:quantum-plane-homogeneous-ring}
\end{equation}
Let $U_0$ and $U_1$ be the two standard projective charts.  In the Ore
localizations define the degree-zero coordinates
\[
  z:=x_0^{-1}x_1,
  \qquad
  w:=x_1^{-1}x_0.
\]
There is only one affine coordinate on each chart, and consequently
\begin{equation}
  \mathcal O_{q,\mathbb P^1}(U_0)=\Bbbk[z],
  \qquad
  \mathcal O_{q,\mathbb P^1}(U_1)=\Bbbk[w],
  \label{eq:qP1-chart-algebras}
\end{equation}
while on the overlap
\begin{equation}
  \mathcal O_{q,\mathbb P^1}(U_0\cap U_1)
  =\Bbbk[z,z^{-1}],
  \qquad w=z^{-1}.
  \label{eq:qP1-overlap-algebra}
\end{equation}
All three algebras are commutative.  Therefore, after fixing the ground field and specializing $q$, the base ringed
space on the topology generated by this two-chart cover is isomorphic to the
corresponding restriction of the ordinary algebraic projective line:
\begin{equation}
  (\mathbb P^1,\mathcal O_{q,\mathbb P^1})
  \simeq
  (\mathbb P^1,\mathcal O_{\mathbb P^1})
  \quad\text{on the standard two-chart topology.}
  \label{eq:qP1-base-is-classical}
\end{equation}
Over the parameter ring $\Bbbk[q,q^{-1}]$ the same statement holds relatively.

Equation \eqref{eq:qP1-base-is-classical} is easy  to overlook.  It does not say
that the quantum projective-line construction is trivial.  The parameter $q$
remains visible in the following data:
\begin{enumerate}[label=\textup{(\roman*)},leftmargin=2.3em]
\item the noncommutative multiplication of the homogeneous coordinate ring
      $A_q$;
\item the total sheaf $\mathcal F_q$ of $\mathcal O_q(P)$-comodule algebras and
      its local trivializations;
\item the tensor and equivariance data of associated sheaves;
\item the differential calculus induced from the quantum group.
\end{enumerate}
Indeed \cite{AschieriFioresiLatiniWeber} exhibits calculi for which the base
structure sheaf is commutative but the bimodule of one-forms is noncommutative.
Thus the word ``quantum'' in this model is carried by more than the multiplication
on the local base algebras.

This also shows why the sheaf approach does not immediately settle the KLvS
problem.  The KLvS algebra of the standard Podle\'s sphere is globally
noncommutative, and a holomorphic bundle in the KLvS sense includes a
$\bar\partial$-connection.  Passing to the bare sheaf
\eqref{eq:qP1-base-is-classical} forgets this noncommutative differential data.
A comparison theorem must therefore involve enriched sheaves, not merely locally
free modules over $\mathcal O_{q,\mathbb P^1}$.

\subsection{The standard line sheaves and their section ring}
\label{subsec:standard-line-sheaves}

Let $H=\mathcal O_q(P)$ be the quantum Borel or parabolic Hopf algebra appearing
in the principal sheaf $\mathcal F_q$.  Its standard group-like element $t$
gives one-dimensional comodules $\Bbbk_n$, $n\in\mathbb Z$, by powers $t^n$.
Up to a sign convention, define
\begin{equation}
  \mathscr O_q(n)
  :=\mathcal F_q\square_H\Bbbk_n.
  \label{eq:associated-line-sheaf}
\end{equation}
On each standard chart this sheaf is free of rank one.  If $e_0^{(n)}$ and
$e_1^{(n)}$ are local frames obtained from the two cleaving maps, then on the
overlap
\begin{equation}
  e_1^{(n)}=z^n e_0^{(n)}
  \label{eq:line-sheaf-transition}
\end{equation}
after a harmless nonzero scalar normalization and, if necessary, replacing $n$
by $-n$ according to the character convention.  Consequently, as a sheaf of
modules over the base,
\begin{equation}
  \mathscr O_q(n)\simeq\mathcal O_{\mathbb P^1}(n).
  \label{eq:associated-line-is-classical}
\end{equation}
For the underlying sheaves, regarded as central modules over the
commutative base, the ordinary tensor and duality laws have the form
\begin{equation}
  \mathscr O_q(m)\otimes\mathscr O_q(n)
  \simeq\mathscr O_q(m+n),
  \qquad
  \mathscr O_q(n)^\vee\simeq\mathscr O_q(-n).
  \label{eq:line-sheaf-tensor-law}
\end{equation}

\begin{proposition}[Classification of line bundles in the bare sheaf category]
\label{prop:line-bundles-bare-sheaf}
Every finite locally free rank-one
$\mathcal O_{q,\mathbb P^1}$-module is isomorphic to exactly one
$\mathscr O_q(n)$, $n\in\mathbb Z$.  Thus
\[
  \operatorname{Pic}(\mathbb P^1,\mathcal O_{q,\mathbb P^1})
  \simeq\mathbb Z.
\]
\end{proposition}

\begin{proof}
By \eqref{eq:qP1-base-is-classical}, the descent problem on this cover is the
usual computation of the Picard group of $\mathbb P^1$.  The characters $t^n$ realize every integer class through
\eqref{eq:associated-line-sheaf}.
\end{proof}

The proposition concerns one-sided locally free sheaves, or equivalently central
line bundles over the commutative base sheaf.  If one enlarges the category to
noncentral invertible bimodules, or if one remembers distinct equivariant
structures on the same underlying line sheaf, additional objects may occur.  The
category must therefore be specified whenever one asks whether the
$\mathscr O_q(n)$ exhaust all ``line bundles.''

There is no contradiction between
\eqref{eq:associated-line-is-classical} and the noncommutativity of the section
ring.  As graded vector spaces,
\[
  \bigoplus_{n\geq0}H^0\bigl(\mathbb P^1,\mathscr O_q(n)\bigr)
  \simeq
  \bigoplus_{n\geq0}H^0\bigl(\mathbb P^1,\mathcal O(n)\bigr),
\]
but multiplication is inherited from the ordered multiplication of
semi-coinvariant elements in the total quantum algebra.  It is therefore a
twisted multiplication, and one recovers
\begin{equation}
  \bigoplus_{n\geq0}H^0\bigl(\mathbb P^1,\mathscr O_q(n)\bigr)
  \simeq
  \Bbbk\langle x_0,x_1\rangle/(x_0x_1-qx_1x_0)
  \label{eq:q-section-ring}
\end{equation}
as graded algebras. This multiplication uses the  bimodule data
of the semi-coinvariants. It must be distinguished from the ordinary tensor
product of the underlying central line sheaves in
\eqref{eq:line-sheaf-tensor-law}; that ordinary tensor law by itself does
not specify the quantum multiplication.

This phenomenon is naturally expressed by twisted homogeneous coordinate rings.
Let
\[
  \sigma_q:\mathbb P^1\longrightarrow\mathbb P^1,
  \qquad
  \sigma_q^*(z)=qz,
\]
and set
\[
  \mathcal M_n
  =\mathcal O(1)\otimes\sigma_q^*\mathcal O(1)
   \otimes\cdots\otimes(\sigma_q^{n-1})^*\mathcal O(1).
\]
Each $\mathcal M_n$ is isomorphic to $\mathcal O(n)$ as a line bundle, but the
multiplication
\[
  H^0(\mathcal M_m)\otimes H^0(\mathcal M_n)
  \longrightarrow H^0(\mathcal M_{m+n})
\]
twists the second factor by $\sigma_q^m$.  The resulting ring
$B(\mathbb P^1,\mathcal O(1),\sigma_q)$ is the quantum plane
\eqref{eq:quantum-plane-homogeneous-ring}. Thus the sheaves can be classical
as objects while their twisted graded section algebra is noncommutative.

\subsection{Algebraic vector bundles on the affine line}
\label{subsec:affine-line-admissibility}

Before gluing the two projective charts, it is useful to state precisely what is
meant by an algebraic vector bundle on an affine chart.  For
$U=\operatorname{Spec}\Bbbk[z]$, a rank-$r$ algebraic vector bundle is a finite
locally free sheaf of $\mathcal O_U$-modules.  By the affine correspondence it is
equivalent to a finitely generated projective $\Bbbk[z]$-module.  Since
$\Bbbk[z]$ is a principal ideal domain, every such module is free.  More
generally, the Quillen--Suslin theorem gives the same conclusion over affine
space $\Bbbk[z_1,\ldots,z_N]$ \cite{Quillen,Suslin}.

Thus there is only one rank-$r$ algebraic vector bundle on the affine line up to
bundle isomorphism:
\begin{equation}
  \mathcal E\simeq\mathcal O_{\mathbb A^1}^{\oplus r}.
  \label{eq:affine-line-bundles-trivial}
\end{equation}
Here equivalence means an $\mathcal O_U$-linear sheaf isomorphism, or equivalently
an isomorphism of the corresponding finitely generated projective modules.  For
a chosen trivialization, changes of frame form
$\operatorname{GL}_r(\Bbbk[z])$.  Stable equivalence, Morita equivalence of base
algebras, and equivalence after adding a connection are different relations and
should not replace ordinary bundle isomorphism without an explicit reason.

A connection is extra structure.  On the trivial algebraic bundle one may write
\begin{equation}
  \nabla=d+\omega(z)\,dz,
  \qquad \omega(z)\in M_r(\Bbbk[z]).
  \label{eq:algebraic-connection-affine-line}
\end{equation}
Since there are no two-forms on the affine line, every such connection is flat.
Nevertheless, two connections are algebraically gauge equivalent only when they
are related by a polynomial gauge transformation
$g\in\operatorname{GL}_r(\Bbbk[z])$:
\begin{equation}
  \omega \longmapsto g \omega g^{-1}-(dg)g^{-1}.
  \label{eq:algebraic-gauge-affine-line}
\end{equation}
For $r=1$, the only units of $\Bbbk[z]$ are nonzero  constants, so the polynomial
one-form $A(z)\,dz$ cannot in general be gauged away algebraically.  Over the
holomorphic affine line, by contrast, a flat holomorphic connection can be
trivialized by solving a holomorphic matrix differential equation.  The
underlying bundle is trivial in both cases; the distinction lies in the gauge
category.  This elementary example is the algebraic analogue of the regularity
issue already encountered for $\bar\partial$-connections in Section~\ref{sec:klv-module-picture}.

For the purposes of the sheaf-theoretic splitting problem, the minimal notion of
an \emph{admissible algebraic vector bundle} will therefore mean a finite locally
free module over the chosen base sheaf, locally trivialized on the projective Ore
cover, with algebraic restriction maps.  Equivariance, a lift to the principal
sheaf, and a compatible differential or Dolbeault connection are additional
admissibility conditions that define strictly finer categories.

\subsection{Transition matrices and the Hazewinkel--Martin factorization}
\label{subsec:transition-matrix-factorization}

Let
\[
  U_0=\operatorname{Spec}\Bbbk[z],
  \qquad
  U_1=\operatorname{Spec}\Bbbk[w],
  \qquad w=z^{-1},
\]
be the standard cover of $\mathbb P^1$.  Since vector bundles on both charts are
trivial, a rank-$r$ algebraic vector bundle is represented by a transition matrix
\begin{equation}
  G(z)\in\operatorname{GL}_r\bigl(\Bbbk[z,z^{-1}]\bigr).
  \label{eq:transition-matrix-loop}
\end{equation}
We use coefficient columns with $s_0=G(z)s_1$ on the overlap;
thus a scalar transition $G(z)=z^n$ agrees with
\eqref{eq:line-sheaf-transition} and represents $\mathcal O(n)$.
Changing the frame on $U_0$ and $U_1$ gives the equivalence relation
\begin{equation}
  G(z)\sim A(z)G(z)B(z^{-1}),
  \qquad
  A(z)\in\operatorname{GL}_r(\Bbbk[z]),
  \quad
  B(z^{-1})\in\operatorname{GL}_r(\Bbbk[z^{-1}]).
  \label{eq:transition-double-coset}
\end{equation}
Thus vector bundles are classified by a double-coset problem, not by conjugacy
inside a single group.

\begin{theorem}[Hazewinkel--Martin form of Birkhoff--Grothendieck]
\label{thm:hazewinkel-martin-matrix}
For every matrix
$G(z)\in\operatorname{GL}_r(\Bbbk[z,z^{-1}])$ there exist
\[
  A(z)\in\operatorname{GL}_r(\Bbbk[z]),
  \qquad
  B(z^{-1})\in\operatorname{GL}_r(\Bbbk[z^{-1}]),
\]
and uniquely determined integers
$n_1\geq\cdots\geq n_r$ such that
\begin{equation}
  A(z)G(z)B(z^{-1})
  =\operatorname{diag}(z^{n_1},\ldots,z^{n_r}).
  \label{eq:hazewinkel-martin-diagonal}
\end{equation}
Equivalently,
\[
  \mathcal E_G
  \simeq
  \mathcal O(n_1)\oplus\cdots\oplus\mathcal O(n_r).
\]
\end{theorem}

Hazewinkel and Martin give an elementary matrix proof, valid over an arbitrary
field, by reducing the transition matrix to a polynomial matrix and applying
normal-form arguments over the principal ideal domain $\Bbbk[z]$
\cite{HazewinkelMartin}.  In rank two the result reads
\begin{equation}
  A(z)G(z)B(z^{-1})
  =\begin{pmatrix}z^m&0\\0&z^n\end{pmatrix},
  \qquad m\geq n.
  \label{eq:rank-two-HM}
\end{equation}

\begin{corollary}[Splitting in the bare AFL sheaf category]
\label{cor:AFL-bare-splitting}
Every finite locally free
$\mathcal O_{q,\mathbb P^1}$-module is a direct sum of the standard line sheaves
$\mathscr O_q(n)$, uniquely up to permutation of the summands.
\end{corollary}

\begin{proof}
The local rings and the overlap ring are those in
\eqref{eq:qP1-chart-algebras}--\eqref{eq:qP1-overlap-algebra}, so
Theorem~\ref{thm:hazewinkel-martin-matrix} applies without alteration.
\end{proof}

Thus a rank-two ``quantum Hazewinkel--Martin theorem'' is already true in this
specific sheaf category, but for a simple reason: the transition-matrix problem
is the classical one.  A genuinely new theorem can arise only after one retains
additional quantum data, or after one replaces the overlap algebra by a genuinely
noncommutative algebra for which no classical Smith or Birkhoff factorization is
available.

\subsection{Classical frame groups versus quantum structure groups}
\label{subsec:GL-versus-GLq}

The preceding discussion also answers the question of whether a rank-$r$ vector
bundle should use $\operatorname{GL}_r$ or $\operatorname{GL}_q(r)$ as its
structure group.  These two choices express different notions.

A locally free right module of rank $r$ has local frames, and changes of frame are
by definition elements of
\begin{equation}
  \operatorname{GL}_r\bigl(\mathcal O_M(U)\bigr),
  \label{eq:ordinary-frame-group}
\end{equation}
the group of invertible matrices over the local coordinate algebra.  This is not
an additional choice: it is the automorphism group of the free module.  When the
local algebra is noncommutative one must fix left or right modules and preserve
the order in the cocycle and gauge formulas, but the matrix group
\eqref{eq:ordinary-frame-group} remains the natural frame group.

By contrast, choosing the coordinate Hopf algebra
$\mathcal O_q(\operatorname{GL}_r)$ means choosing a quantum principal bundle or
a reduction of the frame data to a particular quantum group.  In the strict
``point-valued'' interpretation, a $\operatorname{GL}_q(r)$-valued transition
function on an overlap algebra $R$ would be an algebra homomorphism
\begin{equation}
  \mathcal O_q(\operatorname{GL}_r)\longrightarrow R,
  \label{eq:quantum-group-valued-point}
\end{equation}
not an arbitrary invertible matrix in $M_r(R)$.  The entries of its matrix must
satisfy the Manin relations.  This is normally a much stronger condition.
More general gluing data for a quantum principal bundle are expressed through
convolution-invertible maps from the structure Hopf algebra to the overlap
algebra; such maps need not be algebra homomorphisms.  They still constitute
additional principal-bundle data and should not be confused with the full frame
group of an arbitrary locally free module.

The distinction is especially sharp for the quantum projective line.

\begin{proposition}[Quantum-group-valued loops are too restrictive]
\label{prop:GLq2-diagonal}
Assume that $q^2\neq1$ and let $R=\Bbbk[z,z^{-1}]$ be a commutative domain.
For every algebra homomorphism (that is, every strict $R$-valued point)
\[
  \rho:\mathcal O_q(\operatorname{GL}_2)\longrightarrow R,
\]
the matrix of coordinate functions has the form
\begin{equation}
  \begin{pmatrix}
    \rho(a)&\rho(b)\\
    \rho(c)&\rho(d)
  \end{pmatrix}
  =
  \begin{pmatrix}
    \alpha z^m&0\\
    0&\beta z^n
  \end{pmatrix},
  \qquad
  \alpha,\beta\in\Bbbk^\times,
  \quad m,n\in\mathbb Z.
  \label{eq:GLq2-diagonal-loop}
\end{equation}
\end{proposition}

\begin{proof}
In a standard presentation one has
\[
  ab=qba,\quad ac=qca,\quad bd=qdb,\quad cd=qdc,
  \quad bc=cb,
  \quad ad-da=(q-q^{-1})bc,
\]
and the quantum determinant is invertible.  Since $R$ is commutative and a
domain, the first four relations imply
$\rho(a)\rho(b)=\rho(a)\rho(c)=\rho(b)\rho(d)
=\rho(c)\rho(d)=0$.  If $\rho(a)=0$, invertibility of the quantum determinant
would force $\rho(b)\rho(c)$ to be a unit, contradicting
$(q-q^{-1})\rho(b)\rho(c)=0$.  Hence $\rho(a)\neq0$, so
$\rho(b)=\rho(c)=0$, and invertibility of the determinant forces
$\rho(a)$ and $\rho(d)$ to be units.  The units of $R$ are precisely the
monomials $\lambda z^k$.
\end{proof}

Thus, if one insists on strict $R$-point-valued
$\operatorname{GL}_q(2)$ transition functions over the commutative overlap of
the AFL projective line, diagonalization has been built into the definition.
Such transition functions restrict the allowed descent data and gauge
transformations. Nevertheless, every bare rank-two bundle has a diagonal
representative by Theorem~\ref{thm:hazewinkel-martin-matrix}; thus the
restriction does not exclude any bare splitting type. It fails to encode
arbitrary frames and extension data in the ordinary
$\operatorname{GL}_2(\Bbbk[z,z^{-1}])$ descent problem.  There is freedom to choose a quantum structure group, but that freedom
changes the problem: it classifies quantum principal reductions and their
associated bundles, not all locally free sheaves of a fixed rank.

\subsection{The Artin--Zhang and twisted homogeneous coordinate-ring picture}
\label{subsec:Artin-Zhang-quantum-line}

Let $A$ be a right noetherian connected graded algebra.  Write
$\operatorname{gr}A$ for the category of finitely generated graded right
$A$-modules and $\operatorname{tors}A$ for its Serre subcategory of
graded-torsion modules.  Artin and Zhang define the noncommutative projective
scheme associated with $A$ through the quotient category
\begin{equation}
  \operatorname{qgr}A
  :=\operatorname{gr}A/\operatorname{tors}A.
  \label{eq:qgr-definition}
\end{equation}
The analogue of the structure sheaf is $\pi A$, and its twists are
\begin{equation}
  \mathscr O_A(n):=\pi A(n),
  \label{eq:AZ-twists}
\end{equation}
where $\pi:\operatorname{gr}A\to\operatorname{qgr}A$ is the quotient functor.

For the quantum plane $A_q$ in
\eqref{eq:quantum-plane-homogeneous-ring}, this construction is completely
understood.  Let $S=\Bbbk[x_0,x_1]$ and let $\tau$ be the graded automorphism
\[
  \tau(x_0)=x_0,
  \qquad
  \tau(x_1)=q x_1.
\]
With the convention $f\star g=f\,\tau^{\deg f}(g)$ for homogeneous $f$,
we have $x_0\star x_1=q(x_1\star x_0)$, so $A_q$ is a Zhang twist of $S$.  Zhang's equivalence theorem gives
\begin{equation}
  \operatorname{Gr}A_q\simeq\operatorname{Gr}S,
  \label{eq:Zhang-graded-equivalence}
\end{equation}
and the equivalence preserves the torsion subcategories.  Consequently,
\begin{equation}
  \operatorname{qgr}A_q
  \simeq
  \operatorname{qgr}S
  \simeq
  \operatorname{coh}(\mathbb P^1).
  \label{eq:qgr-quantum-plane-classical}
\end{equation}
The same conclusion follows from the Artin--Van den Bergh theorem for the twisted
homogeneous coordinate ring
$B(\mathbb P^1,\mathcal O(1),\sigma_q)$.

\begin{theorem}[Birkhoff--Grothendieck in the Artin--Zhang quantum line]
\label{thm:AZ-BG}
Under the equivalence \eqref{eq:qgr-quantum-plane-classical}, the objects
$\mathscr O_{A_q}(n)$ correspond to the ordinary line bundles
$\mathcal O_{\mathbb P^1}(n)$.  Hence every vector-bundle object in
$\operatorname{qgr}A_q$ is isomorphic to a finite direct sum
\begin{equation}
  \mathscr O_{A_q}(n_1)\oplus\cdots\oplus\mathscr O_{A_q}(n_r),
  \qquad n_1\geq\cdots\geq n_r,
  \label{eq:AZ-splitting}
\end{equation}
and the integers are uniquely determined.
\end{theorem}

Thus, in the Artin--Zhang/Artin--Van den Bergh approach, the splitting question
for the two-generator quantum plane is not open.  The relevant result is a
theorem, ultimately because the quantum plane is a twist of the commutative
polynomial ring.  The paper of Artin, Tate and Van den Bergh on algebras
associated to automorphisms of elliptic curves provides part of the broader
historical framework, but the direct references for the present projective-line
statement are the theory of twisted homogeneous coordinate rings, Artin--Zhang
quotient categories, and Zhang twists
\cite{ArtinTateVdB,ArtinVdB,ArtinZhang,ZhangTwist,StaffordVdB}.

The categorical equivalence
\eqref{eq:qgr-quantum-plane-classical} does not erase the quantum algebra.  It
transports the grading shift to a twisted autoequivalence of
$\operatorname{coh}(\mathbb P^1)$ and retains the noncommutative multiplication
of the section ring.  Nor does it identify the Artin--Zhang category with the
KLvS category of projective modules carrying flat $\bar\partial$-connections.
That comparison is an additional theorem to be constructed.

\subsection{Which splitting problem remains open?}
\label{subsec:remaining-sheaf-problems}

The preceding results lead to a useful hierarchy of categories.

\begin{enumerate}[label=\textup{(\roman*)},leftmargin=2.3em]
\item
$\operatorname{Vect}(\mathbb P^1,\mathcal O_{q,\mathbb P^1})$, the category of
finite locally free modules over the bare AFL base sheaf.  This category is the
ordinary category of algebraic vector bundles on $\mathbb P^1$, and splitting is
Theorem~\ref{thm:hazewinkel-martin-matrix}.

\item
The associated or equivariant category generated by
$\mathcal F_q\square_H V$ for finite-dimensional $H=\mathcal O_q(P)$-comodules.
This category remembers the principal origin of a bundle.  Its forgetful functor
to the first category need not be full: an underlying sheaf morphism need not
preserve the retained equivariant or principal data. It can identify
nonisomorphic enriched objects and forget reductions of structure group.

\item
The Artin--Zhang category $\operatorname{qgr}A_q$.  For the quantum plane it is
equivalent to $\operatorname{coh}(\mathbb P^1)$, but it carries a distinguished
shift and a twisted homogeneous section ring.

\item
A differential sheaf category in which locally free sheaves are equipped with a
connection, or with a $(0,1)$-operator, compatible with the sheaf calculus of
\cite{AschieriFioresiLatiniWeber}.

\item
The KLvS category of finitely generated projective modules over the Podle\'s
sphere equipped with flat $\bar\partial$-connections.
\end{enumerate}

The vector-bundle parts of the first and third categories already satisfy
Birkhoff--Grothendieck. A further problem is to construct a comparison
between the fourth and fifth categories. Such a comparison must retain the
Dolbeault operator and the regularity conditions defining its morphisms.
An equivalence of underlying bare sheaf categories does not supply this
differential comparison.
Equivariance can help construct this comparison, but requiring a splitting
to preserve additional equivariance is a separate, stronger condition.

\begin{remark}[A differential comparison programme]
\label{prog:differential-quantum-GAGA}
A differential quantum GAGA statement requires a construction of its
categories and functor. For comparison with the polynomial KLvS framework
of this paper, take $\Bbbk=\mathbb C$, $0<q<1$, and left modules, with
polynomial holomorphic morphisms. The following are construction and
verification steps, rather than an asserted equivalence.
\begin{enumerate}[label=\textup{(\alph*)},leftmargin=2.3em]
\item Construct a passage from polynomial KLvS modules and their Dolbeault
calculus to sheaf data on the projective cover. Specify the local calculi,
the restriction maps, the descent condition, and the local connection
transformation law. In particular, verify that the construction retains
the noncommutative differential data over the commutative base sheaf.
\item Equip the proposed images of the standard line objects with explicit
sheaf $\bar\partial$-operators. With
$\mathcal O_q(n)=(\mathcal L_{-n},\overline\nabla_{(-n)})$ as in
\eqref{eq:intro-Oq}, the required underlying line sheaf is
$\mathscr O_q(n)$. Compare holomorphic Hom spaces and degree-one cohomology,
including their compositions, with the KLvS Dolbeault complexes.
\item Define the morphisms and admissible short exact sequences in the
resulting sheaf category, and prove that the construction is functorial
and exact on the chosen source category. Full faithfulness and a
description of the essential image are additional goals.
\item For tensor compatibility, specify the bimodule connections and the
twisted interchange maps; ordinary tensor products of central sheaves do
not supply this structure. If covariance is retained, specify its action
on objects and morphisms separately.
\end{enumerate}
An analytic extension would also require domains for the unbounded
Dolbeault operators and continuity and regularity results for the local
constructions. 
\end{remark}

\begin{remark}[The scope of a sheaf splitting question]
\label{rem:scope-sheaf-splitting}
Any comparison with the splitting results of this paper must distinguish
holomorphic isomorphisms from equivariant holomorphic isomorphisms. Already
at $q=1$, consider $\mathbb P^1\times\mathbb C^2$ with the trivial
holomorphic structure and the diagonal $\mathrm{SL}_2$ action, using the
standard representation on $\mathbb C^2$. Its ordinary splitting is
$\mathcal O\oplus\mathcal O$. However, it cannot split equivariantly into
the standard line objects: the underlying splitting forces both degrees
to be zero, each $\mathrm{SL}_2$-linearization of $\mathcal O$ is trivial,
and its global sections form the irreducible standard representation
rather than two trivial representations.

The same obstruction occurs in a  covariant quantum category containing
$\mathcal A_q\otimes V$, where $V$ is the two-dimensional type-1
$\mathcal U_q(\mathfrak{sl}_2)$-module, with diagonal covariance and
$\bar\partial$ acting on the $\mathcal A_q$ factor. Its underlying
holomorphic bundle is a sum of two standard degree-zero lines, but its
holomorphic sections transform as $V$. They cannot be identified
covariantly with the two trivial representations supplied by standard
covariant degree-zero lines. Thus a universal splitting assertion
retaining this equivariance would fail even for a holomorphically
trivial rank-two bundle.

The relevant splitting question is consequently holomorphic after
forgetting any such additional equivariance. In the polynomial category,
Theorem~\ref{thm:line-filtered-splitting} proves decomposition into
standard lines for the stated subcategory of bundles with a finite
filtration by standard line objects. A comparison functor, once
constructed and shown to preserve these objects and isomorphisms, would
transport that result to its image. The theorem does not supply such a
filtration for every holomorphic bundle.

Outside that subcategory, existence of rank-one holomorphic summands and
classification of the holomorphic structures on those summands remain
separate questions. In particular, one cannot infer that all rank-one
structures are standard. Proposition~\ref{prop:polynomial-line-nonrigidity}
exhibits nonstandard structures already in the polynomial gauge category.
Thus an unrestricted standard-line normal form is unavailable in that
category, even in rank one. No blanket enriched standard-line splitting
conjecture is asserted here.
\end{remark}

A comparison theorem must also match the relevant cohomologies. A proof
modeled on the classical splitting argument needs exact localization
and gluing, a satisfactory $H^1$-theory, compatible twists, and control
of the dimensions of twisted global sections. In the bare sheaf and
Artin--Zhang categories these tools are inherited from $\mathbb P^1$.
Their agreement with the KLvS Dolbeault cohomology of
Section~\ref{sec:klv-module-picture} is an additional comparison problem,
not a consequence of the commutativity of the affine base charts.

The algebraic quantum projective line therefore has a classical category
of bare locally free sheaves and a classical Birkhoff--Grothendieck
theorem, while its homogeneous coordinate ring and principal bundle
retain quantum data. A factorization theorem for arbitrary matrices
over the commutative overlap, or a replacement of $\operatorname{GL}_2$
by $\operatorname{GL}_q(2)$, cannot by itself resolve the KLvS problem.
The remaining work concerns the differential comparison and holomorphic
splitting in a precisely defined category, with line-rigidity and any
equivariant refinement treated as separate questions.

\section*{Acknowledgements}
\addcontentsline{toc}{section}{Acknowledgements}
We thank the Banff International Research Station (BIRS) for its hospitality
during a two-week stay in Banff in July 2026, when part of this work was
carried out. We also thank CIRM (Centre International de Rencontres
Mathématiques) for its hospitality during a two-week stay in August 2026,
when this work was completed. 
PA acknowledges partial support from INFN, CSN4, Iniziativa Specifica GSS.
This research received financial support from Universit\`a del Piemonte Orientale. PA is affiliated with INdAM-GNFM.
The work of MK was supported by a grant from the {\it Natural Sciences and Engineering Research Council of Canada} (NSERC). 
GL acknowledges partial support from INFN, CSN4, Iniziativa Specifica GAST.
GL is a member of INdAM-GNSAGA. GL acknowledges support from the PNRR MUR project PE0000023-NQSTI.

\paragraph{Use of AI tools.}
ChatGPT Pro was used to assist with formatting and editorial revisions, and to check details of
several examples in this paper. The authors take full responsibility for
the mathematical content and the final text.

\begingroup
\small
\raggedright
\setlength{\parindent}{0pt}

\begin{samepage}
\textbf{Paolo Aschieri}\\
Universit\`a del Piemonte Orientale,\\
Dipartimento di Scienze e Innovazione Tecnologica,\\
viale T. Michel 11, 15121 Alessandria, Italy,\\
and INFN Torino, via P. Giuria 1, 10125 Torino, Italy.\\
Email address: \href{mailto:paolo.aschieri@uniupo.it}{\nolinkurl{paolo.aschieri@uniupo.it}}
\end{samepage}

\medskip
\begin{samepage}
\textbf{Masoud Khalkhali}\\
Department of Mathematics, The University of Western Ontario,\\
London, ON, Canada.\\
Email address: \href{mailto:masoud@uwo.ca}{\nolinkurl{masoud@uwo.ca}}
\end{samepage}

\medskip
\begin{samepage}
\textbf{Giovanni Landi}\\
Matematica, Universit\`a di Trieste
and INFN, Sezione di Trieste, Trieste, Italy.\\
Email address: \href{mailto:landi@units.it}{\nolinkurl{landi@units.it}}
\end{samepage}

\endgroup

\end{document}